\documentclass[11pt,sort&compress]{elsarticle}

\usepackage{amsmath,amsfonts,bm}

\def\eqref#1{equation~\ref{#1}}

\def\1{\bm{1}}

\DeclareMathAlphabet{\mathsfit}{\encodingdefault}{\sfdefault}{m}{sl}
\SetMathAlphabet{\mathsfit}{bold}{\encodingdefault}{\sfdefault}{bx}{n}

\usepackage[utf8]{inputenc} 
\usepackage[T1]{fontenc}    
\usepackage{url}            
\usepackage{booktabs}       
\usepackage{amsfonts}       
\usepackage{amsmath,amssymb}
\usepackage{textcomp}
\usepackage{nicefrac}       
\usepackage{xfrac}

\usepackage[
    protrusion=true,
    activate={true,nocompatibility},
    final,
    tracking=true,
    kerning=true,
    spacing=true,
    factor=1100]{microtype}
\SetTracking{encoding={*}, shape=sc}{40}

\usepackage{lipsum}
\usepackage{graphicx,color}
\usepackage{algorithm}
\usepackage{algpseudocode}
\usepackage{multirow}
\usepackage{comment}
\usepackage{xcolor}
\usepackage[margin=1in]{geometry}

\definecolor{lightblue}{rgb}{0.63, 0.74, 0.78}
\definecolor{seagreen}{rgb}{0.18, 0.42, 0.41}
\definecolor{orange}{rgb}{0.85, 0.55, 0.13}
\definecolor{silver}{rgb}{0.69, 0.67, 0.66}
\definecolor{rust}{rgb}{0.72, 0.26, 0.06}

\colorlet{lightsilver}{silver!30!white}
\colorlet{darkorange}{orange!75!black}
\colorlet{darksilver}{silver!65!black}
\colorlet{darklightblue}{lightblue!70!black}
\colorlet{darkrust}{rust!85!black}
\colorlet{darkseagreen}{seagreen!85!black}

\usepackage{hyperref}
\hypersetup{
  colorlinks=true,
}

\usepackage[parfill]{parskip}

\usepackage{cleveref}

\usepackage{subcaption}

\usepackage{enumitem}

\usepackage[labelfont=bf,font=small]{caption}

\journal{Journal of Computational Physics}
\begin{document}

\hypersetup{
  linkcolor=darkrust,
  citecolor=seagreen,
  urlcolor=darkrust,
  pdfauthor=author,
}

\begin{frontmatter}
\title{{\bf\large Rational-WENO: A lightweight, physically-consistent three-point weighted essentially non-oscillatory scheme}}

\author{Shantanu~Shahane\fnref{google}}
\author{Sheide~Chammas\fnref{google}}
\author{Deniz~A.~Bezgin\fnref{tum}}
\author{Aaron~B.~Buhendwa\fnref{tum}}
\author{Steffen~J.~Schmidt\fnref{tum}}
\author{Nikolaus~A.~Adams\fnref{tum}}
\author{Spencer~H.~Bryngelson\fnref{gatech}}
\author{Yi-Fan~Chen\fnref{google}}
\author{Qing~Wang\fnref{google}}
\author{Fei~Sha\fnref{google}}
\author{Leonardo~Zepeda-Núñez\fnref{google}}
\fntext[google]{Google Research, 1600 Amphitheatre Pkwy, Mountain View, CA 94043 USA}
\fntext[tum]{Technical University of Munich, School of Engineering and Design, Chair of Aerodynamics and Fluid Mechanics, Boltzmannstr. 15, 85748
Garching bei München, Germany}
\fntext[gatech]{Georgia Institute of Technology, School of Computational Science \& Engineering, 756 W Peachtree St NW, Atlanta, GA 30332 USA}

\begin{abstract}
Conventional WENO3 methods are known to be highly dissipative at lower resolutions, introducing significant errors in the pre-asymptotic regime.
In this paper, we employ a rational neural network to accurately estimate the local smoothness of the solution, dynamically adapting the stencil weights based on local solution features. As rational neural networks can represent fast transitions between smooth and sharp regimes, this approach achieves a granular reconstruction with significantly reduced dissipation, improving the accuracy of the simulation.
The network is trained offline on a carefully chosen dataset of analytical functions, bypassing the need for differentiable solvers.
We also propose a robust model selection criterion based on estimates of the interpolation's convergence order on a set of test functions, which correlates better with the model performance in downstream tasks.
We demonstrate the effectiveness of our approach on several one-, two-, and three-dimensional fluid flow problems: our scheme generalizes across grid resolutions while handling smooth and discontinuous solutions. In most cases, our rational network-based scheme achieves higher accuracy than conventional WENO3 with the same stencil size, and in a few of them, it achieves accuracy comparable to WENO5, which uses a larger stencil.
\end{abstract}



\end{frontmatter}

\section{Introduction}\label{sec:intro}

For many initial value problems (IVP) in fluid mechanics, the solution of their underlying partial differential equation (PDE) may include spatial locations of high-gradients, or even discontinuities, stemming either from the initial condition or developed afterwards. This complexity of the problem requires the use of numerical methods that are tailored to resolve these areas to efficiently produce an overall accurate solution, presenting a significant challenge to the field.

Classical numerical methods for those PDEs often require human input and intuition to make choices when discretizing the PDEs, such as the choice of stencils.
For instance, in finite difference schemes, a smooth solution should favor wide-centered stencils, while a solution with discontinuities should instead favor shorter upwind stencils. These choices are not straightforward and usually involve a trade-off between high-order accuracy in the smooth region of the solution and sharply resolving discontinuities without introducing spurious oscillations, requiring careful decision-making.

This question spurred the development of high-resolution methods starting with~\cite{harten1983upstream}, and afterward, essentially non-oscillatory (ENO) schemes~\cite{harten1997uniformly}: ENO schemes measure the smoothness of the solution on several sub-stencils and then compute the flux based on the smoothest sub-stencil to avoid interpolating through the discontinuity. Further refinements of those ideas lead to weighted essentially non-oscillatory (WENO) schemes~\cite{jiang1996efficient,liu1994weighted,shu2020essentially}. WENO schemes introduce a continuous relaxation, effectively interpolating between different stencils by weighting them.
As the mesh is refined, the weights are designed to provide the optimal order of accuracy over the union of the sub-stencils, resulting in asymptotically high-order accuracy for smooth solutions.
However, in the pre-asymptotic regime, many methods are prone to provide sub-optimal approximations, particularly when the underlying solution is smooth but with high gradients. Such regions might be incorrectly assessed as discontinuous, resulting in overly dissipative behavior.

While methods using higher-order polynomial discretizations~\cite{CockburnShu1998,CockburnShu:01,CockburnKarniadakisShu:00} within cells (i.e., subdivisions of the computational domain) or adaptive meshes~\cite{BALSARA2020109062} can mitigate this issue, WENO schemes remain prevalent in engineering and science due to their simplicity, which can be optimized in modern computing hardware~\citep{radhakrishnan24,bryngelson19_CPC,elwasif2023application}. Unfortunately, computational constraints and the complexity of modern simulations often require using WENO schemes at coarse resolutions, far from their asymptotic ideal.


Reducing the errors of these methods in the \textit{pre-asymptotic} regime has thus become a priority.
Several analysis-based methodologies with enhanced smoothness indicators have been proposed~\cite{hu2010adaptive,borges2008improved,henrick2005mapped,fu2016family}; however, they remain sub-optimal in the pre-asymptotic regime, as they tend to transition slowly between stencils, even as they regain optimality in the asymptotic limit.
Therefore, a new crop of methods has recently emerged that leverages machine learning (ML) to develop smoothness indicators for optimizing stencil selection~\cite{bezgin2022weno3,stevens2020enhancement}.

In this work, we focus on the third-order classical WENO3, whose narrow and localized stencil renders the scheme computationally efficient but hampers its capability to accurately assess the smoothness of the function, resulting in more dissipative, and therefore less accurate, simulation.
We design a \textit{data-driven smoothness indicators} for WENO3 using rational neural networks to address this challenge. While other neural networks have been used before for similar purpose~\cite{bezgin2022weno3,stevens2020enhancement}, rational neural networks provide additional modeling advantages as they can efficiently approximate discontinuous functions~\cite{boulle2020rational}.

Additionally, we propose a novel criterion for model selection. As the stochastic and non-convex nature of training often leads to trained models with vastly different behaviors, a selection criterion that ensures performance on downstream tasks is key. We demonstrate that using reconstruction misfit in a test set is inadequate, as it does not fully correlate with model performance in numerical simulations. Instead, we propose an empirically robust criterion: compute estimates of the order interpolation for a carefully chosen function and select the model whose estimate is closest to the theoretical optimum. We show that such a criterion ensures the model's performance in downstream simulations.

Finally, we showcase the behavior of our methodology in several numerical examples in the coarse resolution regime, including simple one-dimensional PDEs, two-dimensional turbulent flows, a bubble problem, and large scale cloud simulations exhibiting three-dimensional  turbulence. In all of the examples, the new methodology is able to represent phenomenological properties more accurately than other classical and competing ML-based methods.

\subsection{Related Work}

We provide a succinct review of the expansive (and still growing) literature in solving PDEs. We can decompose it in three main axes, fast and high-order methods, purely learned methods, and hybrid ones in which ML plays a role of augmenting traditional numerical methods.

\textbf{Fast PDE solvers} typically aim at leveraging the analytical properties of the underlying PDE to obtain low-complexity and highly parallel algorithms. Despite impressive progress~\citep{Martinsson:fast_pde_solvers_2019}, they are still limited by the need to mesh the ambient space, which plays an important role in the accuracy of the solution, as well as the stringent conditions for stable time-stepping. Other techniques are based on adaptive high-order polynomials, such as Discontinuous Galerkin approximations~\cite{CockburnShu1998,CockburnShu:01,CockburnKarniadakisShu:00} and their Hybridizable variant~\cite{nguyen2011implicit}, and even spectral methods with dissipation~\cite{fontana2020fourier} to avoid the Gibbs phenomenon.

\textbf{Machine-Learning models} use neural networks to represent either the solution or the solution operator. This category can be further divided in two sub-groups.

\textit{Neural Ans\"atze} methods aim to leverage the approximation properties of neural networks~\citep{hornik90}, by replacing the usual linear combination of handcrafted basis functions by a more general neural network ansatz. The physics prior is incorporated explicitly by enforcing the underlying PDE in strong~\citep{raissi_physics-informed_2019,eivazi_physics-informed_2021}, weak~\citep{E:2018_deep_ritz,Gao:2022_Galerkin_PINN}, or min-max form~\citep{ZANG:2021_WAN}. Besides a few exceptions, e.g.,~\citep{hyperpinn,Bruna_Peherstorfer:2022}, these formulations often require solving highly non-convex optimization problems at inference time. 

\textit{Purely Learned Surrogates} fully replace numerical schemes with surrogate models learned from data. A number of different architectures have been explored, including multi-scale convolutional neural networks~\citep{ronneberger_u-net_2015,wang_towards_2020}, graph neural networks~\citep{sanchez-gonzalez_learning_2020}, Encoder-Process-Decoder architectures~\citep{stachenfeld_learned_2022}, and neural ODEs~\citep{Ayed_Gallinari:2019}. Similarly, operator learning seeks to approximate the inverse of the underlying differential operator directly by mimicking the analytical properties of its class, such as pseudo-differential~\citep{Hormander:PseudoDiff}
or Fourier integral operators~\citep{Hormander:FIO},
but without explicit PDE-informed components. These methods often leverage the Fourier transform~\citep{li_fourier_2021,tran_factorized_2021}, the off-diagonal low-rank structure of the associated Green's function~\citep{FanYing:mnnh2019,graph_fmm:2020}, or approximation-theoretic structures~\citep{deeponet:2021}. 

\textbf{Hybrid Physics-ML} methods hybridize classical numerical methods with contemporary data-driven deep learning techniques~\citep{mishra_machine_2019,bar-sinai_learning_2019,kochkov_machine_2021,Sirignano_DPM:2020,bruno_fc-based_2021,Frezat2022-fs,Dresdner:2020ml_spectral}.
These approaches \emph{learn} corrections to numerical schemes from high-resolution simulation data, resulting in fast, low-resolution methods with an accuracy comparable to the ones obtained from more expensive simulations. Some of the hybrid physics-ML models are trained online using differentiable solvers \cite{bezgin2023jax}.

Our proposed method can be classified in this last category, although, instead of learning the correction online, we learn it offline using analytical data. This has several advantages: 
i) there is no need for a differentiable solver, ii) the output of the network can be easily tested against an analytical and exact ground truth, whereas using a solver, the output of the model is tested against a numerical solution which contains biases and errors from the solver, iii) we output the interpolation coefficient, which are themselves constrained to give at least a second order interpolation. Therefore, the resulting model is robust, as it will respect the conservation law, and it will provide at least a low-order answer in the worst case. 
 
\subsection{Organization}

This manuscript is organized as follows in \cref{sec:weno} we briefly introduce the third-order scheme WENO scheme used as basis of our methodology. The rational functions and architecture of the neural network are introduced in \cref{sec:neural_network_architecture}.
In \cref{sec:training}, we provide details on the training pipeline and model selection. The results showcasing the behavior of the network when used for solving different equations are provided in \cref{sec:results}. Further, we summarize our key findings in \cref{sec:conclusion}.

\section{Review of WENO schemes} \label{sec:weno}

Although we focus on WENO3 in this paper, our approach could be used to enhance many other related methods. For the sake of completeness, we briefly describe the algorithmic pipeline of WENO methodology.

For simplicity, consider the hyperbolic scalar transport equation:
\begin{equation}
    \frac{\partial u}{\partial t} + \frac{\partial}{\partial x} f(u) = 0.
    \label{eq:hyperbolic_pde}
\end{equation}

\begin{figure}[t]
     \centering
     \begin{subfigure}[b]{0.4\textwidth}
         \centering
         \includegraphics[width=\textwidth]{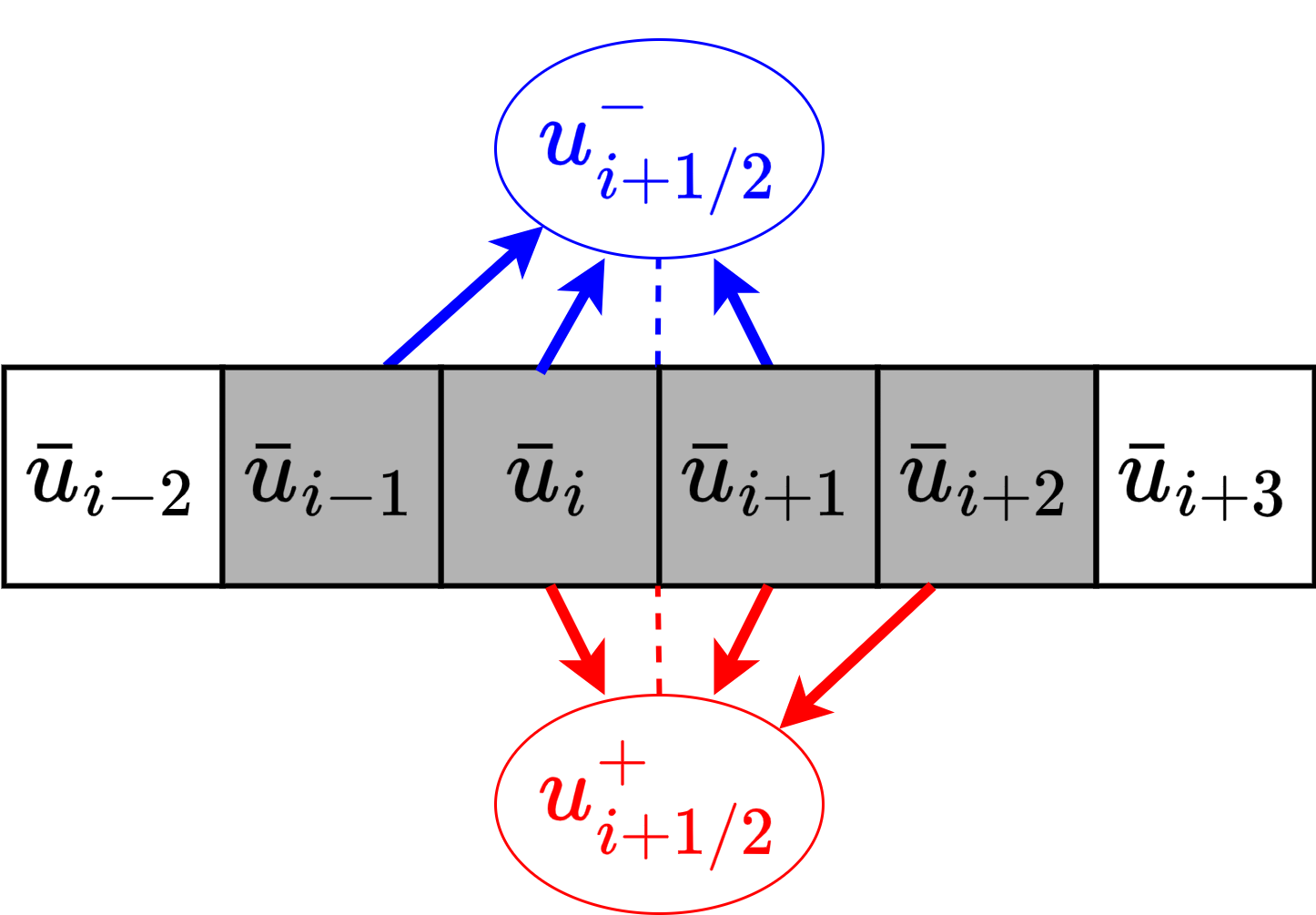}
         \caption{WENO3}
     \end{subfigure}
     \hfill
     \begin{subfigure}[b]{0.4\textwidth}
         \centering
         \includegraphics[width=\textwidth]{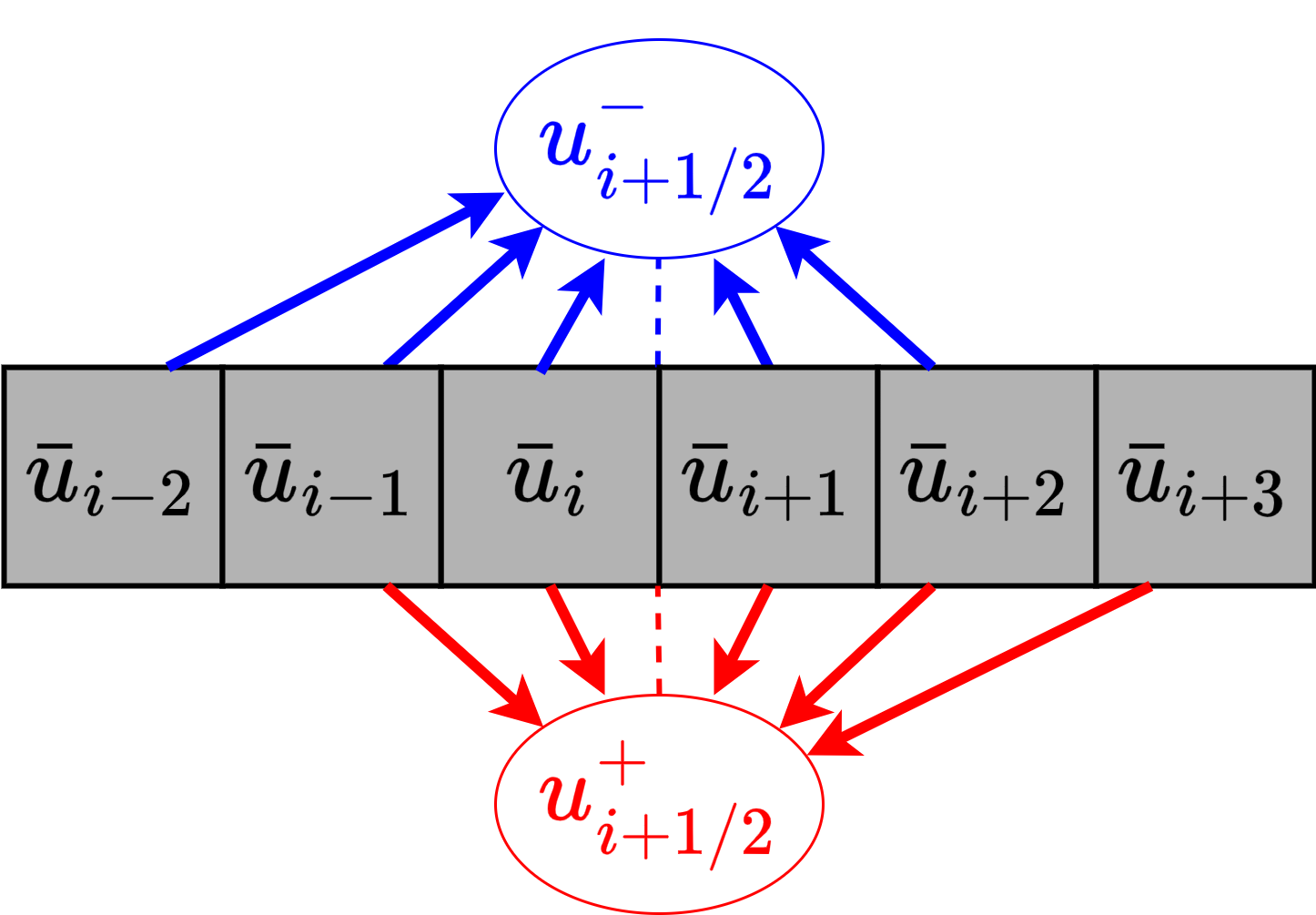}
         \caption{WENO5}
     \end{subfigure}
    \caption{Location of cell centers, faces and stencils. Here $\bar{u}_i$ is the cell average of the solution at the $i$-th cell, $u^{-}_{i+1/2}$ is the left-value of the solution at the interface between the $i$-th and $i+1$-th cells. Similarly, $u^{+}_{i+1/2}$ is the right-value of the solution at the same interface. WENO schemes seek to interpolate the value of the solution at the interfaces using neighboring cell averages.
    WENO3 uses three adjacent cells to compute the interpolation and WENO5 uses five cells.}
    \label{fig:stencil_schematic}
\end{figure}
Consider a uniform spatial discretization in cells, as shown in \Cref{fig:stencil_schematic}, and let $\bar{u}_i$ be the average value of the field $u(x,t)$ over a finite volume cell $i$ centered at $x_i$, i.e., 
\begin{equation} \label{eq:cell_mean}
    \bar{u}_i = \frac{1}{\Delta x} \int_{x_{i-\sfrac{1}{2}}}^{x_{i+\sfrac{1}{2}}} u(x,t) dx, 
\end{equation}
which is then used to discretize \cref{eq:hyperbolic_pde} in space yielding
\begin{equation} \label{eq:discretized_pde}
    \frac{\partial \bar{u}_i}{\partial t} = - \frac{f(u_{i+\sfrac{1}{2}}) - f(u_{i-\sfrac{1}{2}})}{\Delta x}.
\end{equation}
The physical flux $f(u)$ is approximated by its numerical counterpart $\hat{f}$ (usually called the numerical flux) which is estimated from the interface values of $u$ taking the form
\begin{equation} \label{eq:num_flux}
    \hat{f}_{i\pm\sfrac{1}{2}} = \hat{f}(u^-_{i\pm\sfrac{1}{2}}, u^+_{i\pm\sfrac{1}{2}}).
\end{equation}
where, $u^-$ and $u^+$ denote the interpolated values at faces based on the left and right biased stencils respectively (\cref{fig:stencil_schematic}).
Replacing this expression in \cref{eq:discretized_pde} we obtain the system
\begin{equation}
    \frac{\partial \bar{u}_i}{\partial t} = - \frac{\hat{f}_{i+\sfrac{1}{2}} - \hat{f}_{i-\sfrac{1}{2}}}{\Delta x},
\end{equation}
which is discretized in time and solved using time integrators~\cite{Dormand-Prince:1980}.

Even though the exact form of the numerical flux is method dependent (e.g., Roe, Godunov, etc.~\cite{leveque2002finite}), they usually follow the form in \cref{eq:num_flux}, i.e., they are function of the field $u$ at the interface of the cell ($u^{\pm}_{i+\sfrac{1}{2}}$). Thus, such quantity needs to be estimated from cell averages of the field on adjacent cells.
\Cref{fig:stencil_schematic} illustrates a schematic of uniform discretization along the $x$ coordinate, highlighting the stencils employed for $u^{\pm}_{i+\sfrac{1}{2}}$ in both WENO3 and WENO5 schemes.
In what follows we describe the computation of the negative side only $u^-_{i+\sfrac{1}{2}}$ based on the cell values: $[\bar{u}_{i-1}, \bar{u}_i, \bar{u}_{i+1}]$. $u^+_{i+\sfrac{1}{2}}$ is defined similarly.

WENO methods assume that the interface value, $u^-_{i+\sfrac{1}{2}}$, is a convex combination of sub-stencil values:
\begin{equation} \label{eq:convex_combination}
    u^-_{i+\sfrac{1}{2}} = \sum_{k=0}^{r-1} \omega_k u^{(k)}_{i+\sfrac{1}{2}}, \text{ where } \omega_k \geq 0 \text{ and } \sum_{k=0}^{r-1} \omega_k = 1.
\end{equation}
For third-order accuracy, i.e., $r=2$, the interpolants for the sub-stencil are given by
\begin{equation} \label{eq:interpolants}
    u^{(0)}_{i+\sfrac{1}{2}} = \frac{-\bar{u}_{i-1} + 3 \bar{u}_i}{2} \qquad \text{and} \qquad u^{(1)}_{i+\sfrac{1}{2}} = \frac{\bar{u}_i + \bar{u}_{i+1}}{2}.
\end{equation}
The classical WENO3-JS scheme~\cite{jiang1996efficient} calculates the weights $\omega_k$ as:
\begin{equation} \label{eq:weno_features}
    \omega_k = \frac{\alpha_k}{\sum_{k=0}^{r-1} \alpha_k}, \text{ where } \alpha_k = \frac{d_k}{(\beta_k + \epsilon)^2}.
\end{equation}
where, $\epsilon$ is set to a small positive value to avoid division by zero. In this case, the \textit{smoothness indicators} $\beta_k$ are defined as following
\begin{equation}
\beta_0 = (\bar{u}_i - \bar{u}_{i-1})^2 \qquad \text{and} \qquad \beta_1 = (\bar{u}_i - \bar{u}_{i+1})^2,
\end{equation}
The ideal weights ($d_k$) for third-order accurate upwind scheme are given by
\begin{equation} \label{eq:optimal_weights}
    d_0=\frac{1}{3} \qquad \text{and} \qquad d_1=\frac{2}{3}.
\end{equation}
In the upcoming sections, we also make comparisons to the WENO3-Z method, employing an alternative approach for calculating $\omega_k$~\cite{acker2016improved}.

\section{Rational Neural Network Architecture} \label{sec:neural_network_architecture}
In this section, we provide details on augmenting the WENO scheme with rational neural networks, a class of multi-layer perceptron with rational functions as activation function. Here we leverage the capability of rational function to efficiently approximate discontinuous functions such as indicator function of disjoint sets\footnote{This approximation problem is usually called the fourth Zolotarev problem~\cite{zolotarev1877application}.} to create better smoothness indicators for WENO schemes.

Following the notation introduced in \cref{sec:weno}, a traditional WENO3 scheme seeks to compute the weights $ \{\omega_0, \omega_1\}$ for the interpolants in \cref{eq:interpolants} taking the form
\begin{equation}
    \omega_k = \omega_k(\bar{u}_{i+1}, \bar{u}_{i}, \bar{u}_{i-1}),\qquad \text{for } k = 0, 1.
\end{equation}
We seek to replace this function by a neural network with the same input,
\begin{equation}
    \omega_k^{\mathrm{NN}} = \Omega(\bar{u}_{i+1}, \bar{u}_{i}, \bar{u}_{i-1}; \bm{\Theta}), 
\end{equation}
where $\bm{\Theta}$ encodes all the learnable parameters of the network. 
In what follows, we provide details on how such network is built.


Let $p$ and $q$ be polynomials of the form $p(x) = \sum_{i=0}^n c_i x^i$ and $q(x) = \sum_{i=0}^m \tilde{c}_i x^i$, where $n$ and $m$ are their corresponding degrees. We define a $(n,m)$-order rational function as 
\begin{equation}
    \mathcal{R}_{\theta}(x) := \frac{p(x)}{q(x)}.
\end{equation}
where $\theta$ is the set of the coefficients of $p$ and $q$: $ \theta := \left [ \{c_i\}_{i=0}^{n}, \{\tilde{c}_i\}_{i=0}^{m} \right ]$. 

We use these functions as activation functions, and render the coefficients to be learnable. 
We follow~\cite{boulle2020rational}, and consider $(3,2)$-order rational function: a third-order polynomial in the numerator, and second order in the denominator. 

\paragraph{Featurization}
We also leverage rational networks to create the input features. We take inspiration from the traditional WENO scheme whose features $\alpha_k$ are given by simple rational functions (as shown in \cref{eq:weno_features}) which are then soft-maxed (with a polynomial instead of an exponential). Analogously, we build features that are also rational functions with parameters that are learnt from data, and whose outputs are also softmaxed.

Following \cref{eq:weno_features} and~\cite{bezgin2022weno3} we use a first layer of features, which are given by the finite differences features from the input local averages $[\bar{u}_{i+1}, \bar{u}_{i}, \bar{u}_{i-1}]$
\begin{equation} \label{eq:delta_features}
    \Delta_1 = |\bar{u}_i - \bar{u}_{i-1}|, \,\,  \Delta_2 = |\bar{u}_{i+1} - \bar{u}_{i}|, \,\,  \Delta_2 = |\bar{u}_{i+1} - \bar{u}_{i-1}|, \,\, \Delta_4 = |\bar{u}_{i+1} - 2\bar{u}_{i} + \bar{u}_{i-1}|.
\end{equation}
Each of these features, which are designed to be Galilean invariant, are then fed to a rational network, whose output is then normalized using a soft-max. The final features, which we denote by $a^{0}$, are given by: 
\begin{equation} \label{eq:rational_feats}
\alpha([\bar{u}_{i+1}, \bar{u}_{i}, \bar{u}_{i-1}]) := \left (
\begin{array}{c}
    \mathcal{R}_{\theta_1} (\Delta_1) \\
    \mathcal{R}_{\theta_2} (\Delta_2) \\
    \mathcal{R}_{\theta_3} (\Delta_3) \\
    \mathcal{R}_{\theta_4} (\Delta_4)
\end{array}
\right ),  \,\, a^{0} := \frac{\alpha}{ \| \alpha \|},
\end{equation}
where $a^0 \in \mathbb{R}^4$, i.e., the first hidden layer has four neurons, corresponding to the rational features.
Each feature has its own set of weights $\theta_i$, for $i = 1,...,4$. We point out that the authors in \citet{bezgin2022weno3} use a different featurization preserving Galilean invariance: where they normalize the features in \cref{eq:delta_features} then feed them directly to an multi-layer perceptron with ReLU activation units.

\paragraph{Rational Layers} 
As usual, we construct the rational network by alternating the application of the rational function and the linear transformation. The recursive expression for each layer (with the first layer given by $a^0$ defined in \cref{eq:rational_feats}) is as follows 
\begin{equation} \label{eq:rational_layers}
    a^{\ell}_{i} = \mathcal{R}_{\theta^{\ell}} \left ( \sum_{j=1}^{N^{\ell -1}}  W^{\ell-1}_{i, j} a_{j}^{\ell-1} + b_i^{\ell-1} \right).
\end{equation}
where, $\ell$ is the index for the layer, $N^{\ell -1}$ is the number of neurons in the $\ell$-th layer,  $W^{\ell} \in \mathbb{R}^{N^{\ell} \times N^{\ell+1}}$ is the weight matrix linking layers $\ell$ and $\ell+1$, and $b^{\ell} \in \mathbb{R}^{N^{\ell+1}}$ is the bias vector. Each rational function has its own set of weights $\theta^{\ell}$, which are shared among neurons in the same layer. 

At the last layer of the network we fix the dimension of the output to two, and we apply a softmax function, i.e., 
\begin{equation} \label{eq:last_layer}
    \omega^{\mathrm{NN}} = \text{softmax}\left (  W^{L-1} \cdot a^{L-1} + b^{L-1} \right),
\end{equation}
where $W^{L-1} \in \mathbb{R}^{2 \times N^{L-1}}$ and $b^{L-1} \in \mathbb{R}^2$. 

\paragraph{ENO Layer} 

The ENO layer introduced in \citet{bezgin2022weno3} is used during inference time to guide the networks to use the asymptotically correct weights, as suitable WENO weights should be close to the ideal weights in smooth flow regions while the stencil with discontinuity should be assigned effectively zero weight. However, due to saturation effect of the softmax function, the output of the network lies in $(0, 1)$, thus the extrema points are never reached. To solve this issue, the ENO layer introduces a hard-thresholding function so that the network recovers the ENO property. I.e., the outputs of the networks are given by 
\begin{equation} \label{eq:eno_layer}
    \widehat{\omega}^{\text{NN}}_k = \frac{\phi(\omega^{\text{NN}}_k)} {\sum_j \phi(\omega^{\text{NN}}_j)},
\end{equation}
where $\phi$ is the hard-thresholding function with threshold $c_{\text{eno}}>0$, given by
\begin{equation}
    \phi(x) = \left \{ \begin{array}{lc} 0, & \text{if } x < c_{\text{eno}}, \\ 
                                x, & \text{if } x \geq c_{\text{eno}}.
                       \end{array}
              \right.
\end{equation}
The normalization in \cref{eq:eno_layer} ensures that resulting interpolation resides in the convex hull of the inputs. Similar to~\cite{bezgin2022weno3}, we set $c_{\text{eno}}=2e-4$.
During training, the ENO layer is deactivated to avoid issues in the backpropagation of gradients. 

\begin{figure}
    \centering
    \includegraphics[width=0.7\textwidth, clip, trim={0mm 80mm 0mm 0mm}]{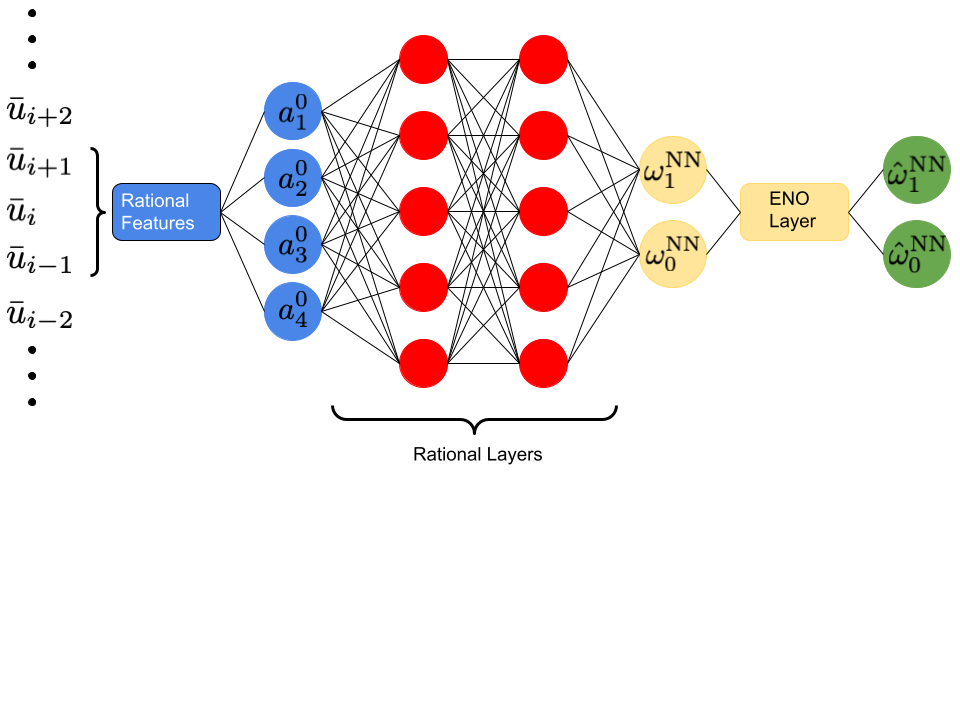}
    \caption{Sketch of the architecture, the network takes cell averages $[\bar{u}_{i+1}, \bar{u}_i, \bar{u}_{i-1}]$ as input and applies rational featurization (\cref{eq:rational_feats}). The resulting features are passed to a rational MLP (\cref{eq:rational_layers}). During inference, an ENO layer (\cref{eq:eno_layer}) generates weights used in a convex combination (\cref{eq:convex_combination}) to approximate boundary values $u_{i\pm \sfrac{1}{2}}$. These values are inputs to the numerical flux calculation (\cref{eq:num_flux}).}
    \label{fig:architecture_sketch}
\end{figure}

Finally, we use $\widehat{\omega}^{\text{NN}}_k$ with the interpolants to build the interface values following \cref{eq:convex_combination} to obtain
\begin{equation} \label{eq:weno_nn_reconstruction}
    u^{\text{NN}}_{i+\sfrac{1}{2}} = \sum_{k=0}^{1} \widehat{\omega}^{\text{NN}}_k u^{(k)}_{i+\sfrac{1}{2}}.
\end{equation}

\section{Training and regularization} \label{sec:training}

\subsection{Loss}

We follow the training loss used in \citet{bezgin2022weno3} where, the total loss ($\mathcal{L}$) is decomposed in three components: reconstruction loss ($\mathcal{L}_r$), deviation loss ($\mathcal{L}_d$) and the $l_2$-regularization loss, which are then weighted resulting in
\begin{equation}
     \mathcal{L} = \mathcal{L}_r + \beta_d \mathcal{L}_d + \beta_W ||\bm{\Theta}||_2^2,
     \label{eq:training_loss}
\end{equation}
where
\begin{align} \label{eq:losses}
    \mathcal{L}_r &= \frac{1}{N_b} \sum_{[s]=1}^{N_b} \left(\gamma^{[s]} \right)^\alpha \left(u^{\text{NN},[s]}_{i+\sfrac{1}{2}} - u^{[s]}_{i+\sfrac{1}{2}}\right)^2,\\
    \mathcal{L}_d &= \frac{1}{N_b} \sum_{[s]=1}^{N_b} \left(1- \left(\gamma^{[s]} \right)^\alpha\right) \sum_{k=0}^{r-1} \left(\omega^{\text{NN},[s]}_k - d_k \right)^2.
\end{align}
The $l_2$-regularization in \cref{eq:training_loss} is applied over all the weights of the rational and fully connected layers ($\bm{\Theta}$). 
Here, the expressions are averages taken over samples indexed by $[s]$ where $N_b$ is the total number of training samples. In addition,  $\{d_k\}_{k = 0}^1$ are the optimal weights defined in \cref{eq:optimal_weights} and $\gamma^{[s]} \in [0,1]$ seeks to quantify the local smoothness of the function from it samples as
\begin{equation}
    \gamma^{[s]} = \frac{|\bar{u}_{i-1}^{[s]} - 2 \bar{u}^{[s]}_i + \bar{u}_{i+1}^{[s]}|}{ |\bar{u}^{[s]}_i - \bar{u}^{[s]}_{i-1}| + |\bar{u}^{[s]}_i - \bar{u}^{[s]}_{i+1}| + \epsilon_\gamma},
\end{equation}
where $\epsilon_\gamma=10^{-15}$ is a small positive number to avoid division by zero.

For smooth functions, $\gamma$ is close to zero (with the exception of critical points) and hence, the deviation loss ($\mathcal{L}_d$) is dominant. Thus, the neural network is expected to predict WENO weights ($\omega^\text{NN}_k$) closer to the ideal weights $d_k$ giving optimal convergence. On the other hand, for stencils with discontinuities, $\gamma$ is close to one, thus, the reconstruction loss ($\mathcal{L}_r$) is dominant and the loss nudges the network to reproduce the interpolated value ($u^{[s]}_{i+\sfrac{1}{2}}$). Hence, the $\omega^\text{NN}_k$ are away from the $d_k$ and this avoids any oscillatory behavior near the discontinuity.
$\alpha$, $\beta_d$ and $\beta_W$ are non-negative scalar hyper-parameters which we obtain using the hyper-parameter sweep. Higher values of $\alpha$ and $\beta_d$ promote the discovery of ideal upwind weights $d_k$.

\paragraph{Initialization}
We initializate the parameters of the rational function with pre-computed weights $ \theta$ such that $R_{\theta}$ approximates a ReLU function, as described in \citet{boulle2020rational}.
The rest of the network is initialized using standard the LeCun normal initialization~\cite{klambauer2017self}.

\subsection{Training data}

Each data pair is composed of the cell averages of a function, $(\overline{u}_{i-1},\overline{u}_{i}, \overline{u}_{i+1})$, and the corresponding value of the function at the interface, $u_{i + \sfrac{1}{2}}$. The training data is generated from a collection of pre-chosen discretized analytical functions. For each discretization, we \textit{exactly} compute the cell averages $\overline{u}_i$ following \cref{eq:cell_mean} and the values at the interface $u_{i + \sfrac{1}{2}}$ using the analytical expression of the functions and their integrals. Furthermore, the data is post-processed to enforce that $u_{i + \sfrac{1}{2}}$ is always a convex combination of $(\overline{u}_{i-1},\overline{u}_{i}, \overline{u}_{i+1})$, as the network itself is constrained accordingly to ensure that the overall scheme is total variation diminishing.

\begin{table}[t]
\centering
\small
\begin{tabular}{c c}
Function $f(x)$     & Random Parameters \\ \midrule
$\sum_{k=0}^n c_k x^k$    & $c_k \in \mathcal{U}(-1,1)$ $\forall k$, $n=3$      \\ 
$u_l$ if $x<0.5$ and $u_r$ otherwise & $u_l, u_r \in \mathcal{U}(-1,1) $   \\ 
$(-1)^a x + \delta (x>0.5)$ & $a \in \mathcal{B}(0.5)$, $\delta \in \mathcal{U}(0.5,1) $   \\ 
$\sin(k \pi x)$ & $k \in \mathcal{U}(2,20) $   \\ 
$\tanh(k x)$ & $k \in \mathcal{U}(5,30)$
\end{tabular}
\caption{Analytical functions used for generating training dataset.}
\label{tab:train_functions}
\end{table}

The family of functions in \cref{tab:train_functions} and their integrals are sampled over one of the following one-dimensional domains
\begin{itemize}
    \item $x \in [-1,1]$, for polynomials and hyperbolic tangent functions, and
    \item $x \in [0,1]$ for other functions.
\end{itemize}

The number of discrete points ($n_x$) along the $x$ coordinate ranges from 16 to 1024 following geometric progression with ratio $2$ (i.e., $2^i$). Each grid size $n_x$ gives $n_x$ data pairs for training. To ensure a constant number of data pairs (16384) across all grid sizes, the number of random samples (\cref{tab:train_functions}) is adjusted inversely proportional to $n_x$. For example, when $n_x=16$, we have $16384/16=1024$ random samples, while for $n_x=32$, we have 512 samples and so on. This adjustment maintains that each value of $n_x$ has equal representation in the training set in terms of data pairs.

We point out that one could use simulation data for the training set to bias the training towards physically relevant data. However, such data usually has discretization induced errors, whose statistics are not consistent across resolutions. This hampers our model selection criterion (explained below) which hinges on estimates of the interpolation order of convergence. Using analytical functions bypasses this issue at the cost leaving untapped the underlying statistics of physically relevant solutions. We speculate that using high-accuracy physically relevant data should further improve the accuracy of the methodology.

\subsection{Model selection}\label{sec:model_selection}
Machine learning  requires  selecting best performing models when there are a large number of tuning knobs. Examples of such knobs are training hyper-parameters, such as number of epoch, learning rate schedule, among others; and network parameters, such as number of layers, and number of neurons at each layer. 

Model selection is an key step for obtaining high-performing model. Conventionally, this is done by selecting the best performing model by measuring its  loss on validation data set where the loss is usually the same as used for training.

Here we propose a better criterion based on estimates on the order of convergence on a validation data set, which is generated using the following two functions:
\begin{align}
    g(x) &= \sin^3(\pi x) \,\, x \in [-1,1], \\
    h(x) &= 
    \begin{cases}
    \sin(2 \pi x), & \text{for } x \in [0, 0.5).\\
    1 + \sin(2 \pi x), & \text{for } x \in [0.5, 1].
  \end{cases}
\end{align}
Each of the above functions is discretized in a grid with number of discretization points, $n_x$, ranging from 16 to 1024 increasing by a multiple of 2. During training, we assessed model performance for each grid resolution by calculating interpolation errors and individual loss terms (\cref{eq:training_loss}). We defined the order of accuracy as the slope of the log-log plot relating grid spacing $\Delta x$ to interpolation error.

Three distinct criteria were used for model selection: minimizing reconstruction loss, minimizing deviation loss, and achieving an empirical convergence order closest to the theoretical third-order. \Cref{tab:net_hyper_params} presents the top-performing models identified for each of these criteria. The performance of these WENO-NN models is analyzed in detail in \cref{sec:results}.

\begin{figure}[t]
     \centering
     \begin{subfigure}[b]{0.49\textwidth}
         \centering
         \includegraphics[width=\textwidth]{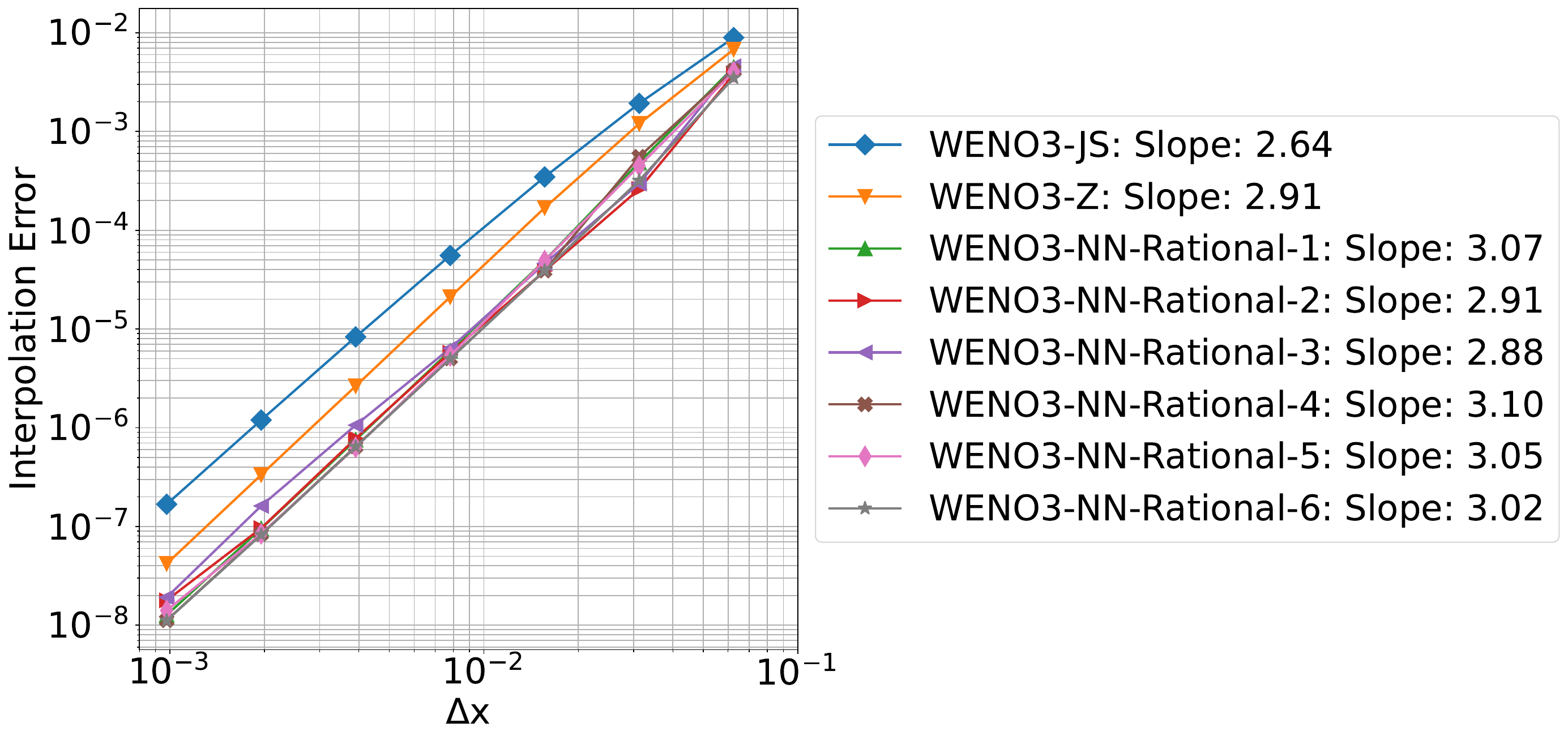}
         \caption{Convergence on the $\sin^3(\cdot)$ function}
     \end{subfigure}
     \begin{subfigure}[b]{0.49\textwidth}
         \centering
         \includegraphics[width=\textwidth]{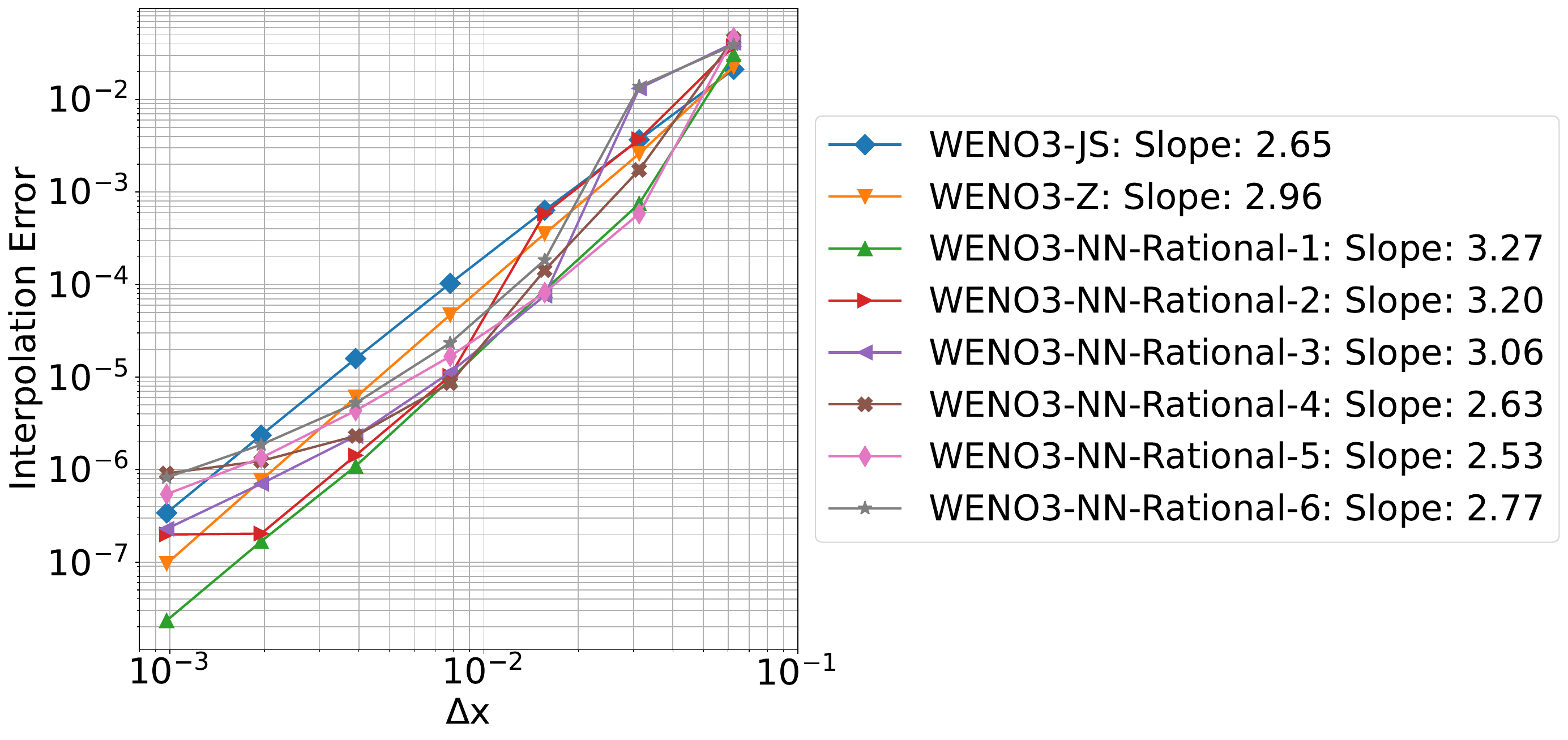}
         \caption{Convergence on the sine-step function}
     \end{subfigure}
    \caption{
        Convergence of rational neural models on the evaluation functions where the order is estimated as the coefficient of linear regression.
        The models with convergence order closest to the theoretical optimal 3 are selected. 
        These are better than the conventional WENO3 schemes.
    }
\end{figure}

\section{Experiments} \label{sec:results}
In this section, we evaluate the performance of our method by comparing it to traditional and ML-based approaches in simulations of varying complexity. A group of WENO3-NN models, chosen for their accuracy, successfully reproduce the flow features with a level of detail comparable to a WENO5 simulation at high resolution.

\subsection{Baselines}
We consider the following baselines for comparison.
\paragraph{Classical methods}
For classical methods, we consider WENO3, WENO5, and QUICK~\cite{versteeg2007introduction,jiang1996efficient}.
WENO3 follows the implementation in \cref{sec:weno}. We point out that for WENO3 and WENO5, we have two different implementations, depending on the experiments. For the small one-dimensional experiment we use a simple implementation in JAX~\cite{jax2018github}, whereas for the more complex numerical simulations we leverage a hand-tuned implementation in TensorFlow~\cite{tensorflow2015}. 

\paragraph{WENO-NN Method} We consider the method introduced in \citet{bezgin2022weno3}, which replaces the rational features by a Delta layer (followed by normalization: $\widetilde{\Delta_j})$ and with Swish activation. The Delta layer calculates the input features of the neural network from the cell averaged values:
\begin{align}
{\Delta_1} &= |\bar{u}_i - \bar{u}_{i-1}|, {\Delta_2} = |\bar{u}_{i+1} - \bar{u}_i|, {\Delta_3} = |\bar{u}_{i+1} - \bar{u}_{i-1}|, {\Delta_4} = |\bar{u}_{i+1} - 2\bar{u}_i + \bar{u}_{i-1}|.\\
\widetilde{\Delta_j} &= {\Delta_j} / \max({\Delta_1}, {\Delta_2}, \epsilon), \forall j= {1,2,3,4},
\end{align}
where $\epsilon=10^{-15}$ is set to avoid division by zero.

\paragraph{Model variants} The networks are trained using a similar pipeline (and parameters) as in \citet{bezgin2022weno3}. We have considered many variants, including random initialization seeds  and network size/architecture configurations. The final list is recorded in \cref{tab:net_hyper_params}).

\paragraph{Computational Cost} We estimate the number of floating point operations (FLOPs) using the JAX internal FLOP estimator (using the XLA~\cite{XLA} compiler) on a CPU device~\cite{jax_aot_webpage}.
Illustrated in \cref{tab:comp_cost}, the FLOPs required for all the WENO-NN methods are higher than the conventional WENO schemes.
This is expected since the neural networks involve multiple matrix--matrix products.
However, our method with the rational network uses six times fewer parameters and four times fewer FLOPs than the approach presented by \citet{bezgin2022weno3}. 
Even though these estimates are hardware specific, the relative difference between the FLOPs of various schemes remain similar.

\begin{table}[t]
\centering
\resizebox{\textwidth}{!}{%
\begin{tabular}{|c|c|c|c|c|l|l|l|}
\hline
\begin{tabular}[c]{@{}c@{}}Model\\  Name\end{tabular} & \begin{tabular}[c]{@{}c@{}}Selection \\ Criteria\end{tabular} & \begin{tabular}[c]{@{}c@{}}Activation \\ Function\end{tabular} & \begin{tabular}[c]{@{}c@{}}Feature \\ Function\end{tabular} & \begin{tabular}[c]{@{}c@{}}No. of \\ Hidden \\ Neurons\end{tabular} & \begin{tabular}[c]{@{}c@{}}Loss \\ $\alpha$\end{tabular} & \begin{tabular}[c]{@{}c@{}}Loss \\ $\beta_d$\end{tabular} & \begin{tabular}[c]{@{}c@{}}Peak \\ Learning \\ Rate\end{tabular} \\ \hline
WENO3-NN-Rational-1 & \multirow{3}{*}{\begin{tabular}[c]{@{}c@{}}Convergence \\ on Sine-step\end{tabular}} & \multirow{12}{*}{Rational} & \multirow{12}{*}{Rational} & \multirow{12}{*}{(4, 4, 4)} & 0.01 & 0.1 & 5E-04 \\ \cline{1-1} \cline{6-8} 
WENO3-NN-Rational-2 &  &  &  &  & 0.03 & 0.03 & 5E-04 \\ \cline{1-1} \cline{6-8} 
WENO3-NN-Rational-3 &  &  &  &  & 0.01 & 0.1 & 1E-04 \\ \cline{1-2} \cline{6-8} 
WENO3-NN-Rational-4 & \multirow{3}{*}{\begin{tabular}[c]{@{}c@{}}Convergence\\  on Sin$^3$\end{tabular}} &  &  &  & 0.1 & 0.3 & 5E-04 \\ \cline{1-1} \cline{6-8} 
WENO3-NN-Rational-5 &  &  &  &  & 0.01 & 0.3 & 5E-04 \\ \cline{1-1} \cline{6-8} 
WENO3-NN-Rational-6 &  &  &  &  & 0.3 & 0.1 & 5E-04 \\ \cline{1-2} \cline{6-8} 
WENO3-NN-R-RMSE-1 & \multirow{3}{*}{\begin{tabular}[c]{@{}c@{}}Least \\ Reconstruction \\ Loss\end{tabular}} &  &  &  & 0.01 & 0.03 & 5E-04 \\ \cline{1-1} \cline{6-8} 
WENO3-NN-R-RMSE-2 &  &  &  &  & 0.3 & 0.03 & 1E-05 \\ \cline{1-1} \cline{6-8} 
WENO3-NN-R-RMSE-3 &  &  &  &  & 0.3 & 0.03 & 1E-05 \\ \cline{1-2} \cline{6-8} 
WENO3-NN-R-RMSE-4 & \multirow{3}{*}{\begin{tabular}[c]{@{}c@{}}Least \\ Deviation \\ Loss\end{tabular}} &  &  &  & 0.1 & 0.3 & 1E-04 \\ \cline{1-1} \cline{6-8} 
WENO3-NN-R-RMSE-5 &  &  &  &  & 0.3 & 0.1 & 1E-04 \\ \cline{1-1} \cline{6-8} 
WENO3-NN-R-RMSE-6 &  &  &  &  & 0.3 & 0.3 & 1E-04 \\ \hline
WENO3-NN-Delta-1 & \citet{bezgin2022weno3} & \multirow{2}{*}{Swish} & \multirow{2}{*}{Delta} & \multirow{2}{*}{(16, 16, 16)} & 0.03 & 0.1 & 1E-05 \\ \cline{1-2} \cline{6-8} 
WENO3-NN-Delta-2 & \citet{bezgin2022weno3} &  &  &  & 0.1 & 0.1 & 1E-05 \\ \hline
\end{tabular}%
}
\caption{Variants of neural network based WENO3 schemes.}
\label{tab:net_hyper_params}
\end{table}

\begin{table}[t]
\centering
\small
\begin{tabular}{l c c c c r}
Source &
  \begin{tabular}[c]{@{}c@{}}Activation\\ Function\end{tabular} &
  \begin{tabular}[c]{@{}c@{}}Feature\\ Function\end{tabular} &
  \begin{tabular}[c]{@{}c@{}}Number of\\ Neurons\end{tabular} &
  \begin{tabular}[c]{@{}c@{}}Number of \\ Parameters\end{tabular} &
  \begin{tabular}[c]{@{}c@{}}Floating Point \\ Operations\end{tabular} \\ \midrule
  WENO3-JS &
  --- &
  --- &
  --- &
  --- &
  19 \\ 
  WENO5-JS &
  --- &
  --- &
  --- &
  --- &
  55 \\ 
\citet{bezgin2022weno3} &
  Swish &
  Delta &
  (16, 16, 16) &
  658 &
  2139 \\ 
Current Work &
  Rational &
  Rational &
  (4, 4, 4) &
  105 &
  508 
\end{tabular}%
\caption{Comparison of models with delta features and rational features.}
\label{tab:comp_cost}
\end{table}

\subsection{One-dimensional problems} 


We show that our methodology exhibits better dispersion and dissipation properties than other classical and ML-based methods with the same stencil width, while rivaling the performance of higher-order schemes. We choose the one-dimensional Burgers' and linear advection equations to showcase the long-time accuracy and shock capturing capabilities of the scheme.

\subsubsection{Linear Advection}

Consider the linear advection problem: 
\begin{align} \label{eq:advection}
     \partial_t u(x, t) + \partial_x u(x, t) &= 0, \qquad \text{for } (x, t) \in [0,1]\times [0, T], \\
     u(x, 0)& = u_0(x) \text{ on } x \in [0, 1].
\end{align}
with periodic boundary conditions, and initial condition, $u_0$, given by
\begin{align}
    \text{Cosine: } u_0(x) &= \cos(2 \pi x), \\
    \text{Sigmoid: } u_0(x) &= \left(1 + \exp(-k (x-x_1)) \right)^{-1} + \left(1 + \exp(k (x-x_2)) \right)^{-1}.
\end{align}
where, $k=100$, $x_1=0.05$ and $x_2=0.2$ in the Sigmoid function.

The first initial condition is designed to showcase that our methodology does not suffer from spurious dissipation at critical points~\cite{henrick2005mapped}, whereas the second one is designed to showcase that our methodology represents high-gradient solution with little spurious dissipation.

\begin{figure}[t]
     \centering
     \begin{subfigure}[b]{0.49\textwidth}
         \centering
         \includegraphics[width=\textwidth]{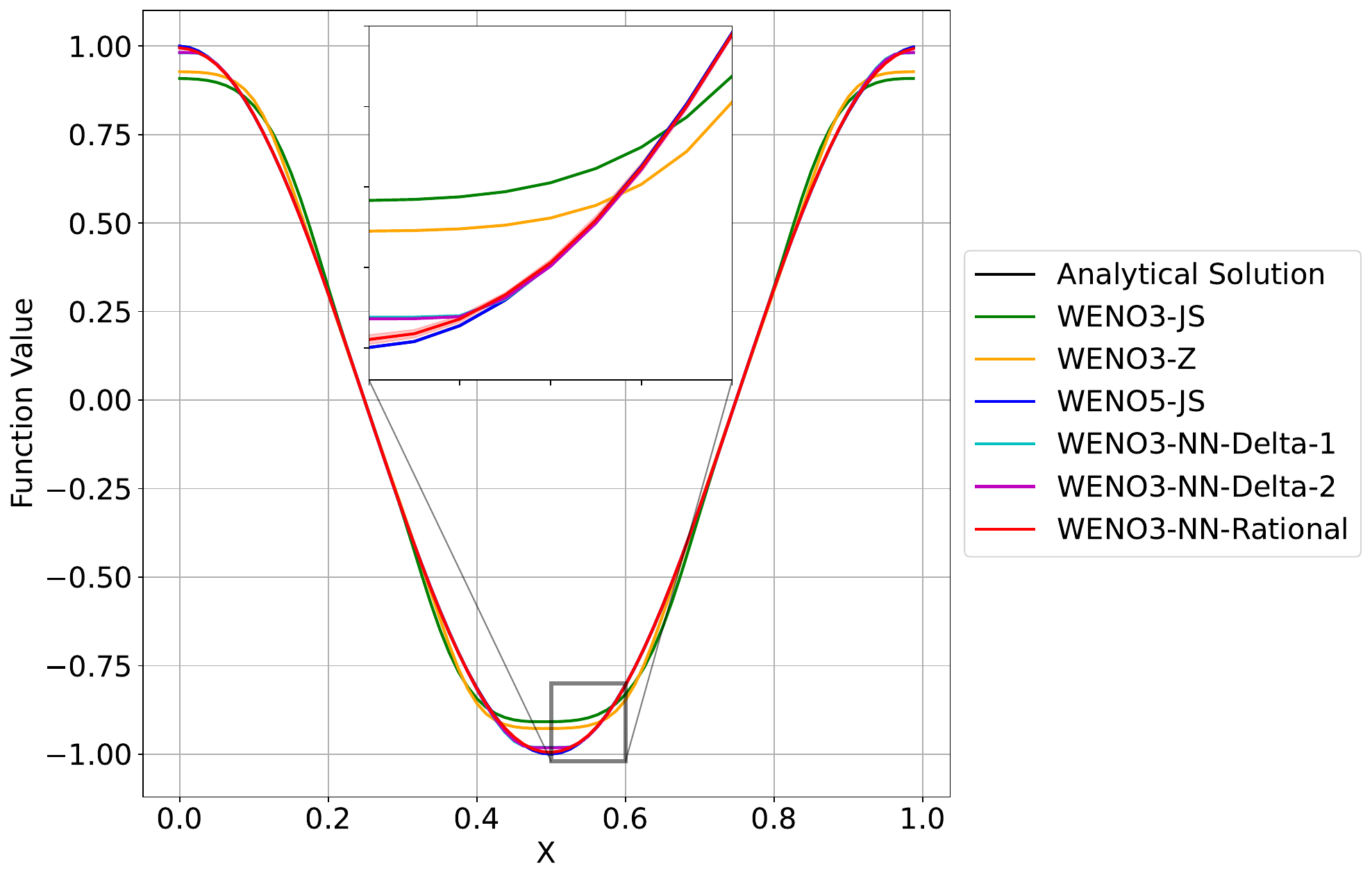}
         \caption{Solution at last time step}
         \label{fig:Advection Cosine Soln}
     \end{subfigure}
     \begin{subfigure}[b]{0.49\textwidth}
         \centering
         \includegraphics[width=\textwidth]{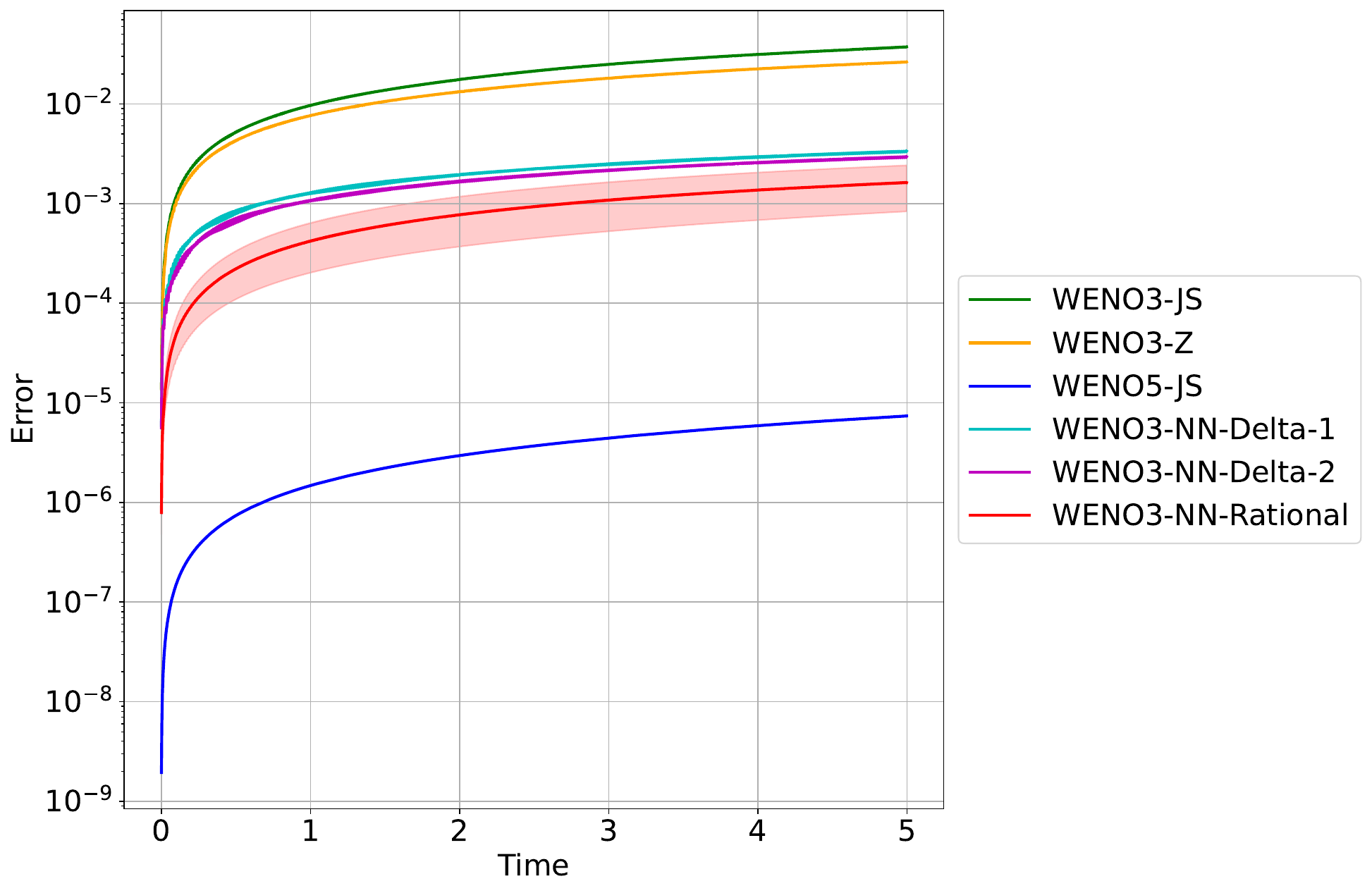}
         \caption{Temporal variation of solution error}
         \label{fig:Advection Cosine Error}
     \end{subfigure}
    \caption{Advection of cosine wave (shaded region depicts the confidence interval of WENO3-NN-Rational-1 to WENO3-NN-Rational-6 in \cref{tab:net_hyper_params}).}
    \label{fig:Advection Cosine}
\end{figure}

For the cosine initial condition, we solve the PDE for a time horizon $T = 5$ with different methods, including several instances of our trained models. \Cref{fig:Advection Cosine Soln} depicts the solution for each method at time horizon $t = T$, i.e., after 5 flow through times. From \Cref{fig:Advection Cosine Soln} we see that all the ML-based models were able to represent the solution accurately, besides both the WENO3 variants, which are overly dissipative, even at this very low frequency. 

We compute the $L^{1}$ error in space of the numerical solution with respect to the analytical solution at each time step which are summarized in \cref{fig:Advection Cosine Error}.
We use a solid red line to represent the average prediction for all the rational networks (WENO3-NN-Rational-1 to WENO3-NN-Rational-6 in \cref{tab:net_hyper_params}), where the red shading shows the confidence interval. \Cref{fig:Advection Cosine Error} shows the evolution of the error for different methods. In summary, WENO5-JS achieves the lowest error due to its larger stencil size. Although WENO3 and WENO3-NN methods share the same stencil size, the WENO3-NN methods produce errors at least 10 times smaller than WENO3.

\begin{figure}[H]
     \centering
     \begin{subfigure}[b]{0.49\textwidth}
         \centering
         \includegraphics[width=\textwidth]{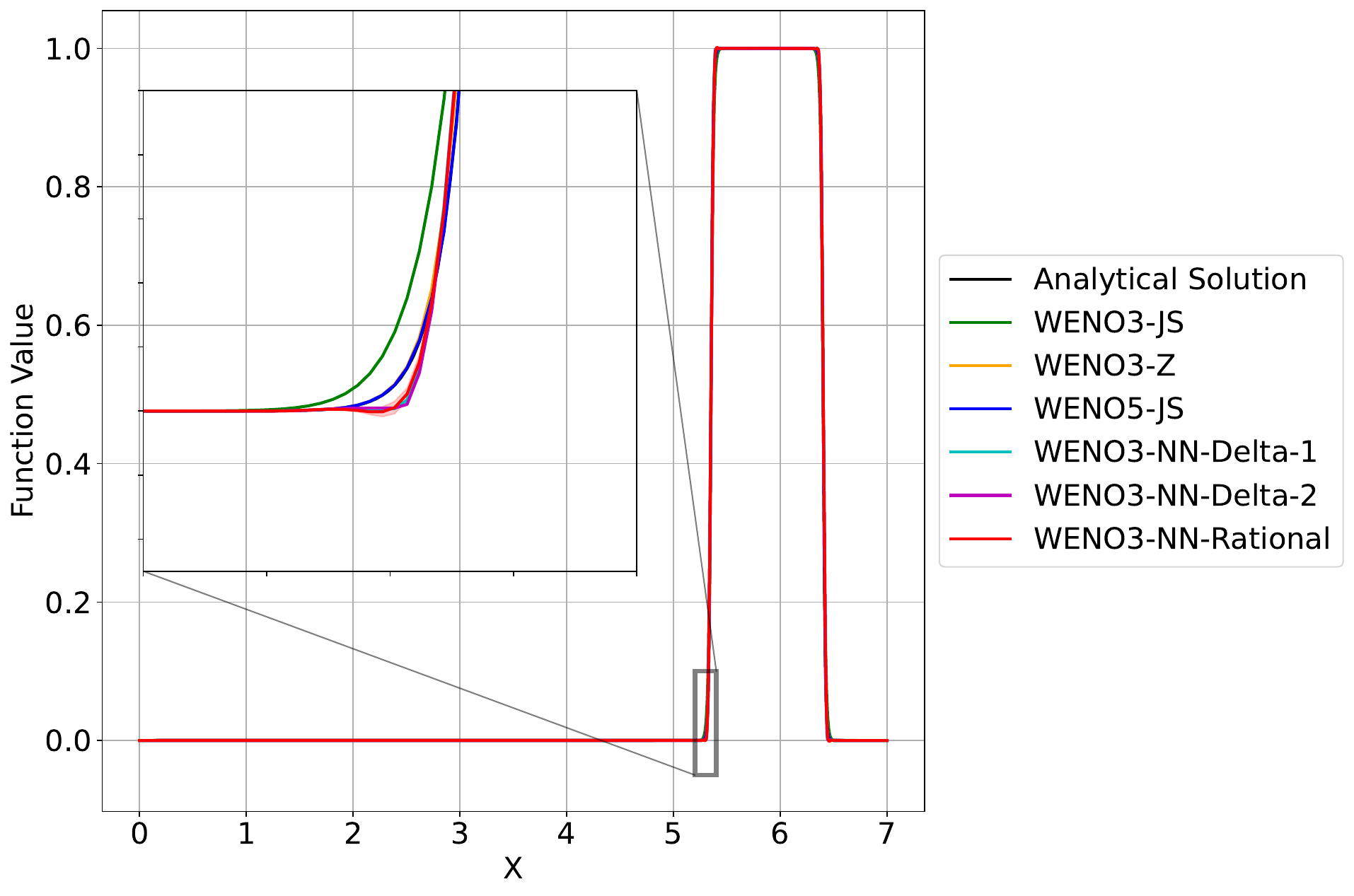}
         \caption{Solution at last time step.}
         \label{fig:Advection Sigmoid Soln}
     \end{subfigure}
     \begin{subfigure}[b]{0.49\textwidth}
         \centering
         \includegraphics[width=\textwidth]{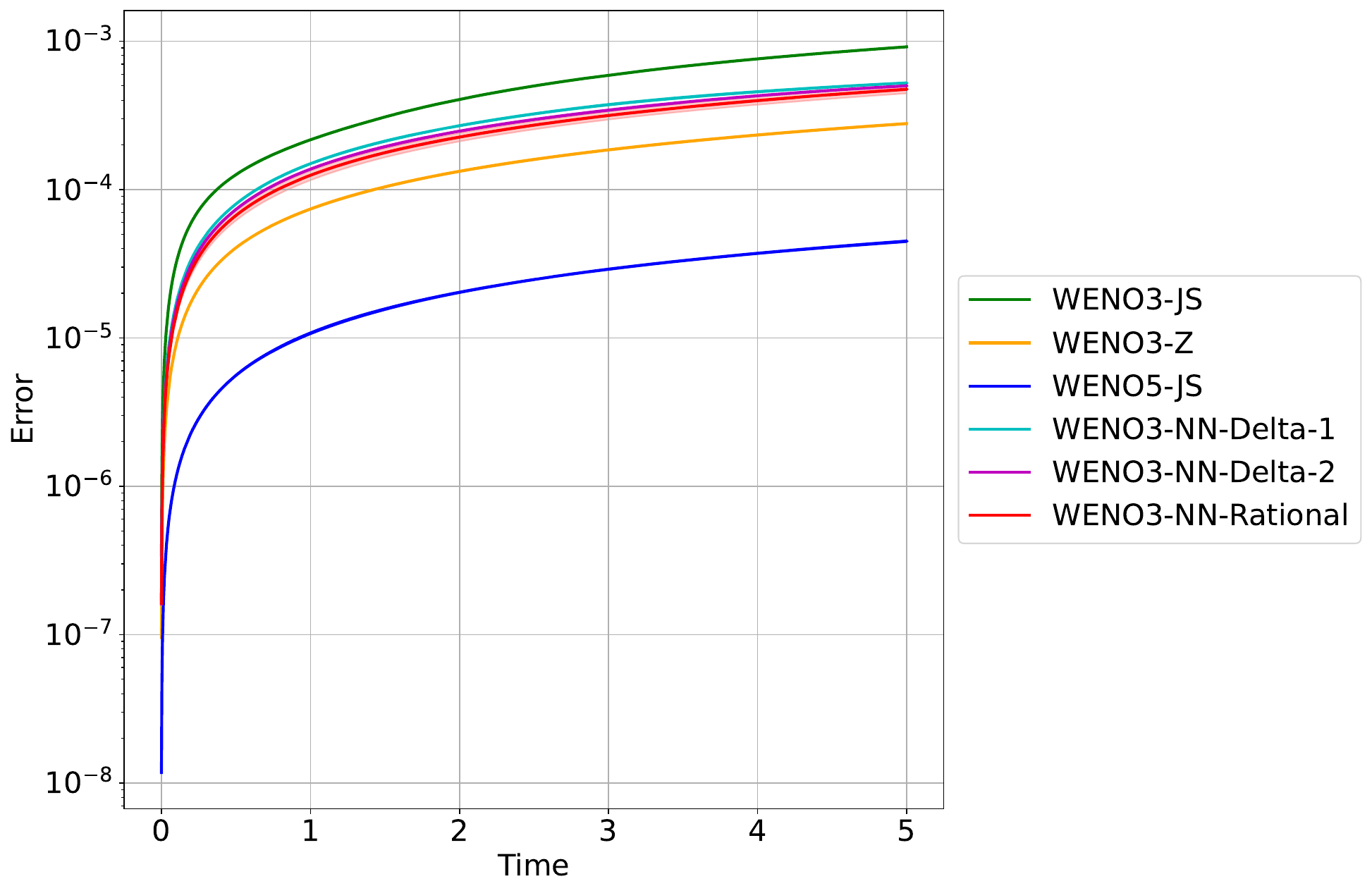}
         \caption{Temporal variation of solution error.}
         \label{fig:Advection Sigmoid Error}
     \end{subfigure}
    \caption{Advection of sigmoid wave (shaded region depicts the confidence interval of WENO3-NN-Rational-1 to WENO3-NN-Rational-6 in \cref{tab:net_hyper_params}).}
    \label{fig:Advection Sigmoid}
\end{figure}

For the sigmod initial condition, we solve the PDE for the same time horizon $T = 5$. \Cref{fig:Advection Sigmoid} shows the solution and error to \cref{eq:advection} at time $T = 5$, where we see that, as in the cosine case, all the ML-based models are able to represent the solution accurately, besides the WENO3-JS, which is also overly dissipative. We observe that the ML-methods outperform WENO3-JS while underperfoming with respect to WENO3-Z and WENO5-JS. Due to a high order global smoothness measure, WENO3-Z scheme is known to be less dissipative than WENO3-JS~\cite{acker2016improved}. In the case of WENO5-JS, the larger stencil allows for better estimation of sharp gradients.

\begin{figure}[H]
     \centering
     \begin{subfigure}[b]{0.49\textwidth}
         \centering
         \includegraphics[width=\textwidth]{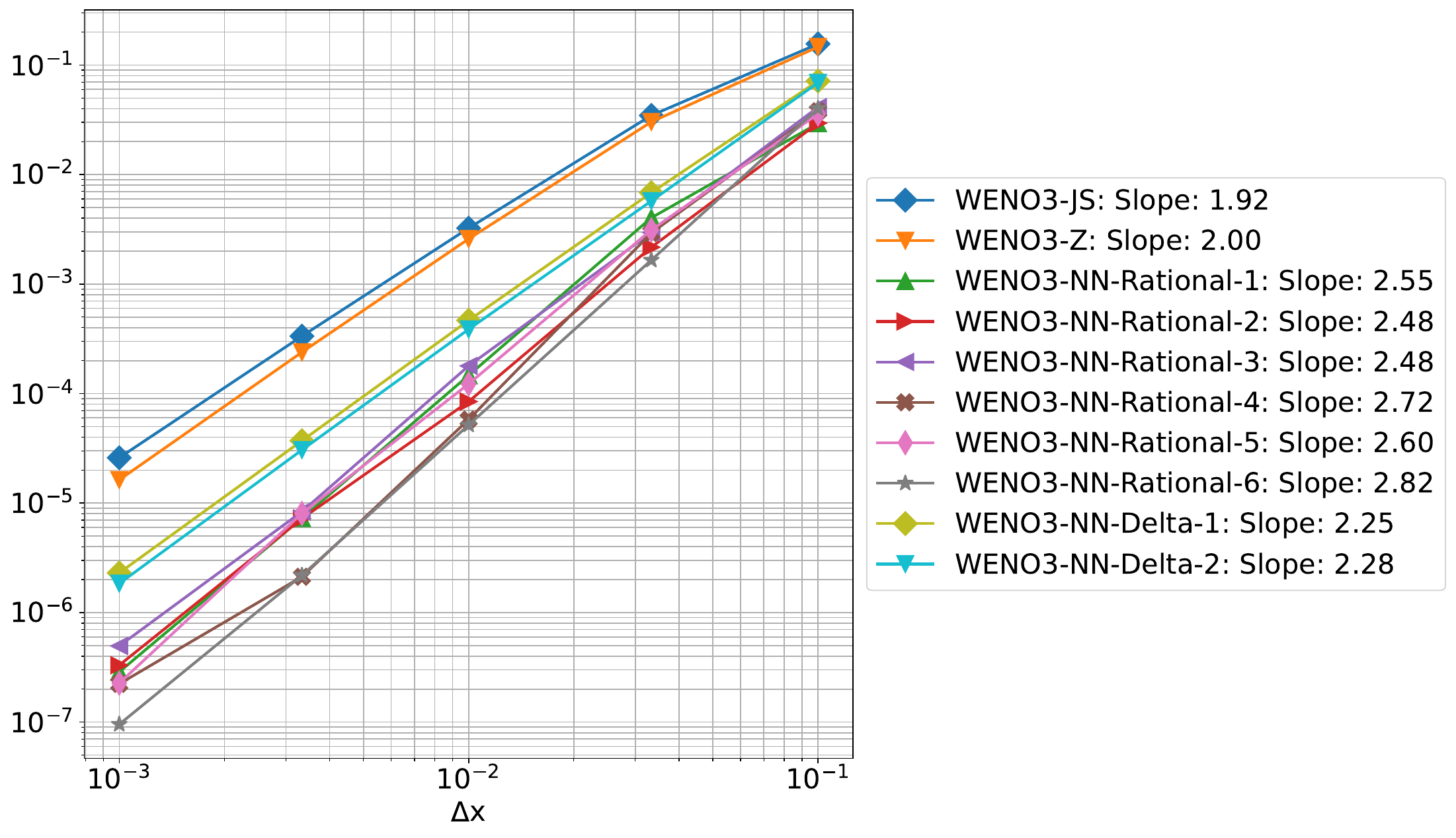}
         \caption{Advection of cosine wave.}
         \label{fig:Scalar Tranport Convergence: Cosine}
     \end{subfigure}
     \begin{subfigure}[b]{0.49\textwidth}
         \centering
         \includegraphics[width=\textwidth]{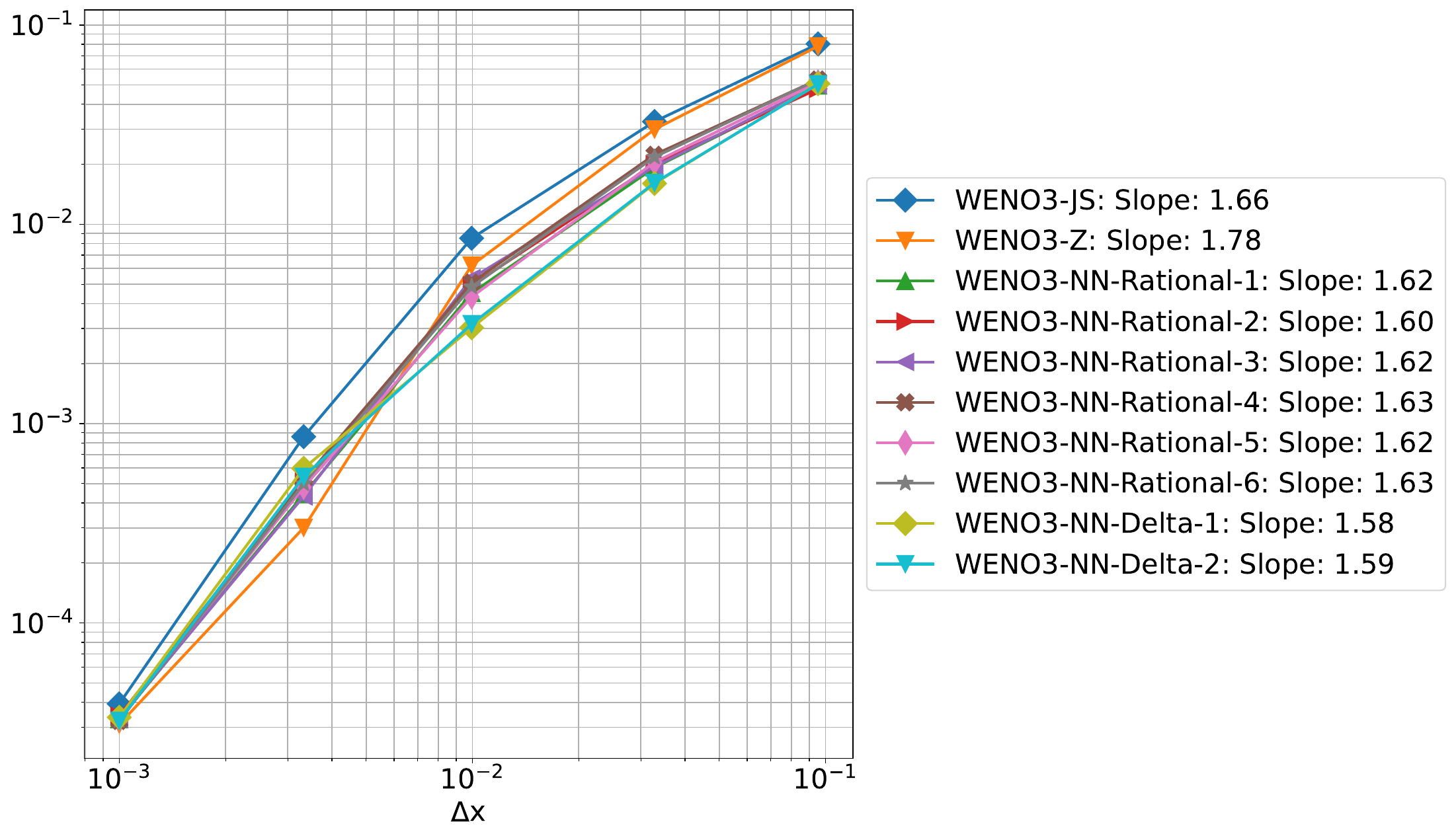}
         \caption{Advection of sigmoid wave.}
         \label{fig:Scalar Tranport Convergence: Sigmoid}
     \end{subfigure}
    \caption{Convergence behavior of different instances of the proposed methods and baselines (slope derived from linear regression).}
    \label{fig:Scalar Tranport Convergence}
\end{figure}

For both initial conditions, we performed a convergence analysis in space with respect to the analytical solution (\cref{fig:Scalar Tranport Convergence}). From \cref{fig:Scalar Tranport Convergence} we see that all WENO3-NN based methods are more accurate than the conventional WENO3-JS method.
From \cref{fig:Scalar Tranport Convergence: Cosine}, we observe that for the cosine wave, WENO3-NN schemes are able to exceed the second order convergence (which is expected for WENO3 schemes~\cite{jiang1996efficient}) while presenting a considerable advantage against traditional methods.
From \cref{fig:Scalar Tranport Convergence: Sigmoid},  we observe that for the sigmoid wave, the order of the convergence for all the methods hovers around $1.6$, which is attributed to the high-gradients in the sigmoid wave. However, in this case the performance gaps narrows, as
 WENO3-Z scheme is slightly more accurate than WENO3-NN alternatives at finer grids.
 However, for coarse meshes, which is our target regime, WENO3-NN schemes are more accurate by a factor 2.

\subsubsection{Inviscid Burgers' Equation}

We showcase the properties of our methodology for nonlinear equations with discontinuous solutions. We consider the inviscid Burgers' equation as a prototypical example of this category, which is given by
\begin{align} \label{eq:burgers}
     \partial_t u(x, t) + \frac{1}{2}\partial_x u^2(x, t) &= 0, \qquad \text{for } (x, t) \in [-6,6]\times [0, T], \\
     u(x, 0)& = u_0(x) \text{ on } x \in [-6, 6],
\end{align}
plus Dirichlet boundary conditions,  where $u_0$ is an initial condition given by the step function\footnote{This type of initial condition is also found in Riemann problems~\cite{leveque2002finite}.}:
\begin{equation}
  u_0(x)=\begin{cases}
    u_l, & \text{if $x<0$},\\
    u_r, & \text{otherwise}.
  \end{cases}
\end{equation}
The boundary conditions are prescribed to be consistent with the initial condition, i.e.,
\begin{align}
     u(-6, t) = u_l, \,\, \text{and }  u(6, t) = u_r.
\end{align}

We consider three representative solution scenarios for the Riemann problem above:
\begin{enumerate}
    \item Shock wave ($u_l>u_r$): $u_l=1$, $u_r=0$.
    \item Rarefaction wave ($0 \leq u_l<u_r$): $u_l=0$, $u_r=1$.
    \item Transonic rarefaction wave ($u_l<0<u_r$): $u_l=-1$, $u_r=1$.
\end{enumerate}
For all cases, we integrate Burgers' equation up to time $T=5$. The exact solution for the shock wave is given by~\cite{leveque2002finite}:
\begin{equation}
  u(x,t)=\begin{cases}
    u_l, & \text{if $x<c_0t$}.\\
    u_r, & \text{otherwise}.
  \end{cases}
\end{equation}
where, $c_0  = (u_l + u_r) / 2$ satisfy the Rankine-Hugoniot conditions. The exact solution for the initial condition of both rarefaction waves is
\begin{equation}
  u(x, t)=\begin{cases}
    u_l, & x<u_l t.\\
    x/t, & u_l t \leq x \leq u_r t.\\
    u_r, & x>u_l t.
  \end{cases}
\end{equation}

\begin{figure}[H]
     \centering
     \begin{subfigure}[b]{0.49\textwidth}
         \centering
         \includegraphics[width=\textwidth]{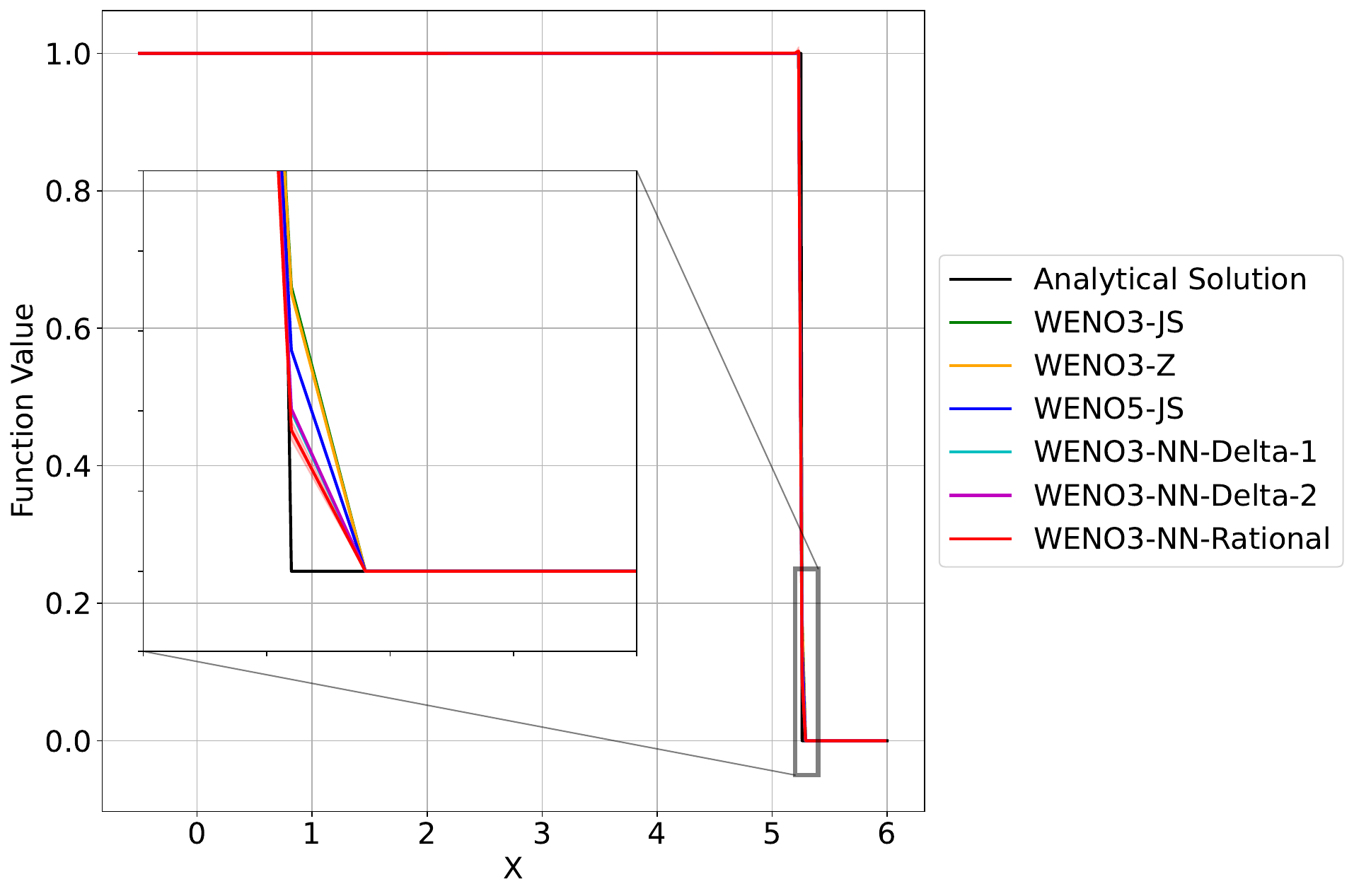}
         \caption{Solution at last time step}
         \label{fig:Burgers Shock Soln}
     \end{subfigure}
     \begin{subfigure}[b]{0.49\textwidth}
         \centering
         \includegraphics[width=\textwidth]{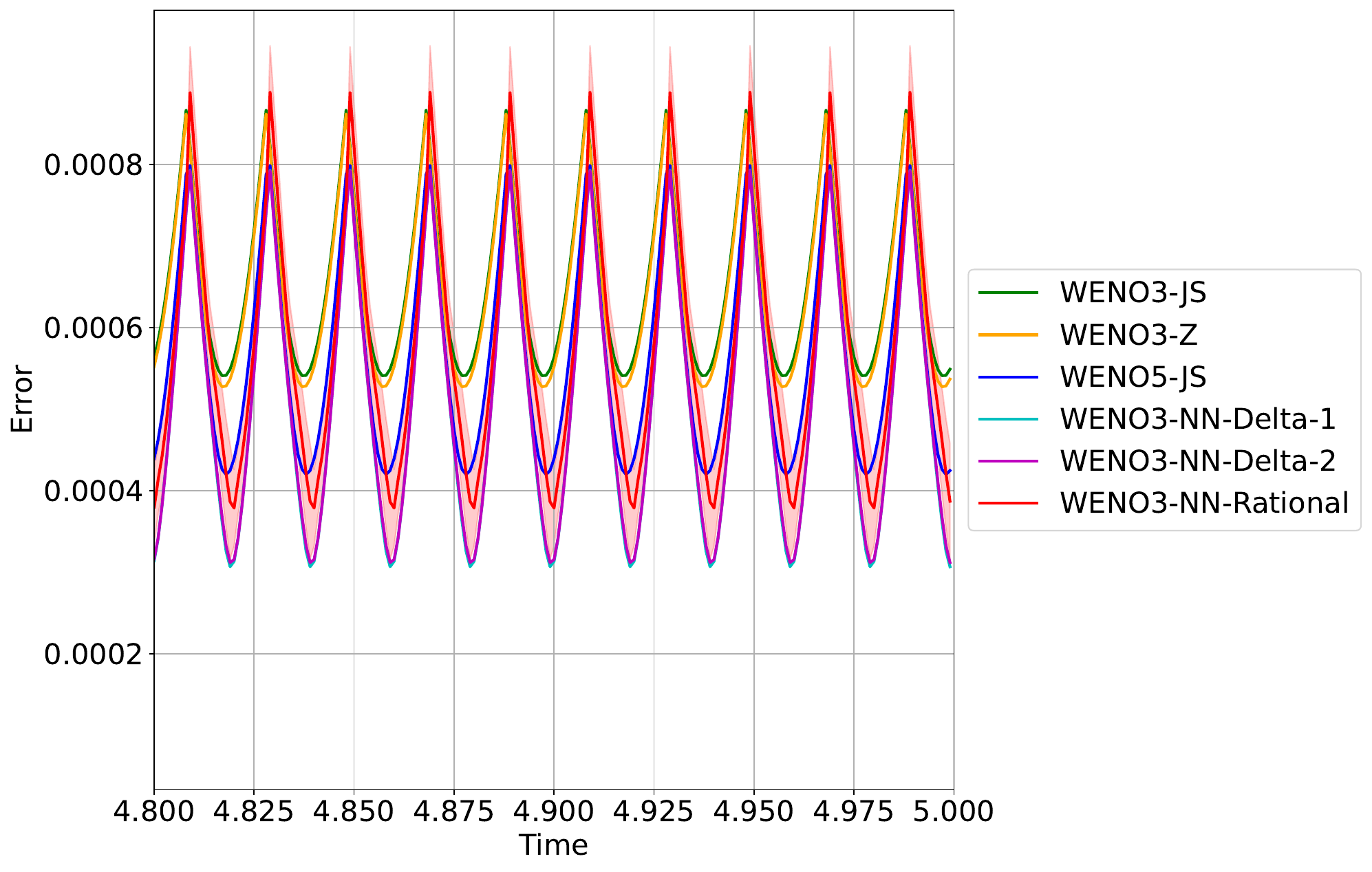}
         \caption{Temporal variation of solution error}
         \label{fig:Burgers Shock Error}
     \end{subfigure}
    \caption{Inviscid Burgers' equation: shock wave (shaded region depicts the confidence interval of WENO3-NN-Rational-1 to WENO3-NN-Rational-6 in \cref{tab:net_hyper_params})}
    \label{fig:Burgers Shock}
\end{figure}

As seen in \cref{fig:Burgers Shock}, WENO3-Z, WENO3-JS, WENO5-JS and WENO-NN methods are able to resolve the sharp discontinuity in the analytical solution.
\Cref{fig:Burgers Shock Error} shows that the WENO-NN models improve the prediction accuracy significantly compared to WENO3. 
We also observe that some of our models have errors smaller than or equal to WENO5-JS, which has a wider stencil.

\begin{figure}[H]
     \centering
     \begin{subfigure}[b]{0.49\textwidth}
         \centering
         \includegraphics[width=\textwidth]{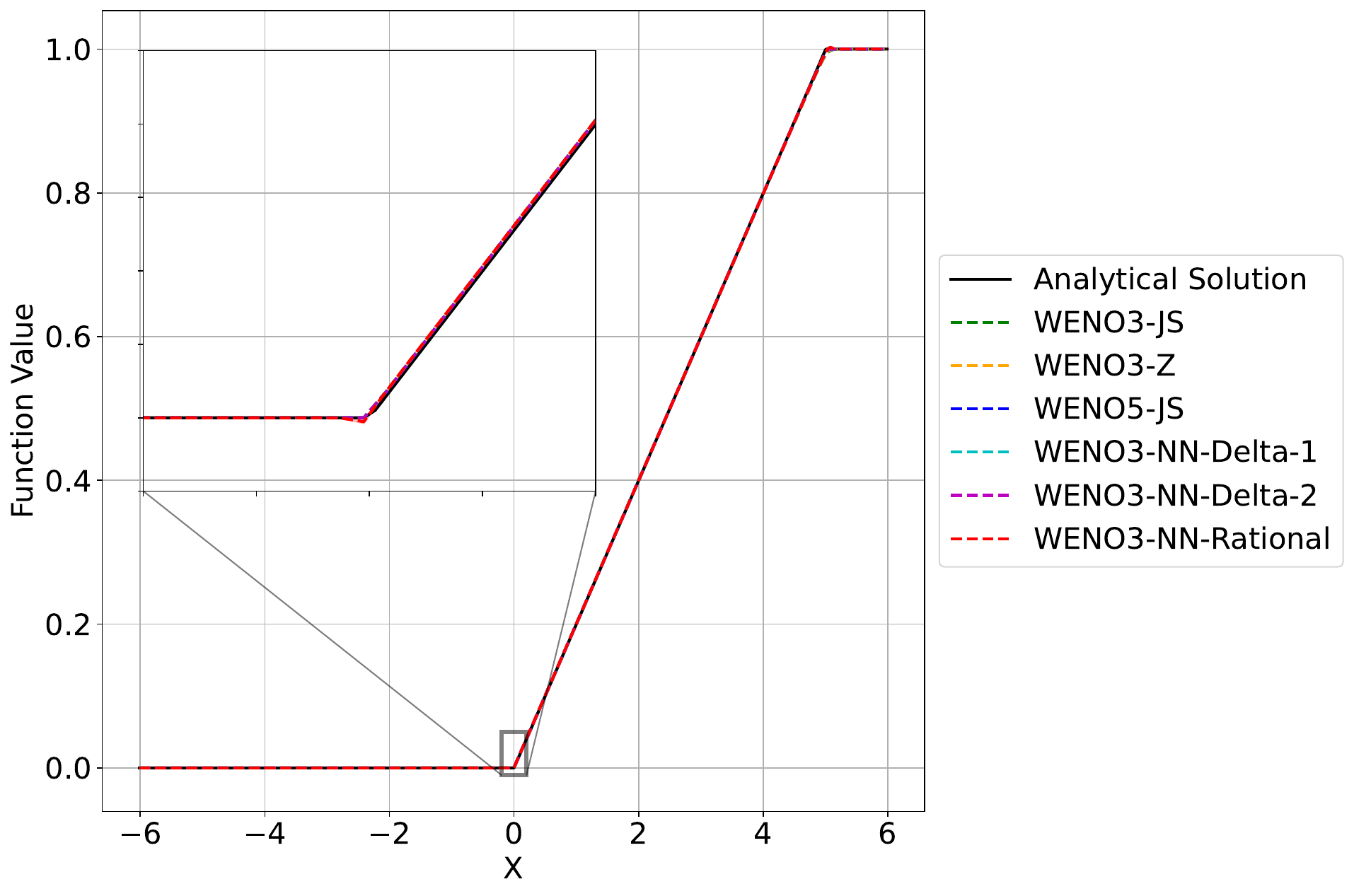}
         \caption{Solution at last time step}
         \label{fig:Burgers Rarefaction Soln}
     \end{subfigure}
     \begin{subfigure}[b]{0.49\textwidth}
         \centering
         \includegraphics[width=\textwidth]{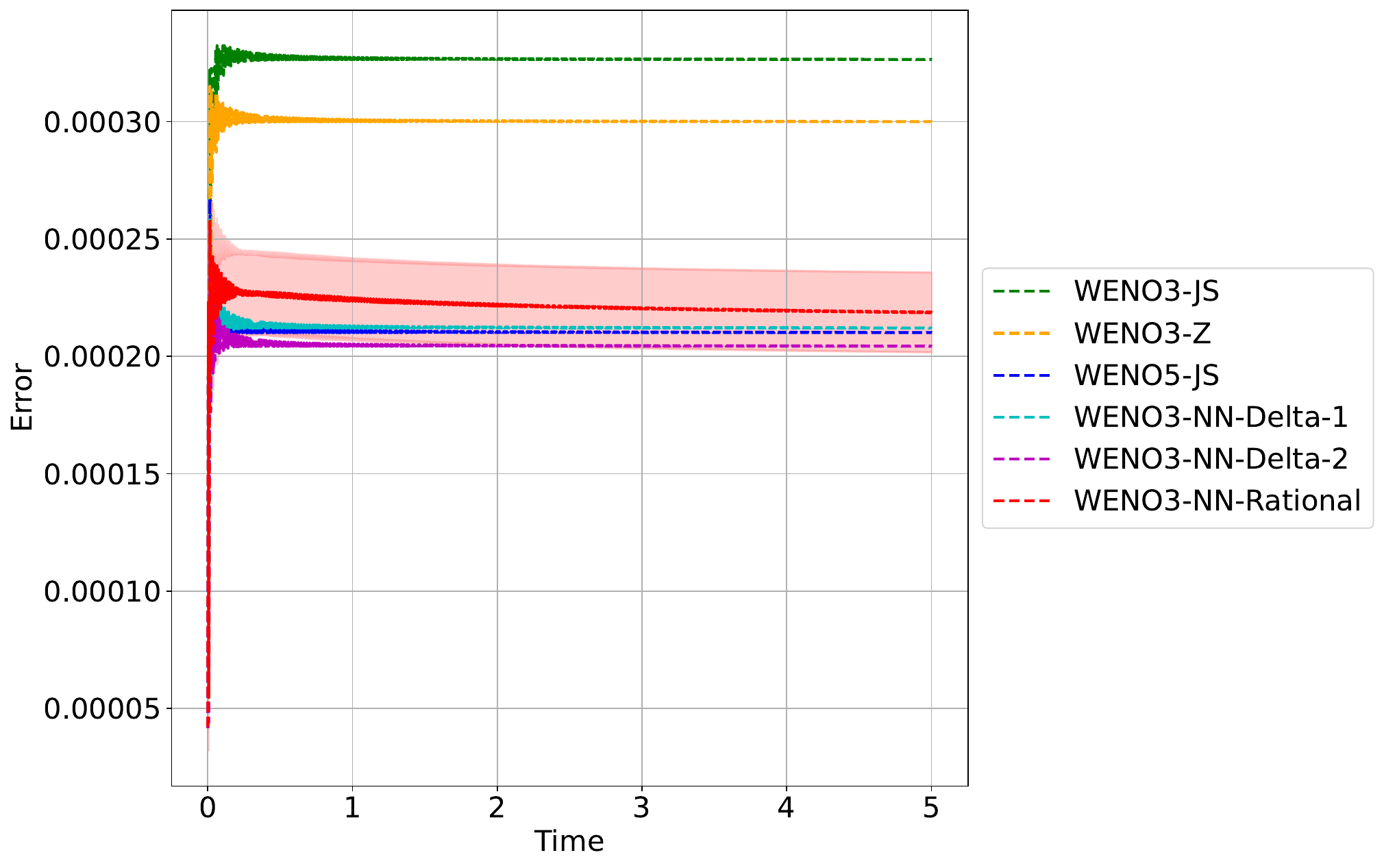}
         \caption{Temporal variation of solution error}
         \label{fig:Burgers Rarefaction Error}
     \end{subfigure}
    \caption{Inviscid Burgers' equation: rarefaction wave (shaded region depicts the confidence interval of WENO3-NN-Rational-1 to WENO3-NN-Rational-6 in \cref{tab:net_hyper_params})}
    \label{fig:Burgers Rarefaction}
\end{figure}

Similar observations carry on for the case of rarefaction waves (\cref{fig:Burgers Rarefaction}). Solutions computed using some of WENO-NN models exhibit errors close to WENO5-JS and both classical WENO3 schemes have approximately 40\% higher error than these models.

\begin{figure}[H]
     \centering
     \begin{subfigure}[b]{0.49\textwidth}
         \centering
         \includegraphics[width=\textwidth]{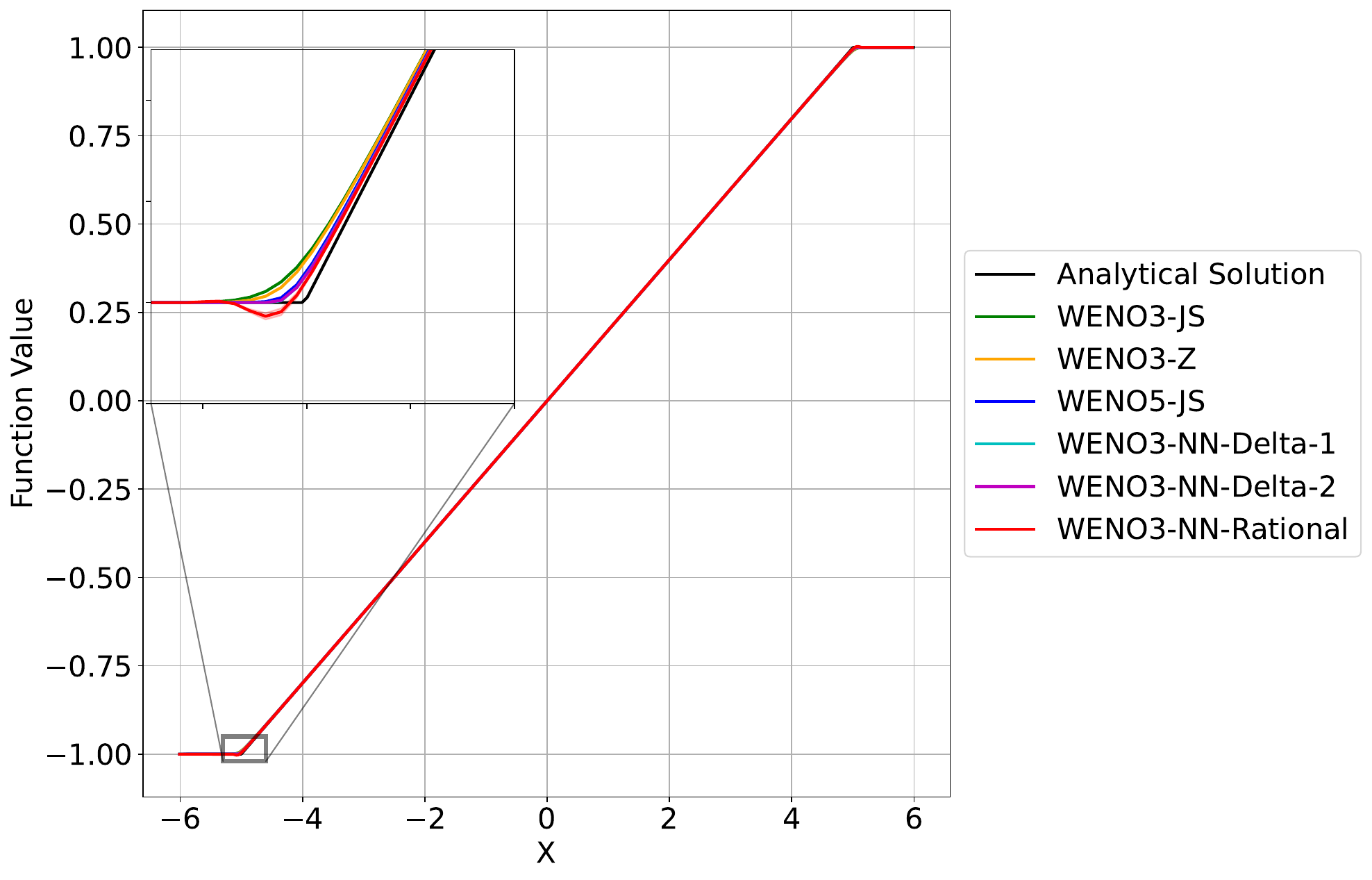}
         \caption{Solution at last time step}
         \label{fig:Burgers Transonic Rarefaction Soln}
     \end{subfigure}
     \begin{subfigure}[b]{0.49\textwidth}
         \centering
         \includegraphics[width=\textwidth]{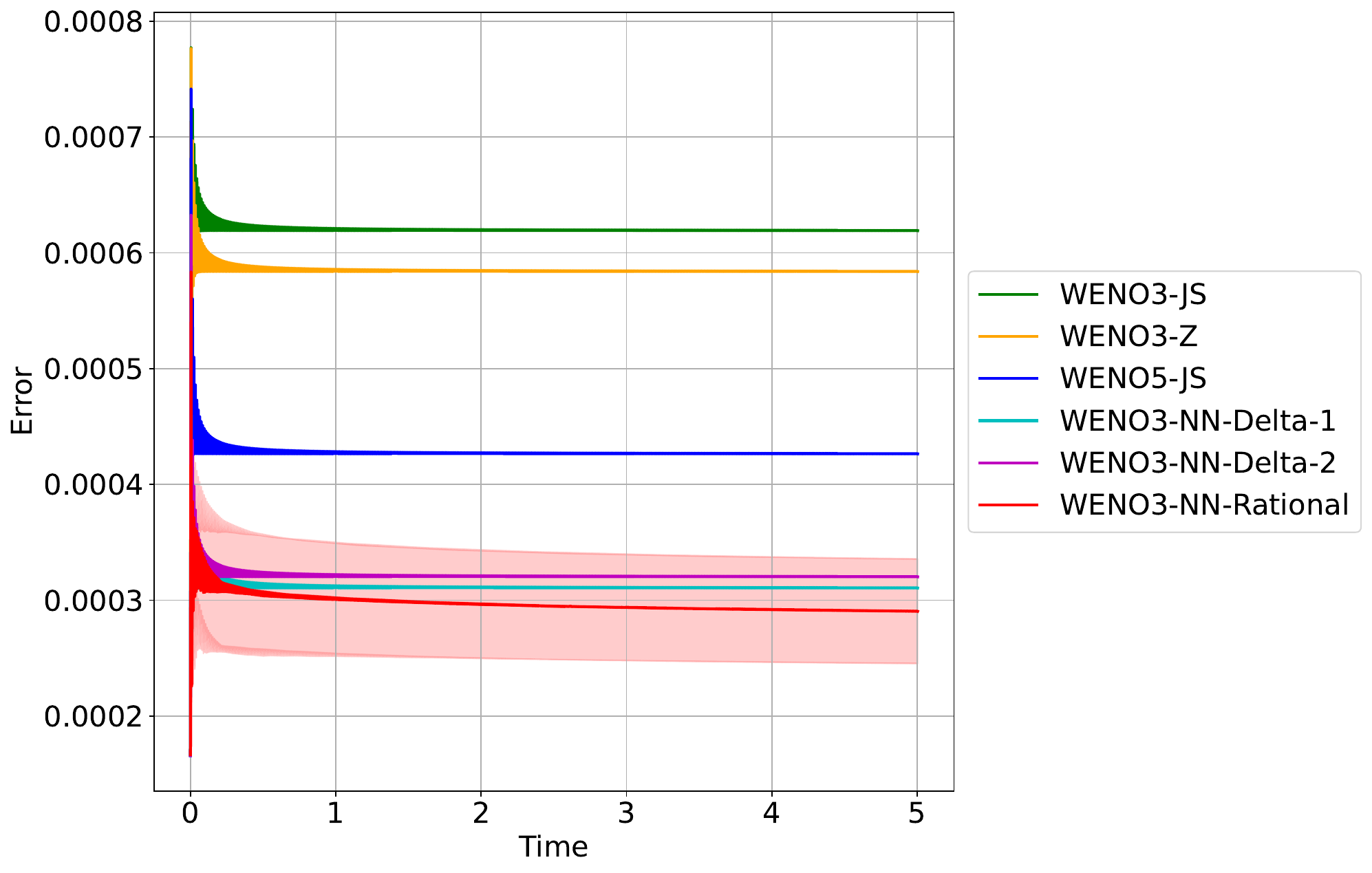}
         \caption{Temporal variation of solution error}
         \label{fig:Burgers Transonic Rarefaction Error}
     \end{subfigure}
    \caption{Inviscid Burgers' equation: transonic rarefaction wave (shaded region depicts the confidence interval of WENO3-NN-Rational-1 to WENO3-NN-Rational-6 in \cref{tab:net_hyper_params})}
    \label{fig:Burgers Transonic Rarefaction}
\end{figure}

\Cref{fig:Burgers Transonic Rarefaction} plots the solution and errors for the transonic rarefaction wave. The error of WENO5-JS is smaller than both classical WENO3 schemes. However, we observe that all the WENO3-NN models are more accurate than WENO5-JS.

\subsubsection{Model selection criterion}
\begin{figure}[H]
    \centering
    \begin{subfigure}[t]{0.28\textwidth}
         \centering
         \includegraphics[width=\textwidth, clip, trim={0mm 0mm 101mm 0mm}]{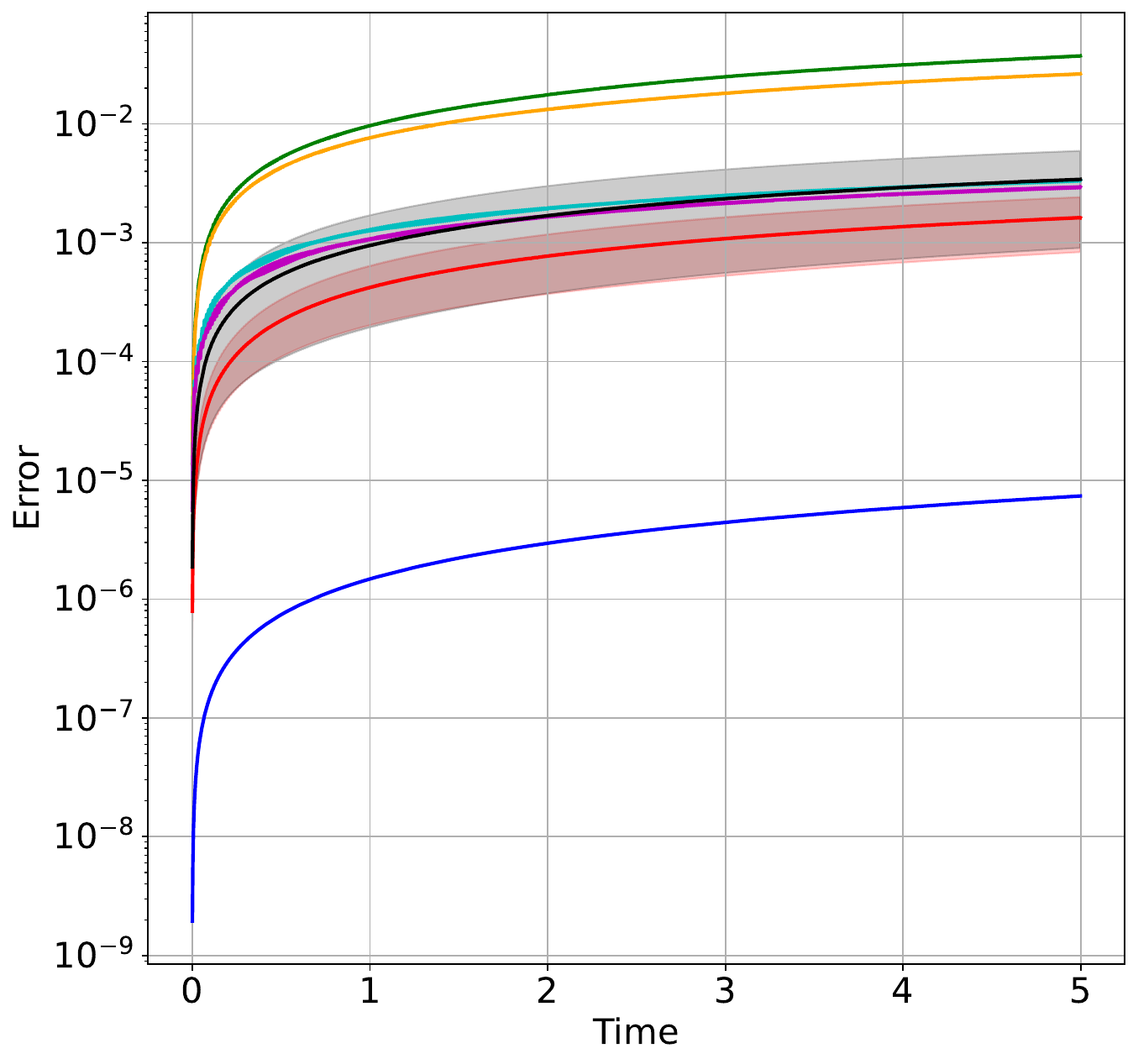}
         \caption{Advection of cosine wave.}
    \end{subfigure}
    \begin{subfigure}[t]{0.28\textwidth}
         \centering
         \includegraphics[width=\textwidth, clip, trim={0mm 0mm 101mm 0mm}]{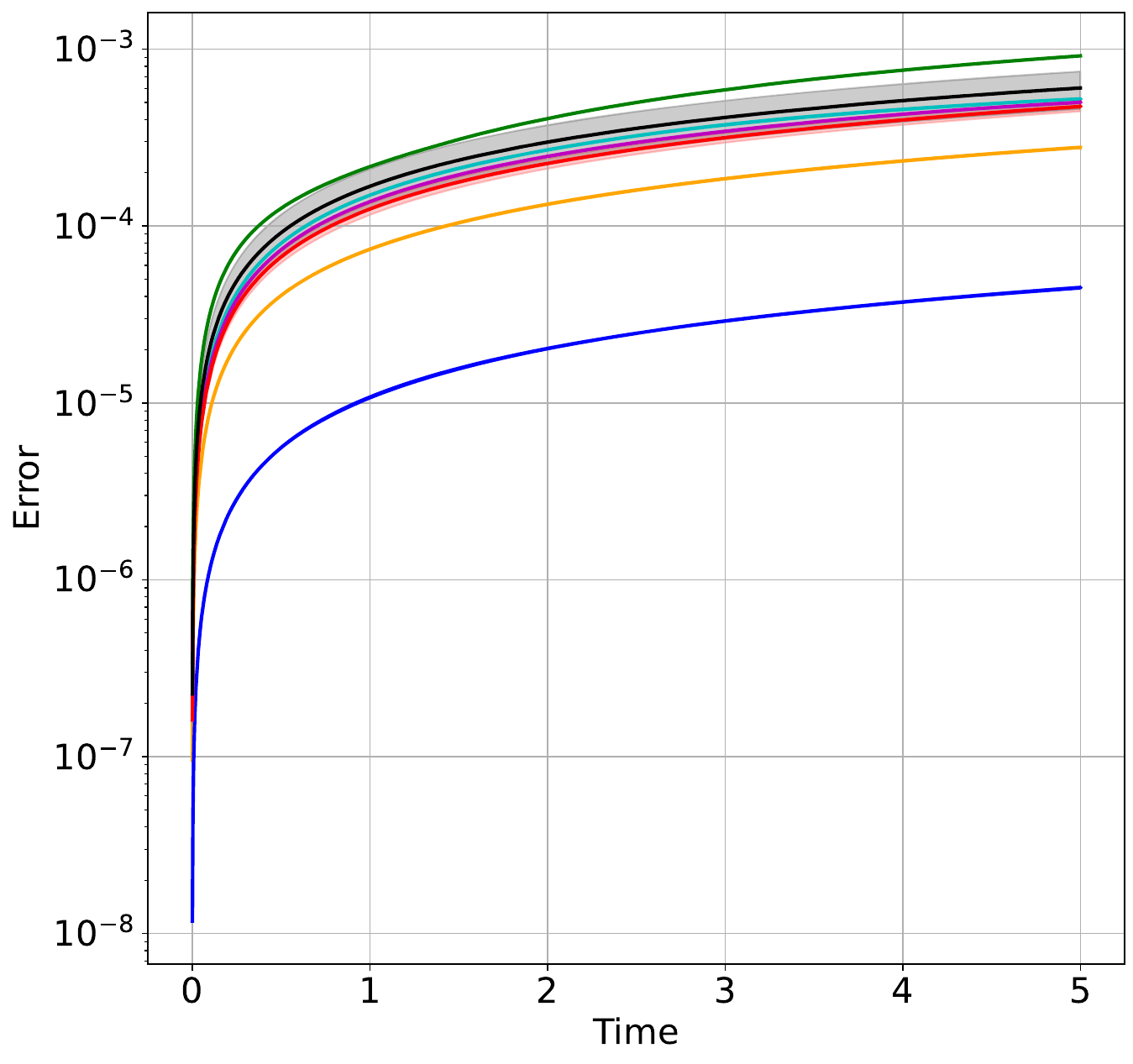}
         \caption{Advection of sigmoid wave.}
    \end{subfigure}
    \begin{subfigure}[t]{0.42\textwidth}
         \centering
         \includegraphics[width=\textwidth]{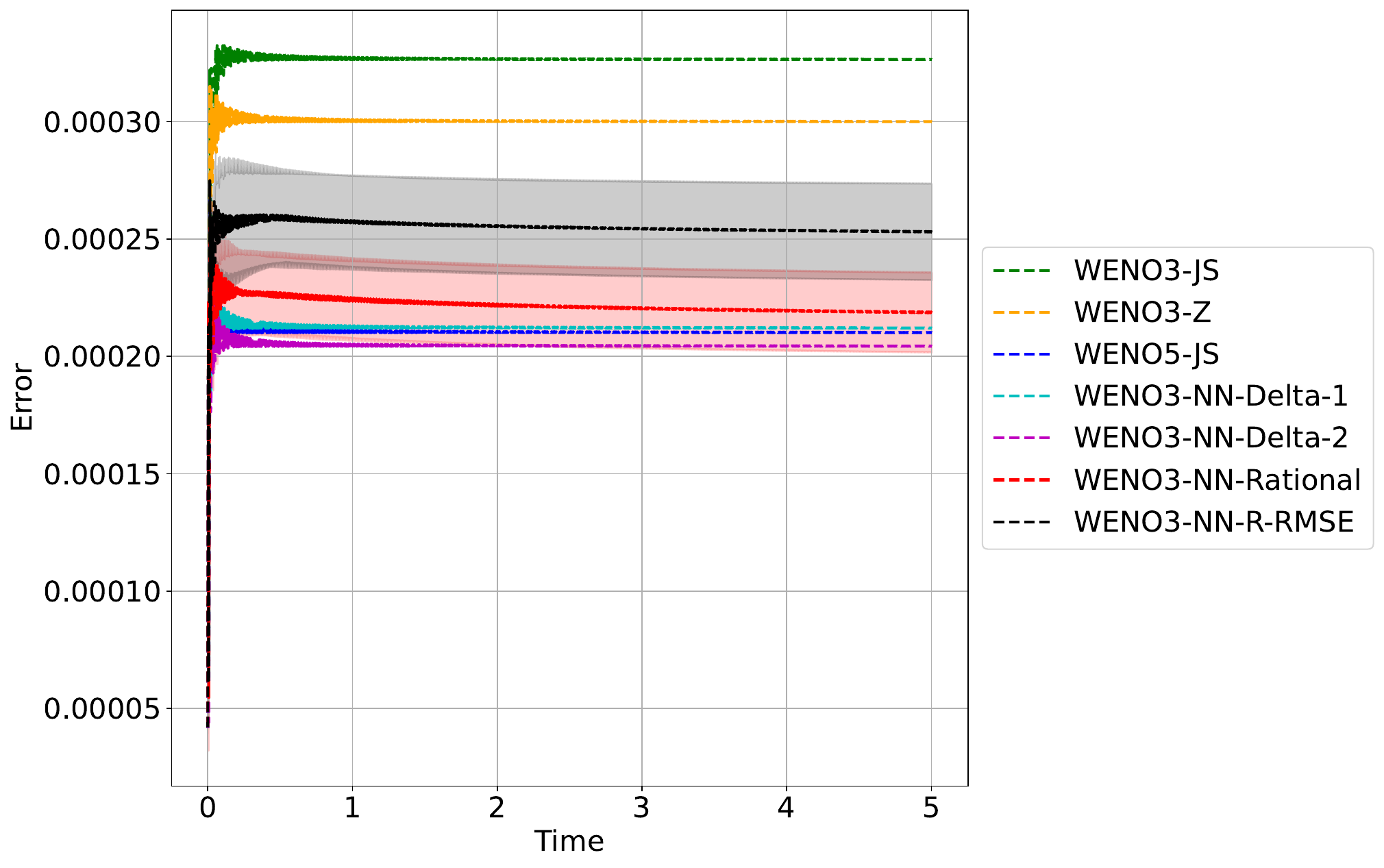}
         \caption{Inviscid Burgers' equation: rarefaction wave.}
    \end{subfigure}
    \caption{Rational neural networks chosen on the basis of convergence criterion have smaller error compared to the ones chosen based on least mean squared error fit on the test functions.}
    \label{fig:Error RMSE models}
\end{figure}

Here, we compare the performance of the models selected with the selection criteria introduced in \cref{sec:model_selection} on the task of solving the Burgers' and linear advection equations. We performed the same experiments as above, using trained models chosen with different criteria. From \cref{fig:Error RMSE models} we see that, while all WENO3-NN models have higher accuracy than WENO3-JS), the models chosen using the convergence criterion outperform the rest (that use the selection criterion listed in \cref{tab:net_hyper_params}). Therefore, for the rest of the paper, we restrict our results to models chosen based on the order of convergence.

\subsubsection{Dispersion-dissipation analysis}
\begin{figure}[H]
    \centering
    \includegraphics[width=0.8\textwidth]{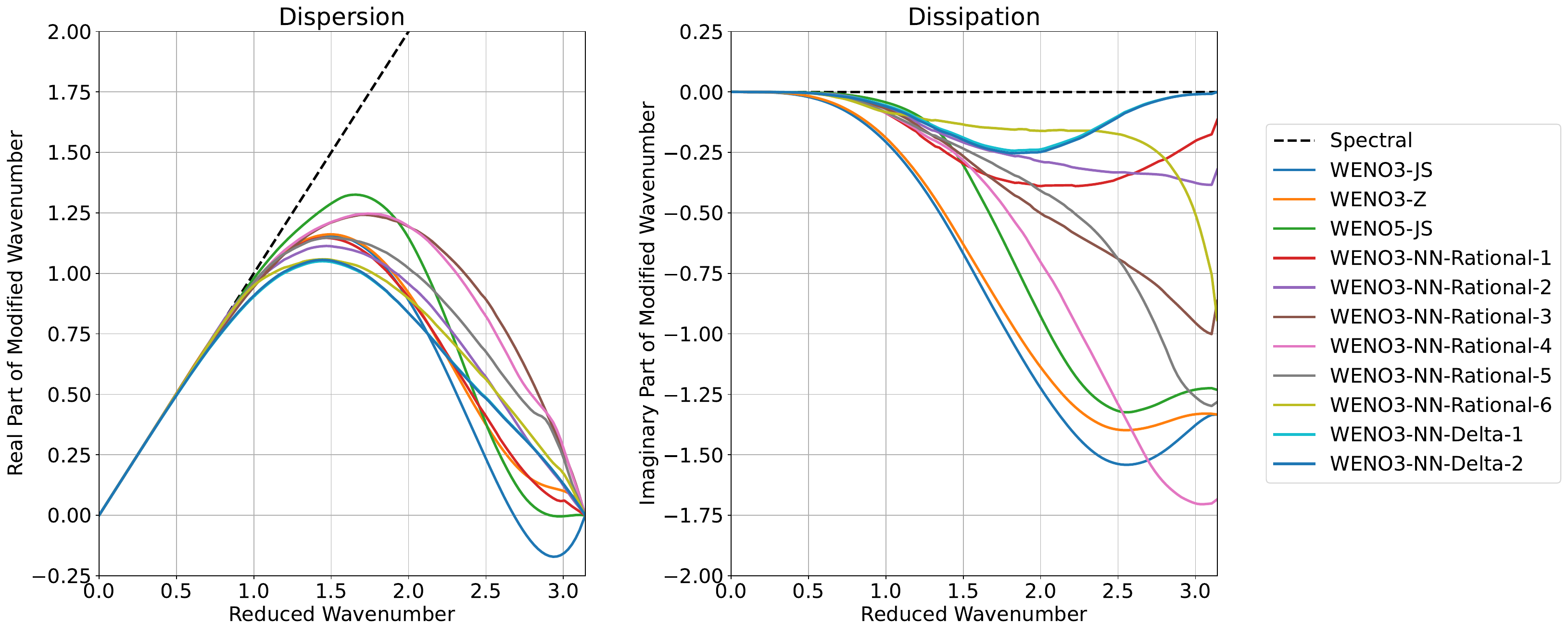}
    \caption{Modified dispersion-dissipation relation for the advection equation for different WENO schemes.}
    \label{fig:WENO-ADR}
\end{figure}
\Cref{fig:WENO-ADR} depicts the dispersion and dissipation analysis for the linear advection equation following~\cite{pirozzoli2006spectral}.
\Cref{fig:WENO-ADR} shows that the modified dispersion relations of all conventional and ML-based methods are closely aligned with the ideal spectral lines for lower wavenumbers. In the mid-range wavenumber region, WENO3-NN schemes exhibit superior spectral properties compared to WENO3 schemes, although they fall slightly short of WENO5. However, at higher wavenumbers, all WENO3-NN schemes outperform WENO5 in terms of spectral properties. Noticeably, WENO3-NN introduces less dissipation near wavenumber cutoff.

\subsection{Two-dimensional problems} 

Here, we apply the WENO-NN models to two-dimensional fluid dynamics simulations.

\subsubsection{Buoyant bubble with smooth initial condition} \label{sec:smooth bubble}

We simulate a thermal bubble rising under buoyancy~\cite{bryan2004reevaluation,chammas2023accelerating}. The objective is to test the developed WENO-NN method on the nonlinear Navier--Stokes equations coupled with scalar transport equations. We show that the predictions from carefully chosen WENO-NN models simulated on a coarse grid of $512^2$ are closer to the fine grid WENO5-JS simulations at $2048^2$.
In contrast, the conventional WENO3-JS method simulated on a coarse grid of $512^2$ is highly dissipative.

Here, we summarize\footnote{For further details, we refer the interested reader to~\cite{chammas2023accelerating}.} the governing equations for the buoyant thermal bubble.
\begin{equation}
    \frac{\partial \rho}{\partial t} + \nabla \cdot (\rho \bm{u}) = 0,
    \label{eq:continuity}
\end{equation}
\begin{equation}
    \frac{\partial (\rho \bm{u})}{\partial t} + \nabla \cdot (\rho \bm{u} \bm{u}) = - \nabla(p') + b \hat{\bm{g}} + \nabla \cdot \tau + \nabla \cdot (\rho \sigma),
    \label{eq:momentum}
\end{equation}
where, $\bm{u}$: velocity vector, $\rho$: density of the fluid, $p'$: the hydrodynamic pressure, $\tau$: stress tensor, $\sigma$: subgrid scale stress per unit mass, $\hat{\bm{g}}$: unit vector in the direction of gravity and $b$: buoyancy term given in \cref{eq:buoyancy}. We solve for scalars: the liquid-ice potential temperature ($\theta_{li}$) and total humidity $q_t$ using the transport equations:
\begin{equation}
    \frac{\partial (\rho \theta_{li})}{\partial t} + \nabla \cdot (\rho \bm{u} \theta_{li}) = \frac{1}{Pr} \nabla \cdot (\rho \nu_t \nabla \theta_{li}),
    \label{eq:theta_l_transport}
\end{equation}

\begin{equation}
    \frac{\partial (\rho q_t)}{\partial t} + \nabla \cdot (\rho \bm{u} q_t) = \frac{1}{Sc_{q_t}} \nabla \cdot (\rho \nu_t \nabla q_t)
    \label{eq:q_t_transport}
\end{equation}
where, $Pr$ is the turbulent Prandtl number, $\nu_t$ is turbulent viscosity given by the model of Lilly~\cite{lilly1962numerical} and Smagorinsky~\cite{smagorinsky1963general} and $Sc_{q_t}$ is the turbulent Schmidt number of water.
The buoyancy is given by
\begin{equation}
    b = g (\rho(\theta_l, q_t, p_0) - \rho_0(z)),
    \label{eq:buoyancy}
\end{equation}
where, $\rho_0$ is a reference density that only depends on $z$, the spatial coordinate in the direction of gravity.

The initial condition of the potential temperature in Kelvin is given by:
\begin{align}
    \theta(t=0) &= 
    \begin{cases}
    300 + 2 \cos(\pi r / 2)^2, & \text{if $r<1$}.\\
    300, & \text{otherwise}.
  \end{cases} \label{eq:smooth bubble} \\
  r &= \sqrt{(x-x_c)^2 + (y-y_c)^2} / r_0.
\end{align}
Center and radius of the bubble are set to $(x_c,y_c) = (10000,2000)$~m and $r_0=2000$~m. We use a square domain with each side of length 20000~m. Gravity is acting in vertically downward direction.

We consider 3 different grid resolutions: $512^2$, $1024^2$ and $2048^2$ with a time step of 0.1, 0.05 and 0.025 seconds. The PDE is then solved up to a time-horizon $T=2000$ s.

We present the resulting simulations of the rising bubble using three conventional methods: QUICK, WENO3-JS, and WENO5-JS. The governing equations include the momentum conservation (\cref{eq:momentum}) and scalar transport equations (\cref{eq:theta_l_transport,eq:q_t_transport}). We simulated the bubble using 7 different combinations of the conventional schemes based on same or different schemes for equations (\cref{tab:bubble-convection-schemes}). In the subsequent figures, `Mtm.' and `Scl.' denote the schemes used for momentum and scalar transport equations respectively.

\begin{figure}
     \centering
     \begin{subfigure}[b]{0.49\textwidth}
         \centering
         \includegraphics[width=\textwidth, clip, trim={0mm 0mm 150mm 0mm}]{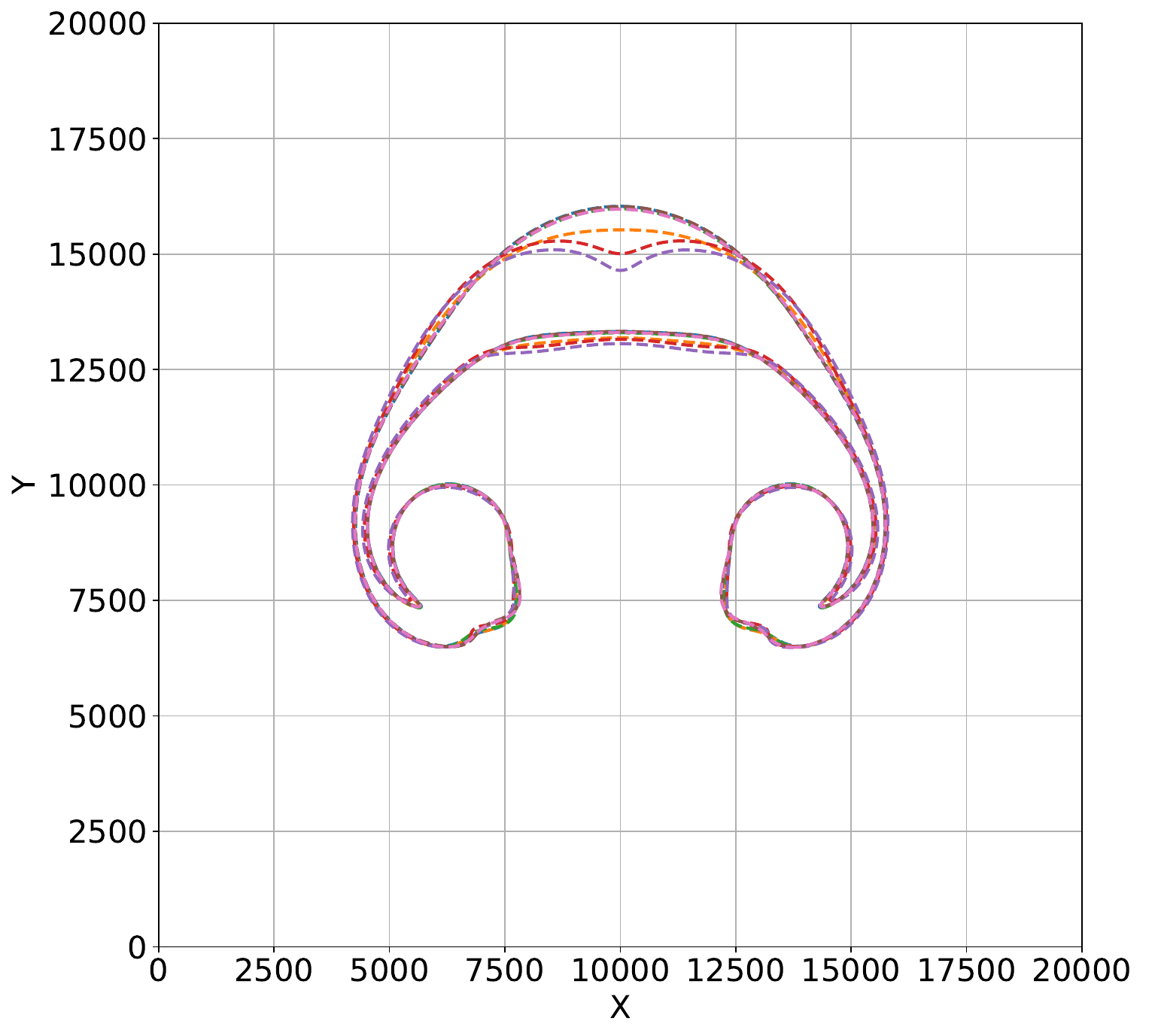}
         \caption{Grid: $512^2$}
         \label{fig:Smooth Bubble conventional 512 300.3}
     \end{subfigure}
     \begin{subfigure}[b]{0.49\textwidth}
         \centering
         \includegraphics[width=\textwidth, clip, trim={0mm 0mm 150mm 0mm}]{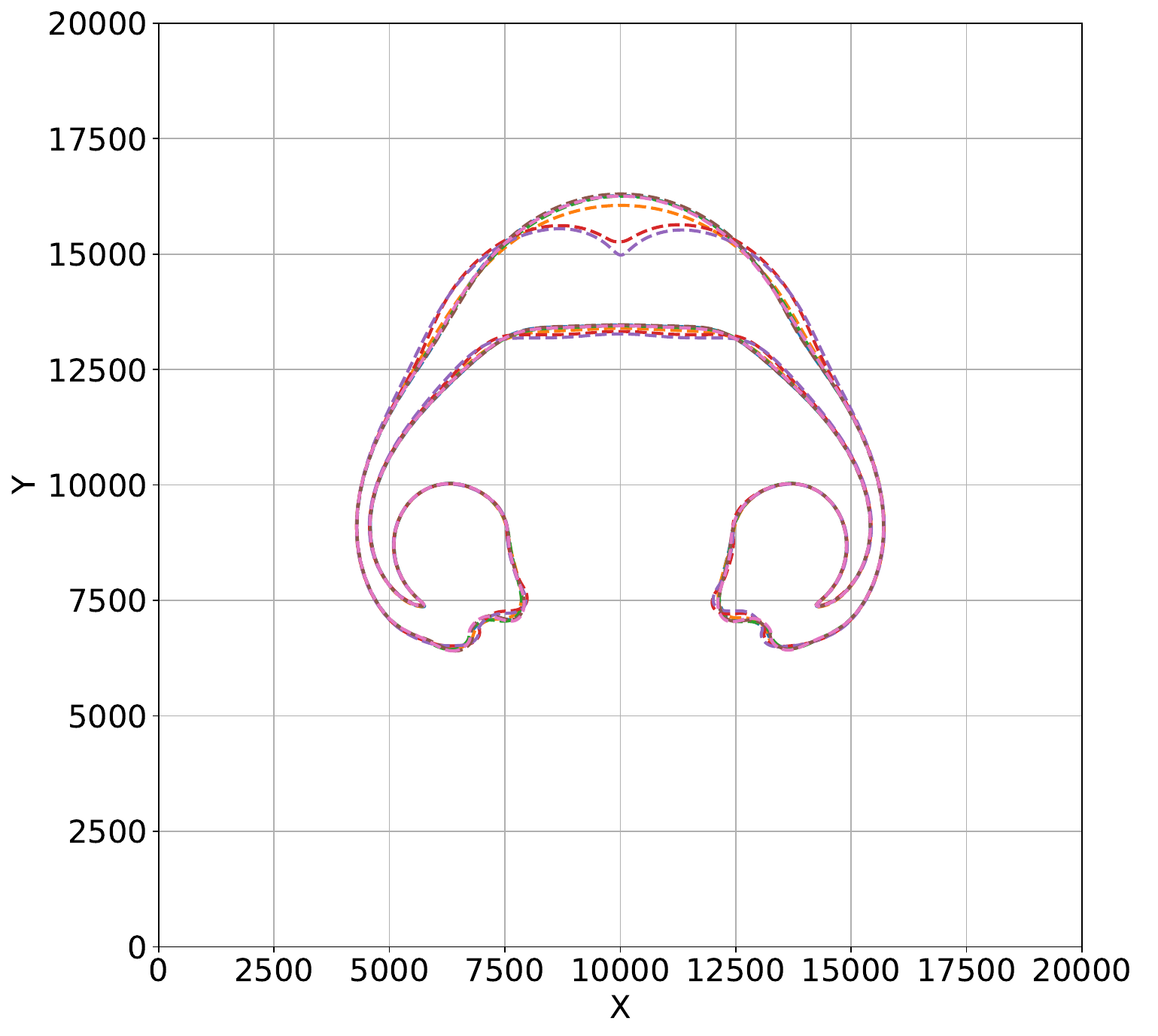}
         \caption{Grid: $1024^2$}
         \label{fig:Smooth Bubble conventional 1024 300.3}
     \end{subfigure}
     \begin{subfigure}[b]{0.75\textwidth}
         \centering
         \includegraphics[width=\textwidth]{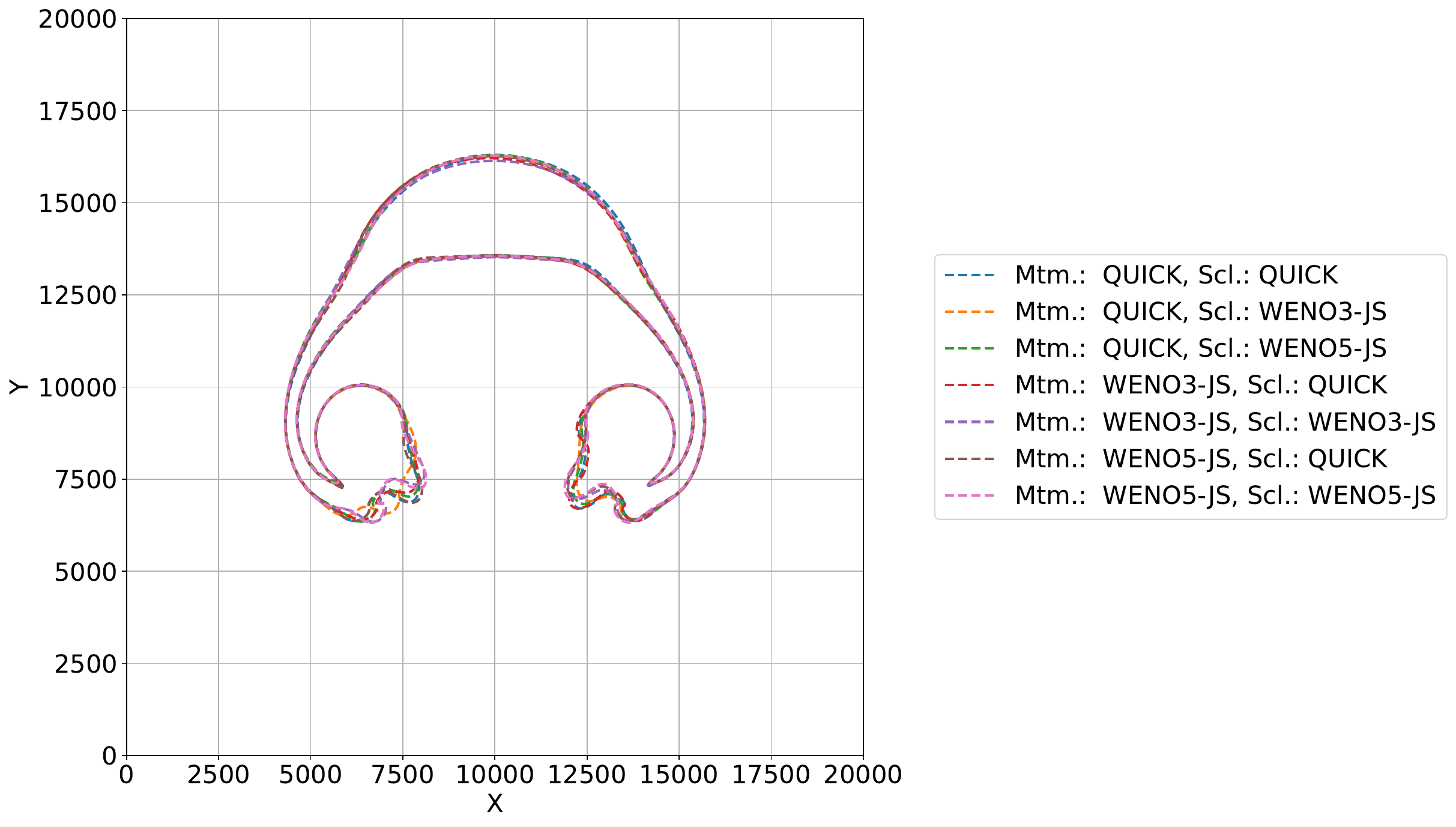}
         \caption{Grid: $2048^2$}
         \label{fig:Smooth Bubble conventional 2048 300.3}
     \end{subfigure}
    \caption{Buoyant bubble with smooth initial condition: potential temperature contour at 300.3~K (conventional schemes).}
    \label{fig:Smooth Bubble conventional 300.3}
\end{figure}

\begin{table}[t]
\centering
\small
\begin{tabular}{l l}
Momentum equations (Mtm.) & Scalar transport equations (Scl.) \\ \midrule
QUICK & QUICK \\
QUICK & WENO3-JS \\
QUICK & WENO5-JS \\ 
WENO3-JS & QUICK \\ 
WENO3-JS & WENO3-JS \\ 
WENO5-JS & QUICK \\ 
WENO5-JS & WENO5-JS 
\end{tabular}%
\caption{Combination of convection schemes for bubble simulations.}
\label{tab:bubble-convection-schemes}
\end{table}
\Cref{fig:Smooth Bubble conventional 300.3} plots the contour lines of the potential temperature at $300.3$~K.  On the one hand, we observe that the simulations using WENO3-JS have higher dissipation and the bubble is deformed near the top edge especially at low resolutions. On the other hand, the solutions produced by QUICK and WENO5-JS schemes do not exhibit this dissipative behaviour at coarser grids of $512^2$ and $1024^2$. As we refine the grid to $2048^2$ in \cref{fig:Smooth Bubble conventional 2048 300.3}, all schemes exhibit an increasingly smoother and convergent contour of the bubble with minimal dissipation, where the top portion of the bubble is well defined.

We present the results of simulations conducted with ML-based models. We consider the simulations using WENO5-JS for both momentum and scalar transport equations with a fine grid of $2048^2$ grid as `ground truth'\footnote{This corresponds to the highest accuracy among our simulations.}. We provide a visual comparison between the contours produced by various WENO3-NN models and the `ground truth'. Additionally, we estimate the adimensionalized Hausdorff distance~\cite{scipy_directed_hausdorff_webpage} between the contours at 300.3~K of each WENO3-NN model and the `ground truth', thus, providing a more quantitative assessment of the discrepancy in prediction between WENO3-NN models and the `ground truth'. In these cases, we use the same scheme (either WENO5-JS or a particular WENO3-NN model) for both the momentum and scalar transport equations.

\begin{figure}
     \centering
     \begin{subfigure}[b]{0.49\textwidth}
         \centering
         \includegraphics[width=\textwidth, clip, trim={0mm 0mm 110mm 0mm}]{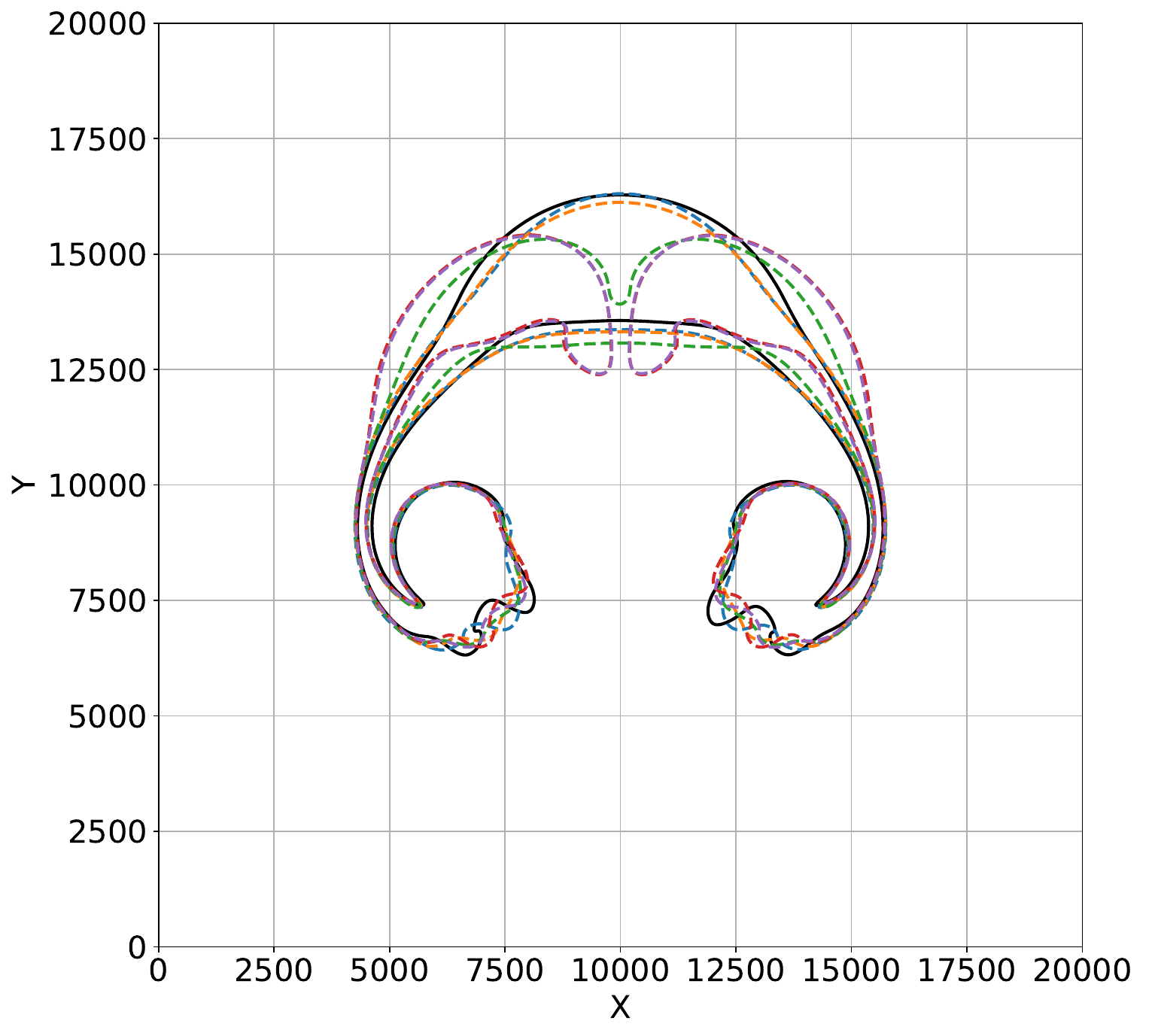}
         \caption{Grid: $512^2$}
         \label{fig:Smooth Bubble WENO-NN 512 300.3}
     \end{subfigure}
     \begin{subfigure}[b]{0.49\textwidth}
         \centering
         \includegraphics[width=\textwidth, clip, trim={0mm 0mm 110mm 0mm}]{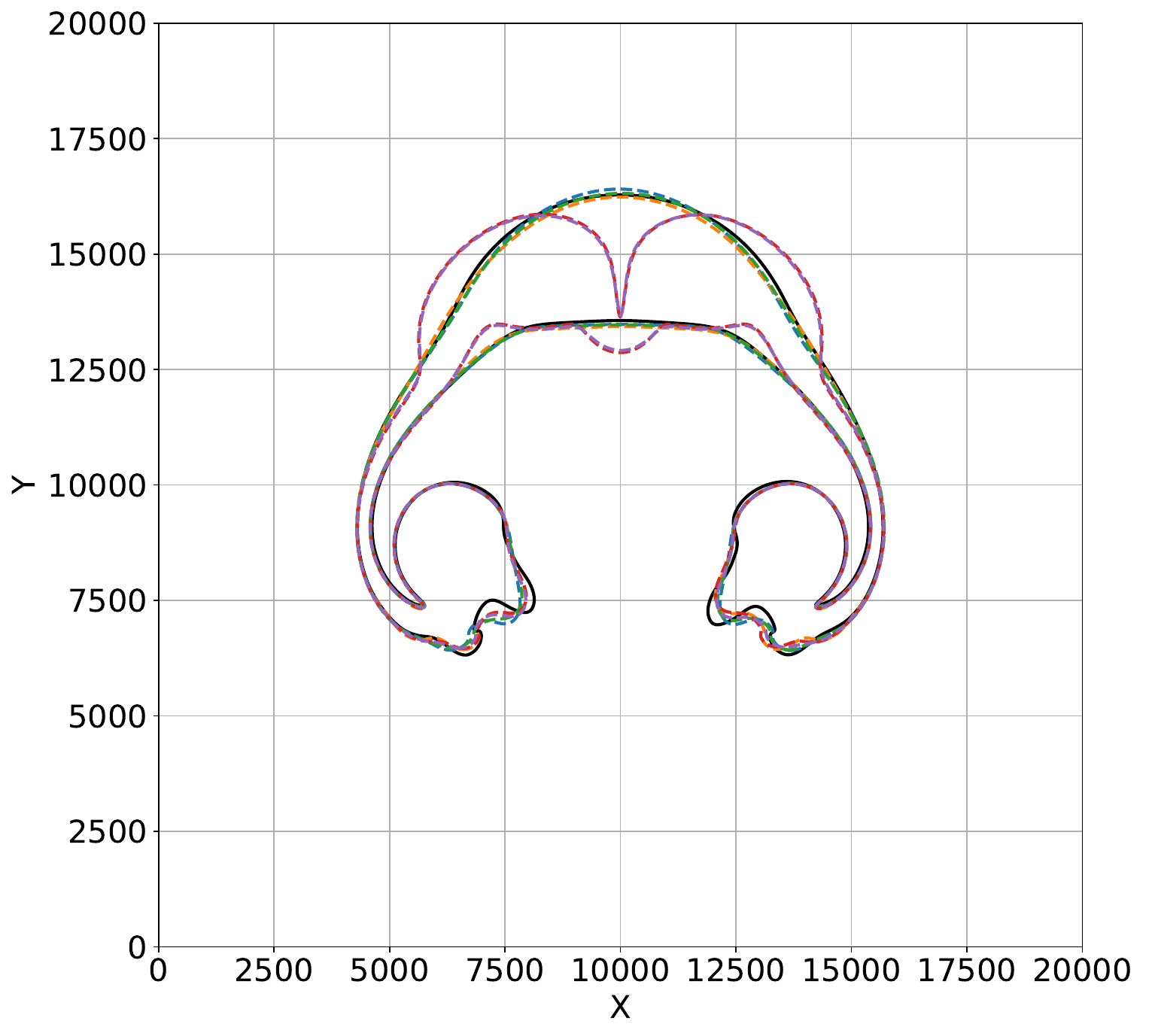}
         \caption{Grid: $1024^2$}
         \label{fig:Smooth Bubble WENO-NN 1024 300.3}
     \end{subfigure}
     \begin{subfigure}[b]{0.75\textwidth}
         \centering
         \includegraphics[width=\textwidth]{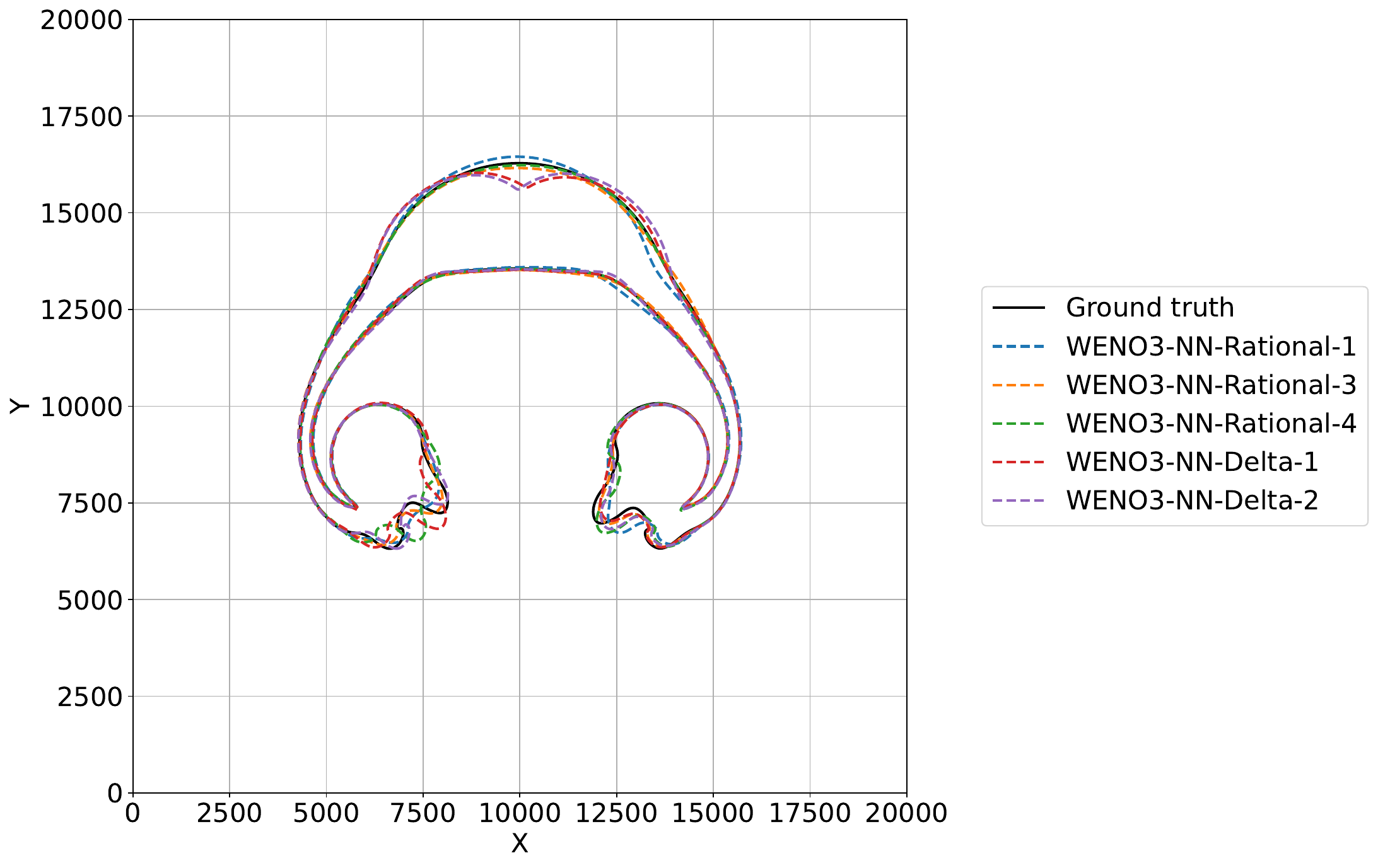}
         \caption{Grid: $2048^2$}
         \label{fig:Smooth Bubble WENO-NN 2048 300.3}
     \end{subfigure}
    \caption{Buoyant bubble with smooth initial condition: potential temperature contour at 300.3 K (WENO-NN schemes).}
    \label{fig:Smooth Bubble WENO-NN 300.3}
\end{figure}

\Cref{fig:Smooth Bubble WENO-NN 300.3} plots the contours of the potential temperature at $300.3$~K of the solutions generated using several WENO-NN models simulated with 3 different grid resolutions. From the behavior of the solutions we observe that some of the models are highly dissipative at smaller resolutions, even more than the baselines, as the deformation of the top of the bubble is more pronounced.
Such is the case for the WENO3-NN-Delta models.
Conversely, the models levering rational networks  demonstrate superior performance, by accurately representing the bubble's shape.

\begin{figure}[H]
    \centering
    \includegraphics[width=0.6\textwidth]{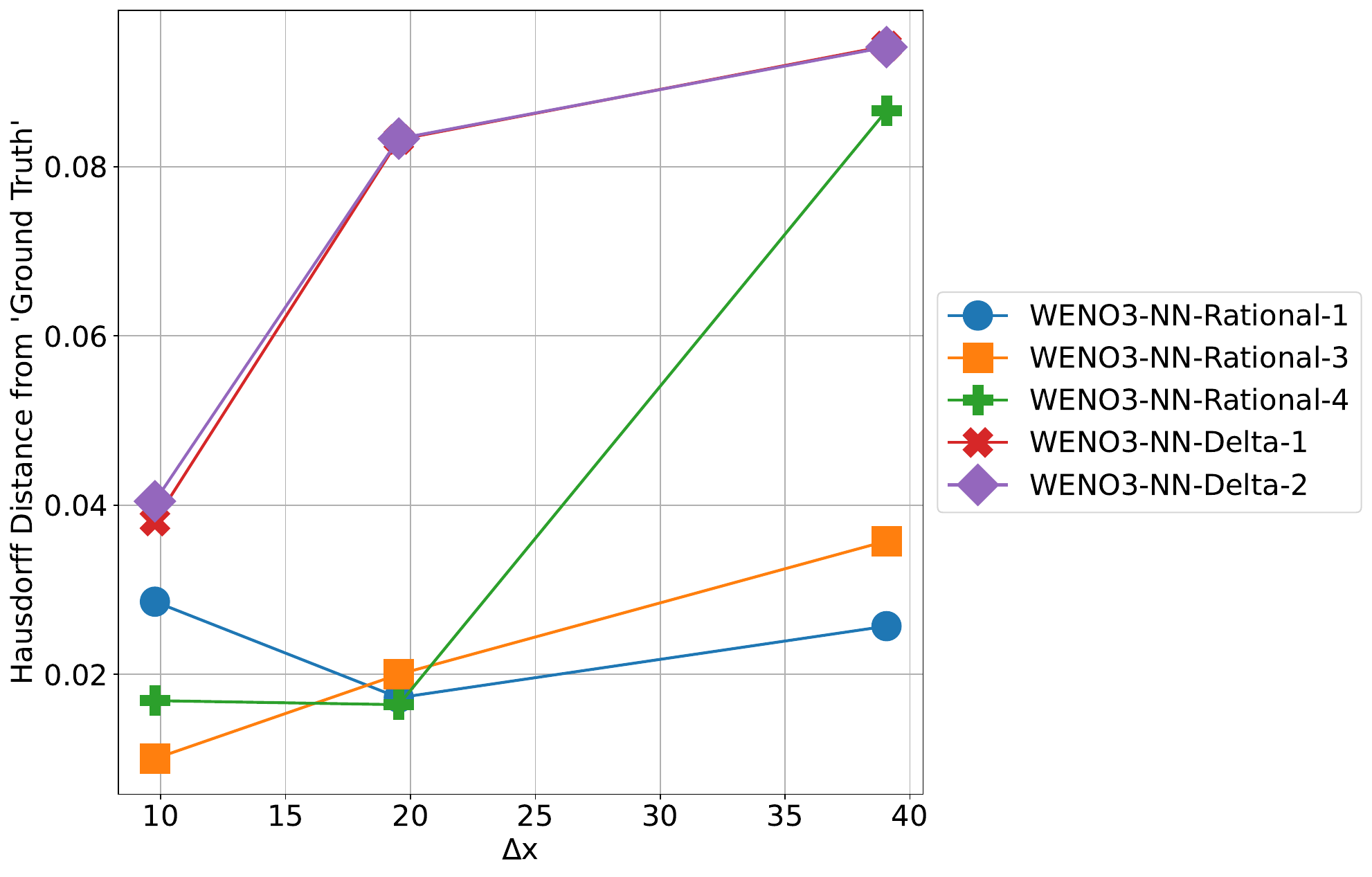}
    \caption{Buoyant bubble with smooth initial condition: Hausdorff distance between potential temperature contour at 300.3 K for WENO-NN schemes and `ground truth'.}
    \label{fig:Smooth Bubble WENO-NN Diff 300.3}
\end{figure}

\Cref{fig:Smooth Bubble WENO-NN Diff 300.3} present the average distance between the generated contours of individual WENO3-NN models and the `ground truth'. As expected, refining the grid resolution leads to a decrease in this distance for most WENO-NN models, indicating convergence towards the `ground truth'. Notably, WENO3-NN-Rational-1 exhibits a low distance even at a coarse resolution of $512^2$, with minimal change at higher resolutions. In summary, the WENO-NN-Rational models are able to resolve the bubble shape with a level of accuracy comparable to the ground truth using much coarser grids.


\paragraph{Computational Cost}
We compare the runtimes\footnote{We point out that even though we use implementations of the classical methods that are hand-tuned for TPUs~\cite{wang2022tensorflow} running in distributed environments, the timings shown in \cref{tab:timing_buble} using different hardware may differ.} for various schemes in a distributed environment using tensor processing units (TPUs) as accelerators.  
For this comparison, we simulated the buoyant bubble problem for a fixed grid size of $512^2$. We used $16$ \texttt{v5litepod} TPU chips with a $4 \times 4$ layout ($2$ hosts with $8$ chips each with high-speed interconnect). The same convection scheme was used for both momentum and scalar transport equations.
The wall-clock timings per time step are summarized in \cref{tab:timing_buble}.
\begin{table}[t]
    \centering
    \small
    \begin{tabular}{r | lllll}
        Model     &  QUICK & WENO3-JS & WENO5-JS & WENO3-NN-Delta-1 & WENO3-NN-Rational-1\\
        Time [ms] &   12   &  12    & 36     & 18            &   17 
    \end{tabular}
    \caption{Average wall clock time of one time step (in milliseconds) of the two-dimensional buoyant bubble simulation using different advection schemes.}
    \label{tab:timing_buble}
\end{table}
From \cref{tab:timing_buble} we observe that QUICK and WENO3-JS schemes are the fastest. Given the distributed nature of the computation, communication is needed in the form of halo exchange at each time step. The WENO3-JS and QUICK schemes require a halo width of 2 grid points whereas the WENO5-JS scheme uses a halo width of 3. WENO5-JS scheme takes thrice the runtime as WENO3-JS due to higher data transfer between TPU cores and the wider stencil, which requires more delocalized memory access. Both the WENO3-NN schemes are slower than WENO3-JS due to increased computational cost in the estimation of WENO weights. However, they are twice faster than WENO5-JS and have a similar accuracy as WENO5-JS in several simulations presented here. Further, the WENO3-NN-R-1 scheme using the rational network developed in this research is slightly faster than the previous WENO3-NN-D-1 scheme~\cite{bezgin2022weno3} as the rational network has four times fewer training parameters.

\subsubsection{Buoyant bubble with sharp initial condition}
\begin{figure}[H]
    \centering
    \includegraphics[width=0.5\textwidth]{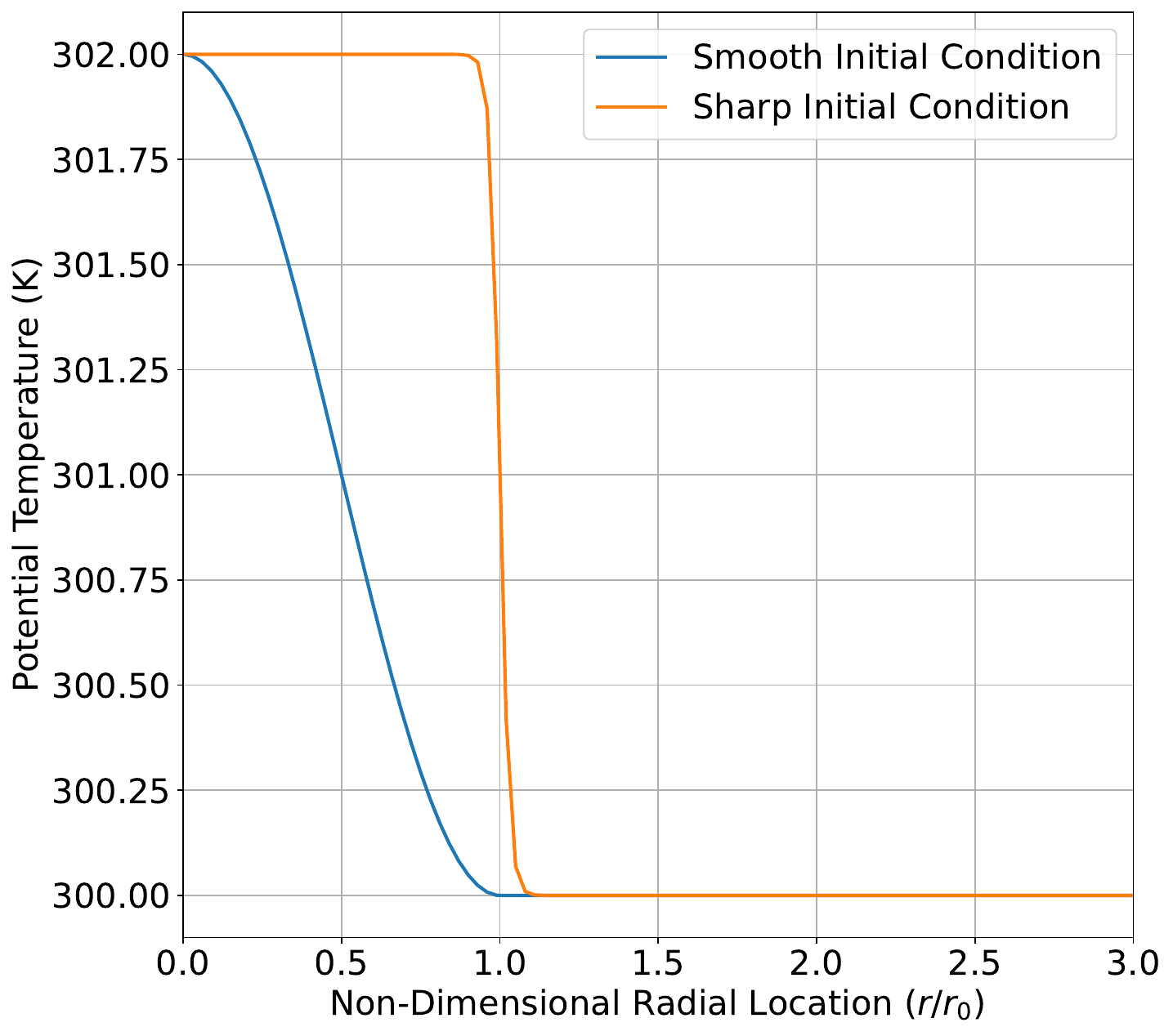}
    \caption{Radial plot of the initial conditions for buoyant bubble.}
    \label{fig:Bubble ICs}
\end{figure}

In \cref{sec:smooth bubble}, we initialized the bubble with a smooth condition given in \cref{eq:smooth bubble}. In order to analyze the performance of the WENO3-NN schemes on sharp fronts, we modified the initial condition as follows: 
\begin{equation}
    T = 300 + \tanh((r+r_0)/ \Delta r) - \tanh((r-r_0)/ \Delta r).
\end{equation}
where, $\Delta r = \sqrt{\Delta x^2 + \Delta y^2 }$, $r=\sqrt{(x-x_c)^2 + (y-y_c)^2}$, $(x_c,y_c) = (10000,2000)$ and $r_0=2000$~m. \Cref{fig:Bubble ICs} plots both the initial conditions for comparison.
The sharp condition is a function of the grid spacing $\Delta x$ and $\Delta y$ such that with refinement, the temperature drops quickly from 302~K to the ambient value of 300~K. All other details are similar to the smooth bubble case described in \cref{sec:smooth bubble}.

\begin{figure}[H]
     \centering
     \begin{subfigure}[b]{0.49\textwidth}
         \centering
         \includegraphics[width=\textwidth, clip, trim={0mm 0mm 150mm 0mm}]{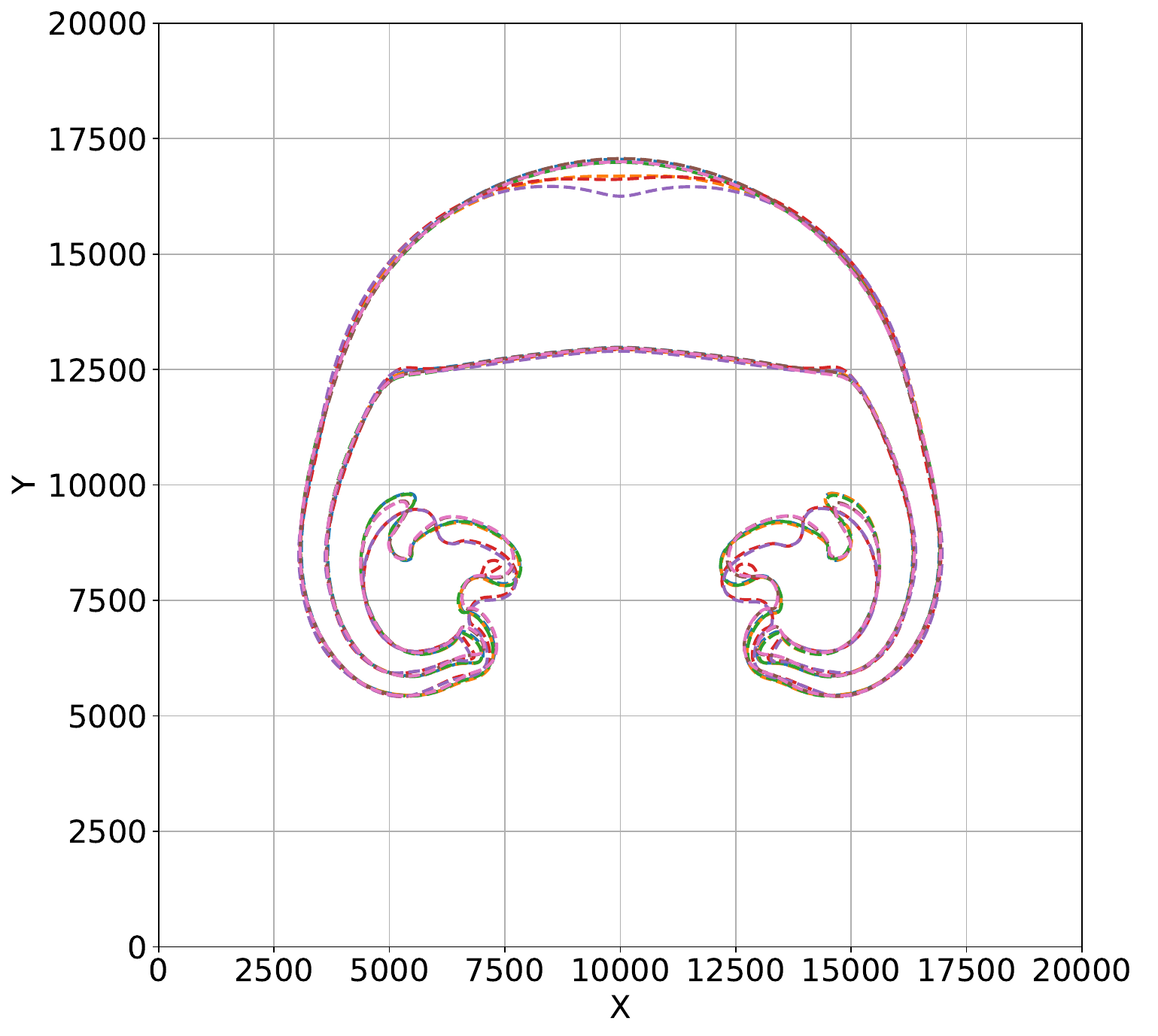}
         \caption{Grid: $512^2$.}
         \label{fig:Sharp Bubble conventional 512 300.5}
     \end{subfigure}
     \begin{subfigure}[b]{0.49\textwidth}
         \centering
         \includegraphics[width=\textwidth, clip, trim={0mm 0mm 150mm 0mm}]{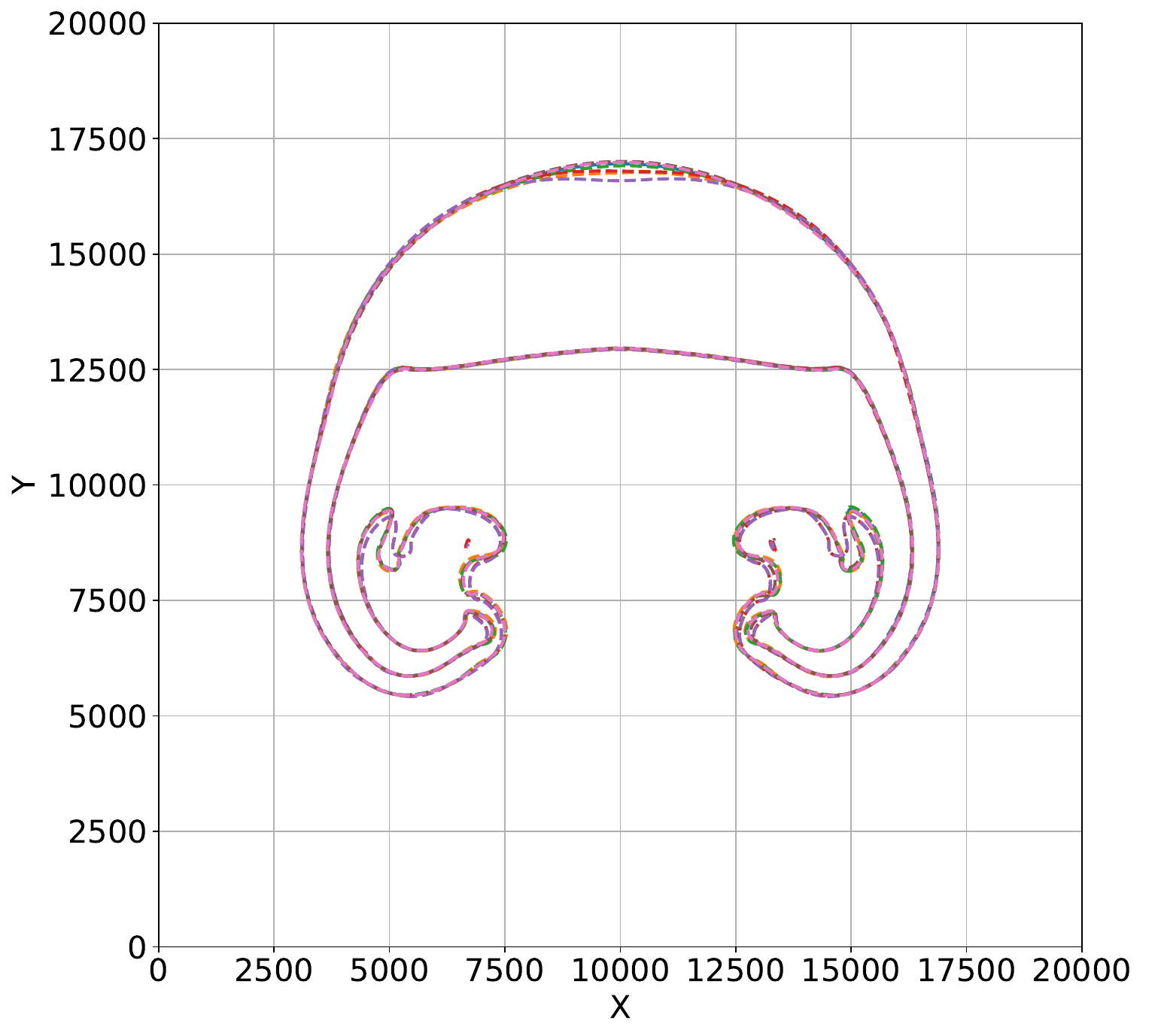}
         \caption{Grid: $1024^2$.}
         \label{fig:Sharp Bubble conventional 1024 300.5}
     \end{subfigure}
     \begin{subfigure}[b]{0.75\textwidth}
         \centering
         \includegraphics[width=\textwidth]{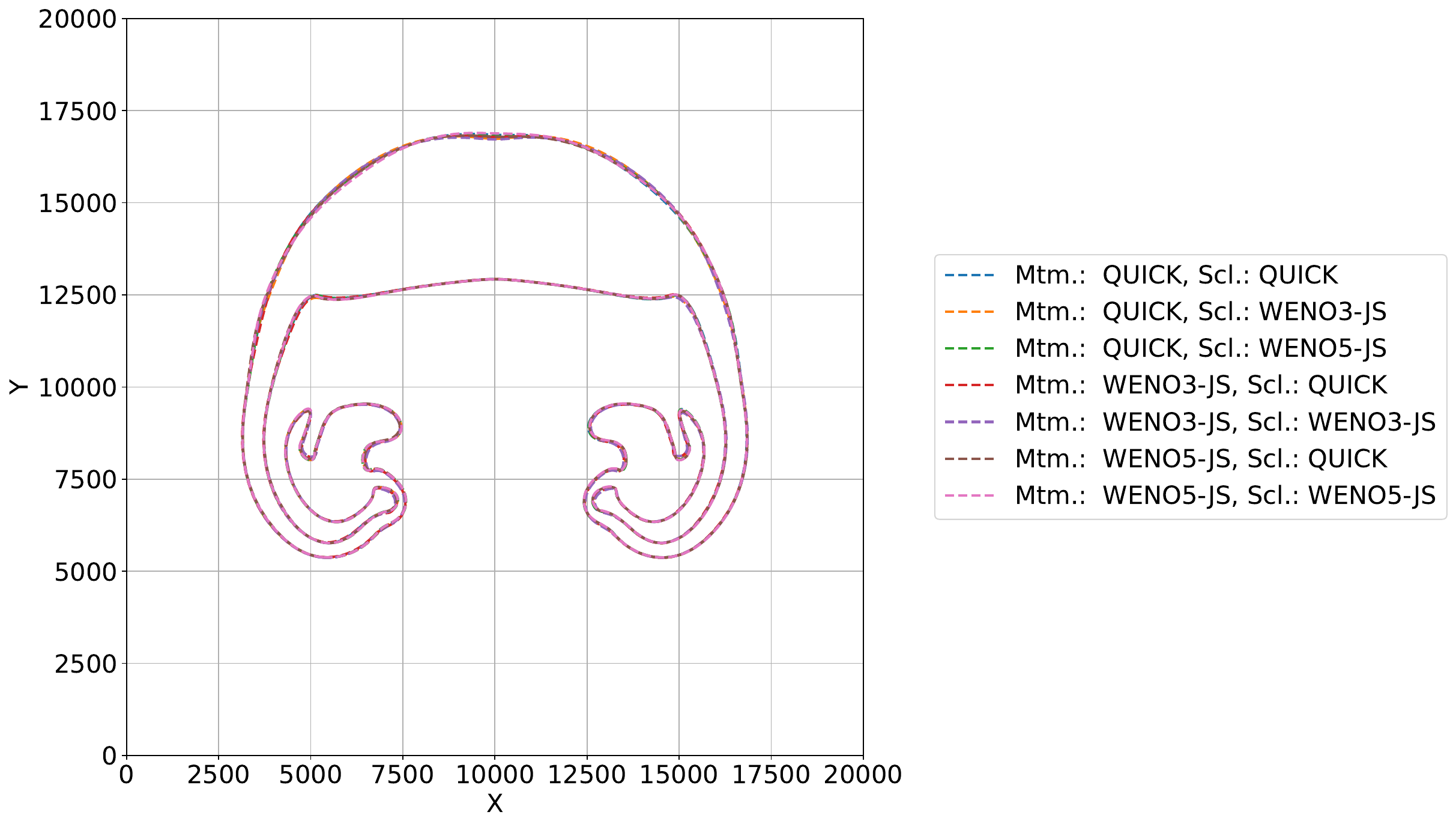}
         \caption{Grid: $2048^2$.}
         \label{fig:Sharp Bubble conventional 2048 300.5}
     \end{subfigure}
    \caption{Potential temperature contour at 300.5 K of the buoyant bubble simulation with sharp initial conditions solved with conventional schemes.}
    \label{fig:Sharp Bubble conventional 300.5}
\end{figure}

\Cref{fig:Sharp Bubble conventional 300.5} plots the contour lines for the value of 300.5 K for the conventional schemes and grid resolution of $512^2$, $1024^2$ and $2048^2$. In contrast to the smooth bubble, we do not see much dissipation with WENO3-JS in this case. For the finest grid of $2048^2$, the contours of all schemes overlap in \cref{fig:Sharp Bubble conventional 2048 300.5}. As before, we use the WENO5-JS simulation with the grid of $2048^2$ as the `ground truth' to compare with WENO3-NN.

\begin{figure}[H]
     \centering
     \begin{subfigure}[b]{0.49\textwidth}
         \centering
         \includegraphics[width=\textwidth, clip, trim={0mm 0mm 110mm 0mm}]{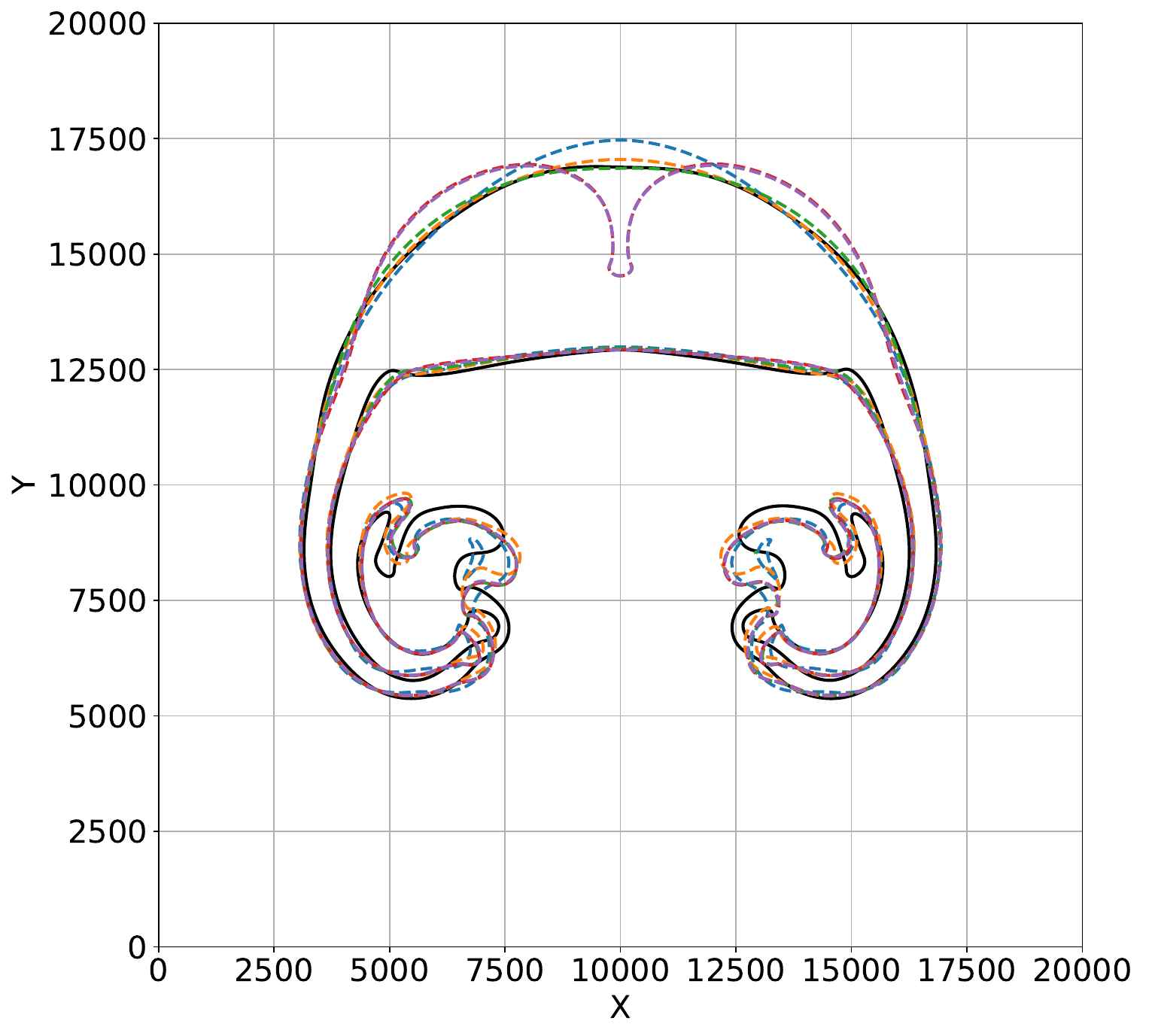}
         \caption{Grid: $512^2$.}
         \label{fig:Sharp Bubble WENO-NN 512 300.5}
     \end{subfigure}
     \begin{subfigure}[b]{0.49\textwidth}
         \centering
         \includegraphics[width=\textwidth, clip, trim={0mm 0mm 110mm 0mm}]{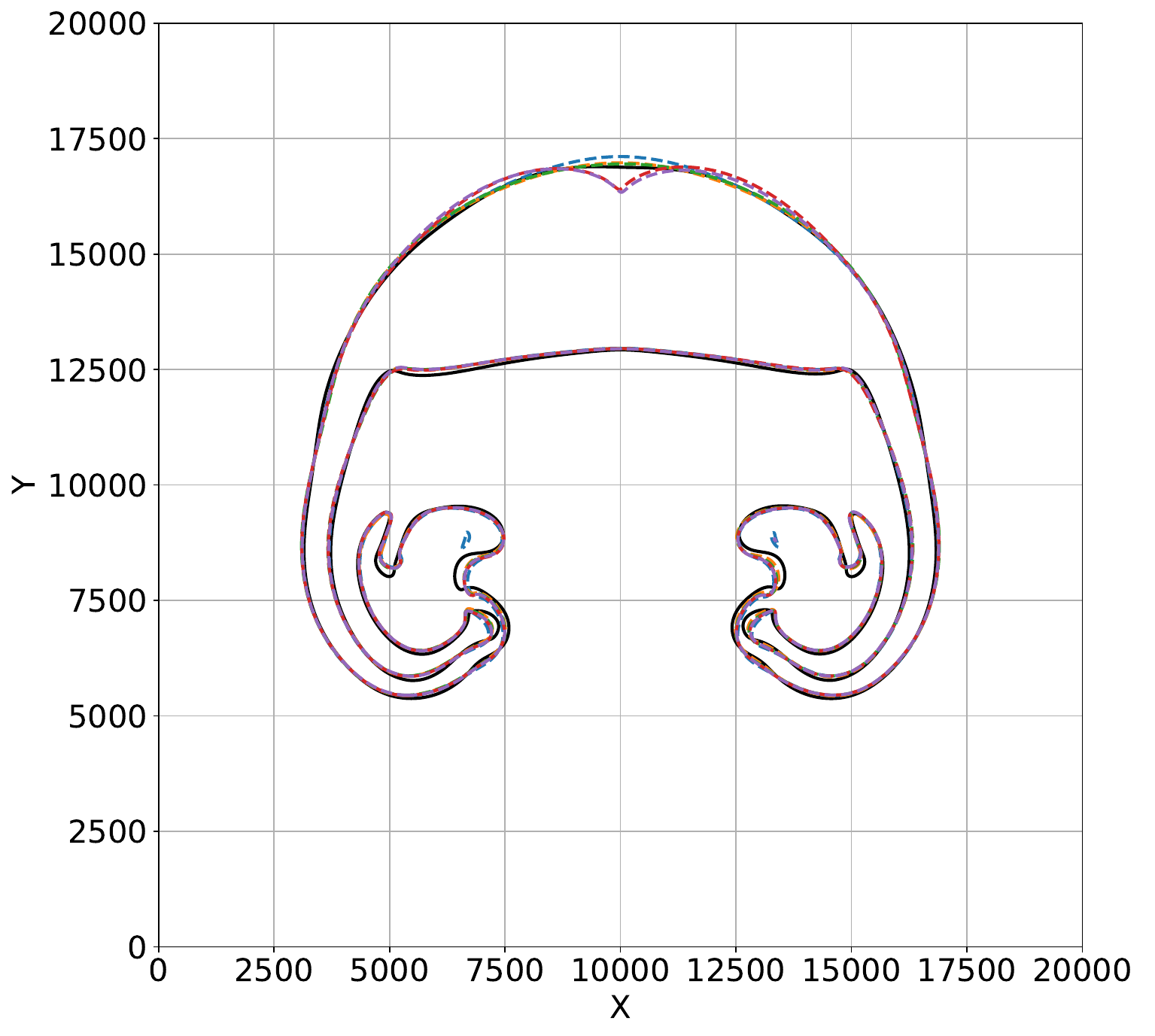}
         \caption{Grid: $1024^2$}
         \label{fig:Sharp Bubble WENO-NN 1024 300.5}
     \end{subfigure}
     \begin{subfigure}[b]{0.75\textwidth}
         \centering
         \includegraphics[width=\textwidth]{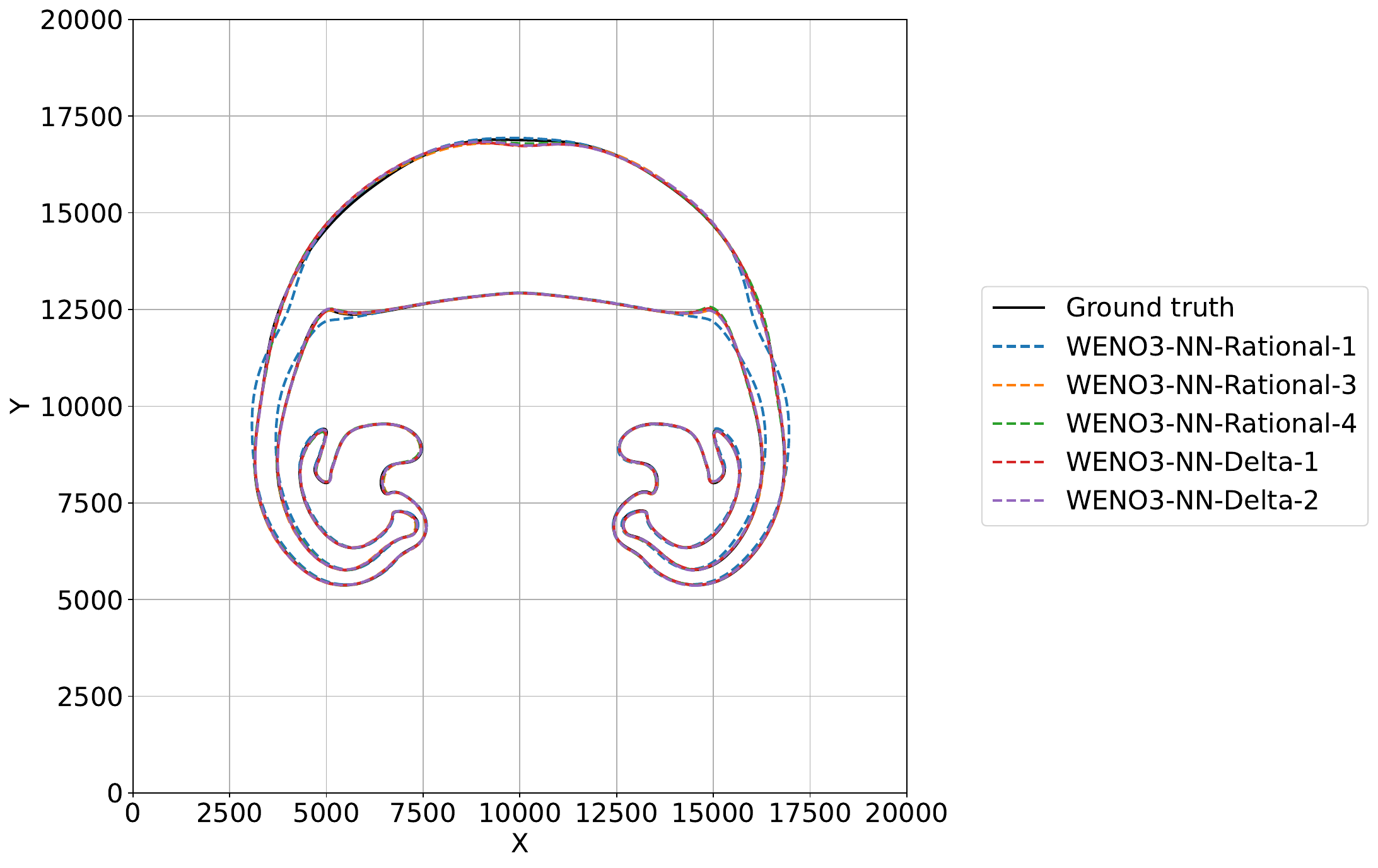}
         \caption{Grid: $2048^2$.}
         \label{fig:Sharp Bubble WENO-NN 2048 300.5}
     \end{subfigure}
    \caption{Potential temperature contour at 300.5 K of the buoyant bubble simulation with sharp initial conditions solved with WENO-NN schemes.}
    \label{fig:Sharp Bubble WENO-NN 300.5}
\end{figure}

\begin{figure}[H]
    \centering
    \includegraphics[width=0.6\textwidth]{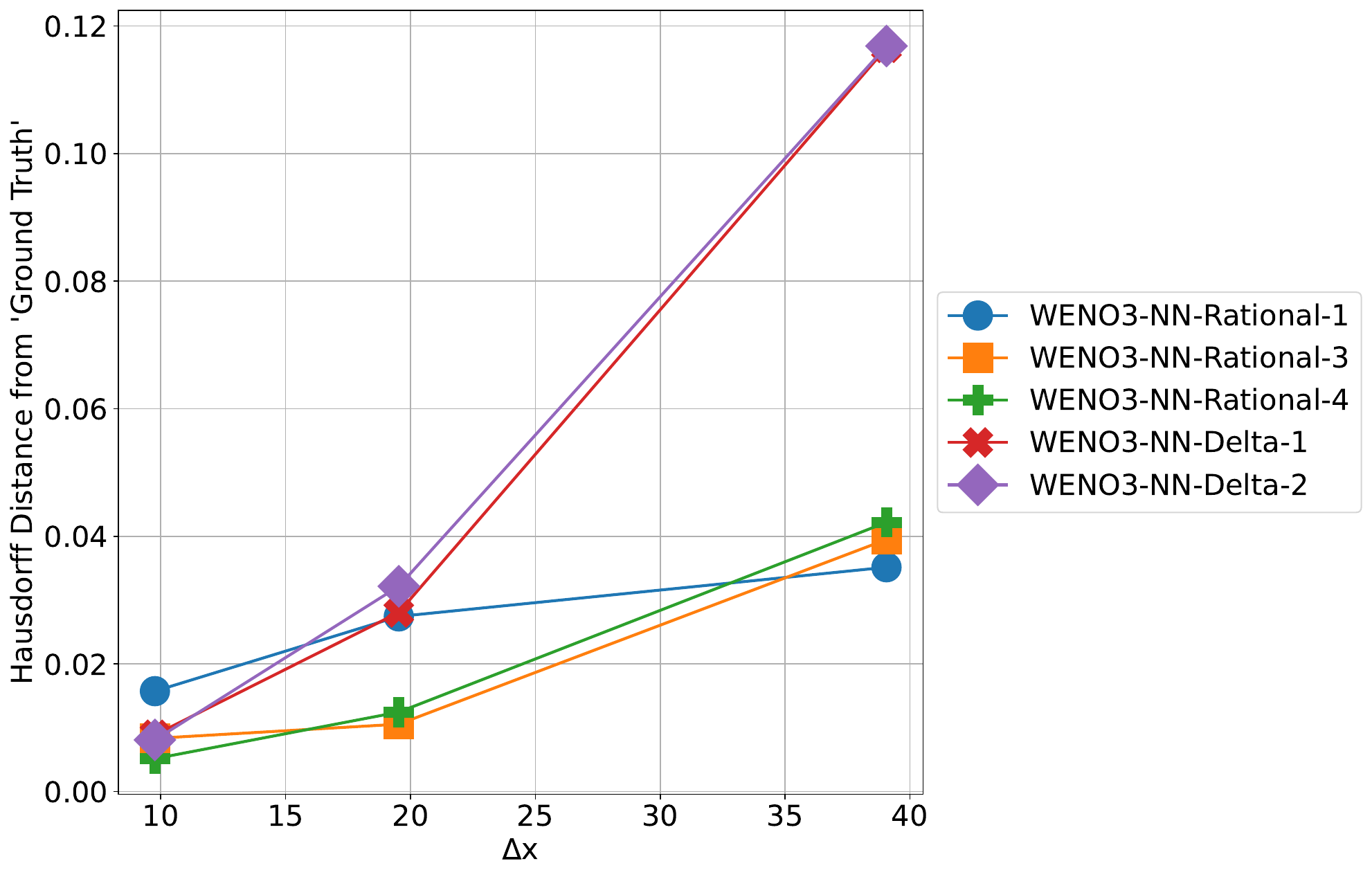}
    \caption{Buoyant bubble with sharp initial condition: Hausdorff distance between potential temperature contour at 300.5 K for WENO-NN schemes and `ground truth'.}
    \label{fig:Sharp Bubble WENO-NN Diff 300.5}
\end{figure}

\Cref{fig:Sharp Bubble WENO-NN 300.5} plots the contour using various WENO-NN models. Similar to the smooth bubble case, we see that some of the networks introduce higher dissipation. 
The networks selected on the basis of the Sine-step function (\cref{tab:net_hyper_params}) have a much better agreement with the fine grid WENO5-JS simulation.
Similar to \cref{sec:smooth bubble}, we tabulate and plot the non-dimensional Hausdorff distance between the contours for the WENO3-NN methods and the `ground truth'.
We also observe from \cref{fig:Sharp Bubble WENO-NN Diff 300.5} that WENO3-NN schemes converge to the `ground truth' as the grid is resolved.


\subsubsection{Kelvin--Helmholtz instability}
\begin{figure}[H]
     \centering
     \begin{subfigure}[b]{\textwidth}
         \centering
         \includegraphics[width=\textwidth]{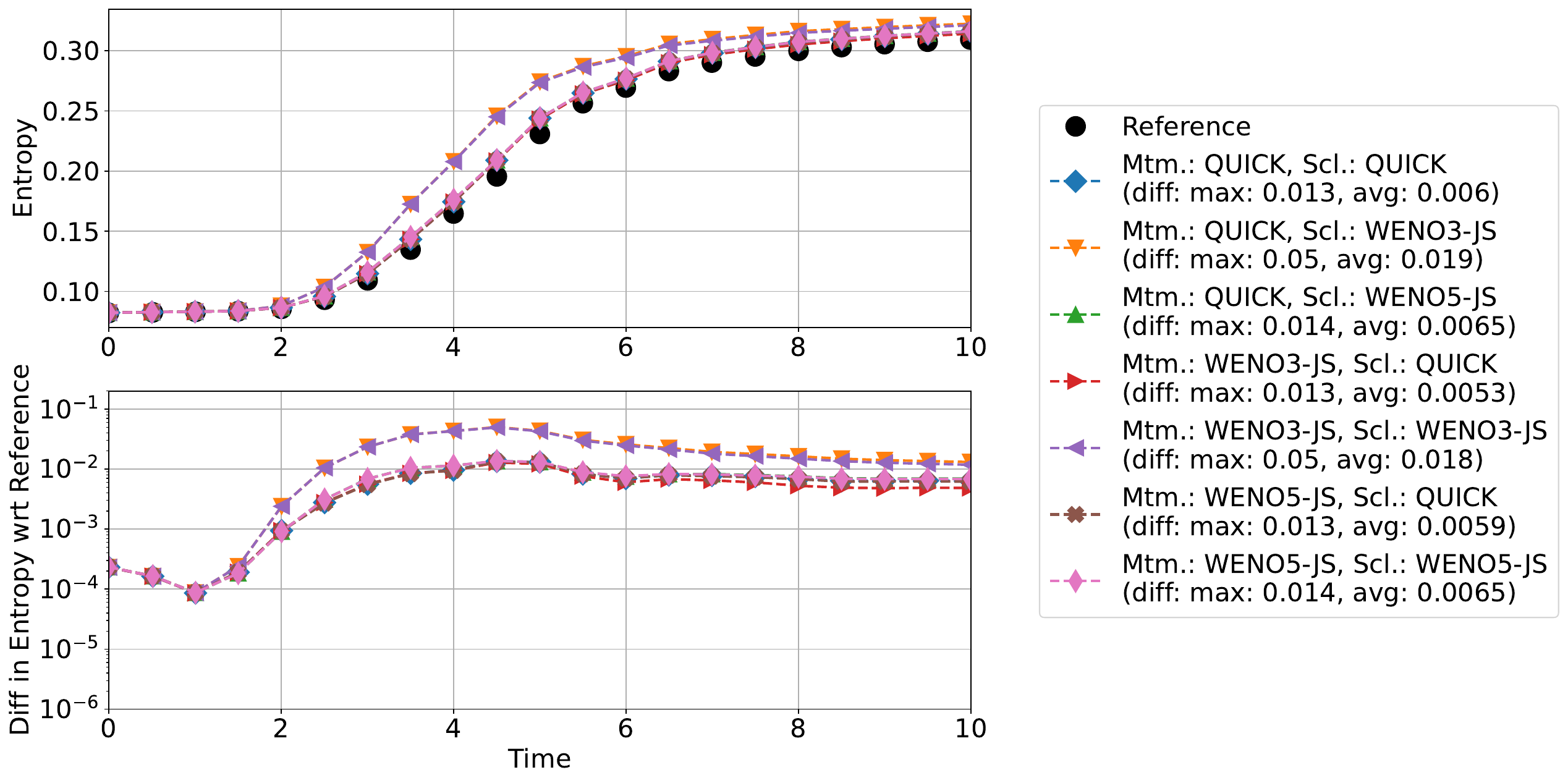}
         \caption{Conventional schemes}
         \label{fig:KH_Instability Entropy Re: 1E5 Grid: 512 X 1024 conv}
     \end{subfigure}
     \begin{subfigure}[b]{\textwidth}
         \centering
         \includegraphics[width=\textwidth]{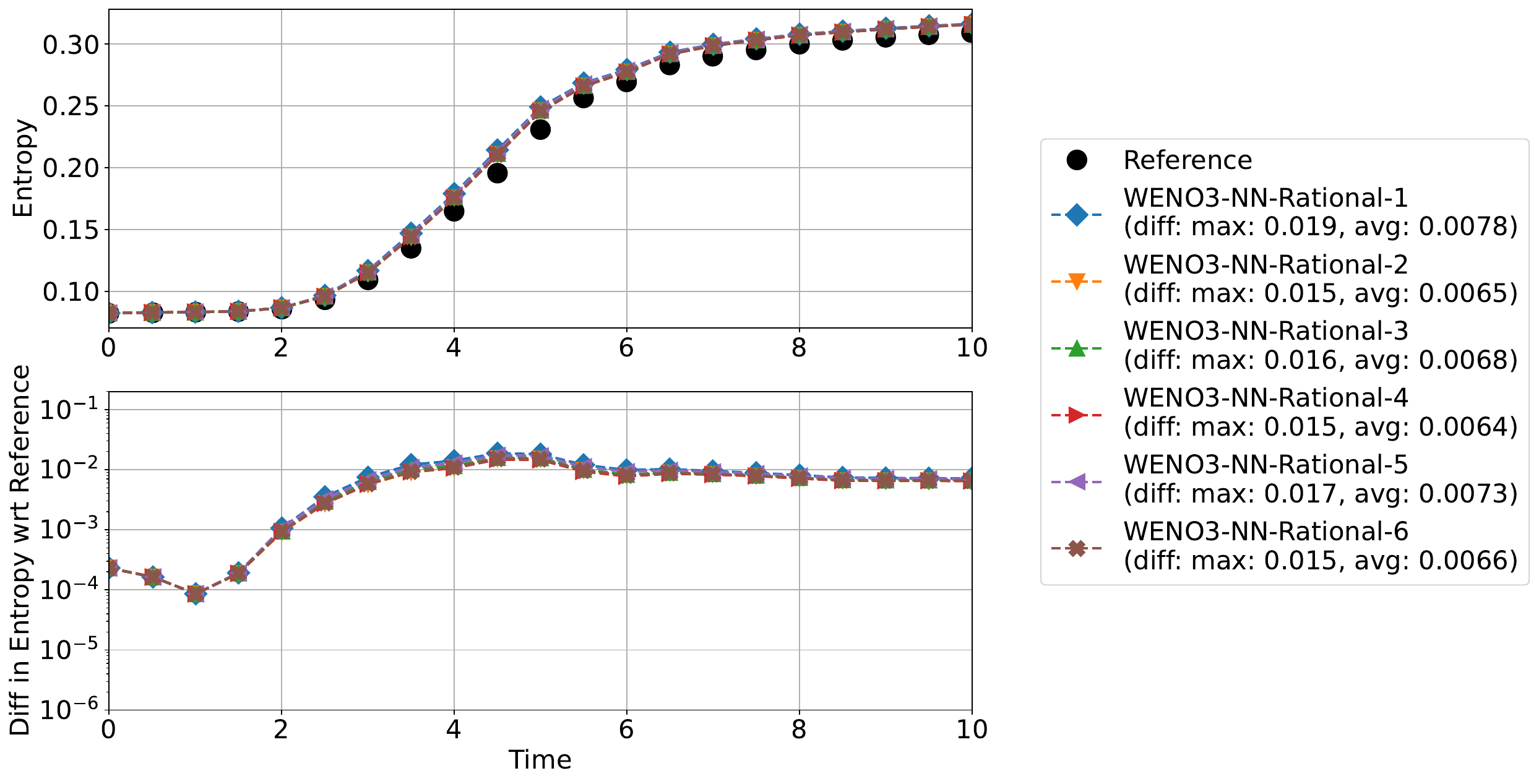}
         \caption{WENO3-NN schemes}
         \label{fig:KH_Instability Entropy Re: 1E5 Grid: 512 X 1024 weno-nn}
     \end{subfigure}
    \caption{Kelvin--Helmholtz instability: volume-integrated entropy (Re: $10^5$, Grid: $512 \times 1024$, reference:~\citet{lecoanet2016validated}).}
    \label{fig:KH_Instability Entropy Re: 1E5 Grid: 512 X 1024}
\end{figure}

\begin{figure}[H]
     \centering
     \begin{subfigure}[b]{\textwidth}
         \centering
         \includegraphics[width=\textwidth]{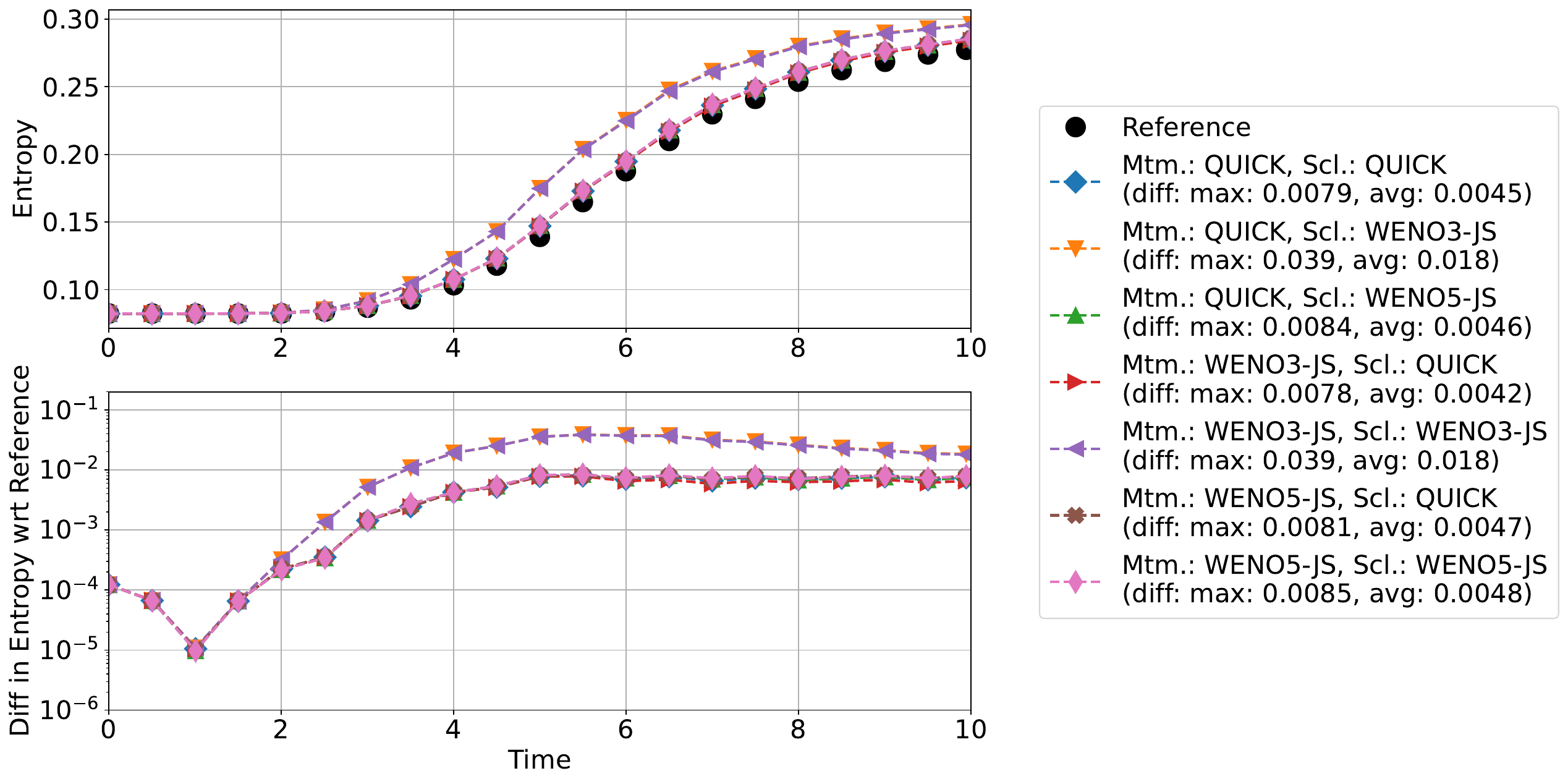}
         \caption{Conventional schemes}
         \label{fig:KH_Instability Entropy Re: 1E6 Grid: 2048 X 4096 conv}
     \end{subfigure}
     \begin{subfigure}[b]{\textwidth}
         \centering
         \includegraphics[width=\textwidth]{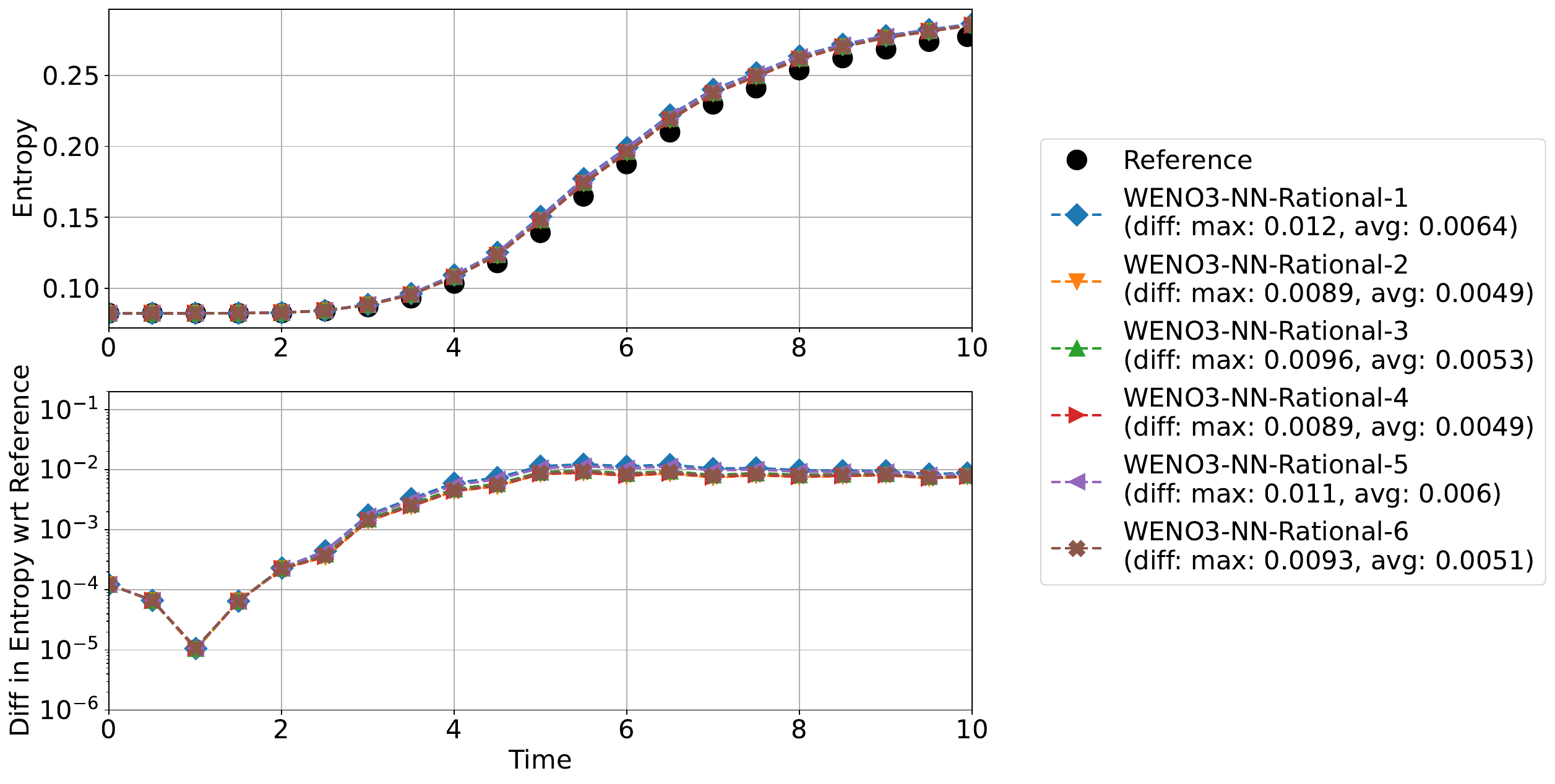}
         \caption{WENO3-NN schemes}
         \label{fig:KH_Instability Entropy Re: 1E6 Grid: 2048 X 4096 weno-nn}
     \end{subfigure}
    \caption{Kelvin--Helmholtz instability: volume-integrated entropy (Re: $10^6$, Grid: $2048 \times 4096$, reference:~\citet{lecoanet2016validated}).}
    \label{fig:KH_Instability Entropy Re: 1E6 Grid: 2048 X 4096}
\end{figure}

In this section, we compare our results with the well known benchmark of \citet{lecoanet2016validated}, consisting of the Kelvin--Helmholtz instability on a rectangular domain of size $1 \times 2$.
The following equations, along with continuity \cref{eq:continuity}, are solved with periodic boundary conditions:
\begin{equation}
    \frac{\partial (\rho \bm{u})}{\partial t} + \nabla \cdot (\rho \bm{u} \bm{u}) = -\nabla p + \nabla \cdot \tau,
\end{equation}
where, $\rho$: density, $\bm{u}$: velocity vector, $\rho$: density of fluid, $\tau$: stress tensor and $p$: pressure. A passive scalar referred to as `dye' is advected following
\begin{equation}
    \frac{\partial (\rho c)}{\partial t} + \nabla \cdot (\rho c \bm{u}) = \nabla \cdot (\rho \nu_{dye} \nabla c),
\end{equation}
where, $c$ is local fraction of dye particles and $\nu_{dye}$ is the viscosity of the dye. 
As outlined in \citet{lecoanet2016validated}, the initial conditions for  $(u,v)$ and $(c)$ are given by
\begin{subequations}
    \begin{equation}
        u(x,y,t=0) = u_0 \left[ \tanh\left(\frac{y-y_1}{a}\right) - \tanh\left(\frac{y-y_2}{a}\right) - 1 \right],
    \end{equation}    
    \begin{equation}
        v(x,y,t=0) = v_0 \sin(2 \pi x) \left[ \exp\left(-\frac{(y-y_1)^2}{\sigma^2}\right) + \exp\left(-\frac{(y-y_2)^2}{\sigma^2}\right) \right],
    \end{equation}
    \begin{equation}
        c(x,y,t=0) = \frac{1}{2} \left[ \tanh\left(\frac{y-y_2}{a}\right) - \tanh\left(\frac{y-y_1}{a}\right) + 2 \right],
    \end{equation}
\end{subequations}
where, $a=0.05$, $\sigma=0.2$, $u_0=1$, $v_0=0.01$, $y_1=0.5$ and $y_2$=1.5. Density and pressure are initialized to 1.0 and 10.0, respectively. The entropy per unit mass is $s=-c \log(c)$, together with its volume integral:
\begin{equation}
    S(t) = \int \rho \,\, s(x, y, t) \, dV.
\end{equation}
The Reynolds number is defined as $\mathrm{Re}=2 u_0 L / \nu$, where $\nu$ is the kinematic viscosity of the fluid. The diffusivity of the dye is set as $\nu_{dye} =\nu$ which corresponds to a Prandtl number of unity. Figures \ref{fig:KH_Instability Entropy Re: 1E5 Grid: 512 X 1024} and \ref{fig:KH_Instability Entropy Re: 1E6 Grid: 2048 X 4096} depict the temporal variation of the volume-integrated entropy for the different schemes together with the reference~\cite{lecoanet2016validated}.

We first present the temporal variation of volume integral of entropy for Reynolds number of $10^5$ for which we used a grid of dimension $512 \times 1024$. Similar to the buoyant bubble problem, we used various combinations of conventional schemes such as QUICK, WENO3-JS and WENO5-JS for the momentum and scalar transport equations.
When compared with \citet{lecoanet2016validated}, \cref{fig:KH_Instability Entropy Re: 1E5 Grid: 512 X 1024 conv} and \cref{fig:KH_Instability Entropy Re: 1E5 Grid: 512 X 1024 weno-nn} show that the WENO3-NN models significantly outperform WENO3-JS in accuracy.

For the simulation with  Reynolds number of $10^6$ we used a refined grid of $2048 \times 4096$. Comparing figures~\ref{fig:KH_Instability Entropy Re: 1E6 Grid: 2048 X 4096 conv} and \ref{fig:KH_Instability Entropy Re: 1E6 Grid: 2048 X 4096 weno-nn} shows that the configurations with WENO3-NN also outperforms WENO3-JS, while producing predictions similar to WENO5-JS in terms of accuracy.

\subsection{Three-dimensional problems} 

To demonstrate the performance of the WENO3-NN method in a three-dimensional setup, we performed a simulation of the low-cloud configuration based on the first nocturnal research of the Dynamics and Chemistry of Marine Stratocumulus (DYCOMS-II) field study. This configuration has a quasi two layer structure with a stable cloud layer as a result of the balance between turbulent mixing and radiative heating. Any excessive dissipation or spurious oscillations will break this balance and thus, mispredict the turbulent statistics or even destroy the cloud structure. The sensitivity of the DYCOMS-II configuration to sub-grid scale mixing makes it desirable for assessing the dissipative behavior of numerical schemes.

We simulate the anelastic equations for moist air, understood to be an ideal admixture of dry air, water vapor, and any condensed water. 
The governing equations are similar to the ones used in the thermal bubble simulation (\cref{eq:momentum,eq:theta_l_transport,eq:q_t_transport}), the continuity and momentum equations coupled with the transport equations for liquid-ice potential temperature $\theta_\text{li}$ and total water specific humidity $q_t$. The main difference from those equations is that the time-dependent density is replaced with a static reference density. Additionally, the scalar transport equations contain additional sources for the subsidence and radiative transfer. As the condensed water is always assumed to be in local thermodynamic equilibrium, the two thermodynamic variables are sufficient to fully represent the thermodynamic processes involved in cloud dynamics. They are initialized with the following profiles along the vertical direction, which is aligned with the $z$ axis:
\begin{equation}
    \theta_\text{li} =
    \begin{cases}
      289\text{ K} &z\leq z_i \\
      297.5 + \left(z - z_i\right)^{\sfrac{1}{3}}\text{ K} &z>z_i,
    \end{cases}
\end{equation}
\begin{equation}
    q_t =
    \begin{cases}
      9\text{ g/kg} &z\leq z_i \\
      1.5\text{ g/kg} &z>z_i,
    \end{cases}
\end{equation}
where $z_i$ is the initial height of the cloud top. A geophysical wind with free stream velocity $\boldsymbol{U} = \left(U, V\right)=(7,-5.5)$~m/s is enforced with a Coriolis force that takes the form
\begin{equation}
    f_g=2\Omega\sin\psi\hat{k}\times\left(\boldsymbol{u}-\boldsymbol{U}\right),
\end{equation}
where $\Omega=7.2921\times10^{-5}$~rad/s is the rotation rate of the earth, $\psi=31.5^\circ$ is the latitude.

The computational domain is of size $4\times4\times1.5$~km$^3$, with $128\times128\times256$ mesh points in each dimension correspondingly. To represent the atmospheric boundary layer, the shear stress and heat flux at the bottom boundary of the domain is modeled with the Monin--Obukhov similarity theory. A free-slip wall is applied on the top boundary, with a Rayleigh damping layer applied on the top 10\% of the vertical domain to eliminate the gravitational wave. Periodic boundary conditions are applied in the horizontal directions. The simulation is performed for 5 physical hours with a time step size of $0.3$~s, which corresponds to a Courant (CFL) number of $0.3$. Flow field statistics are collected for the last hour of the simulation.

\begin{figure}[H]
    \centering
    \begin{subfigure}{0.49\textwidth}
        \includegraphics[width=\textwidth, clip, trim={0mm 0mm 108mm 0mm}]{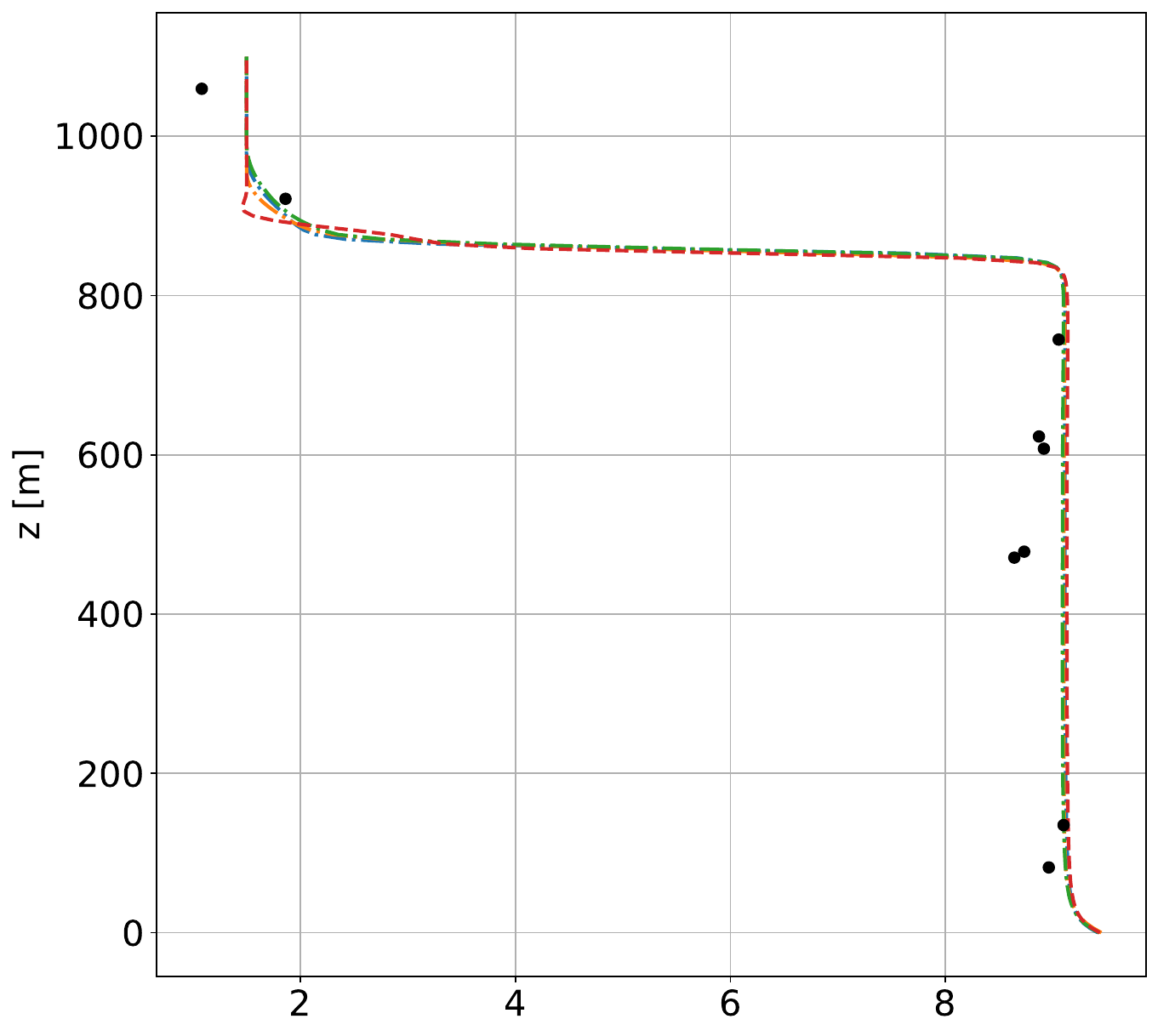}
        \caption{\label{fig:dycoms_q_t}$q_t\text{ [g/kg]}$}
    \end{subfigure}
    \begin{subfigure}{0.49\textwidth}
        \includegraphics[width=\textwidth, clip, trim={0mm 0mm 108mm 0mm}]{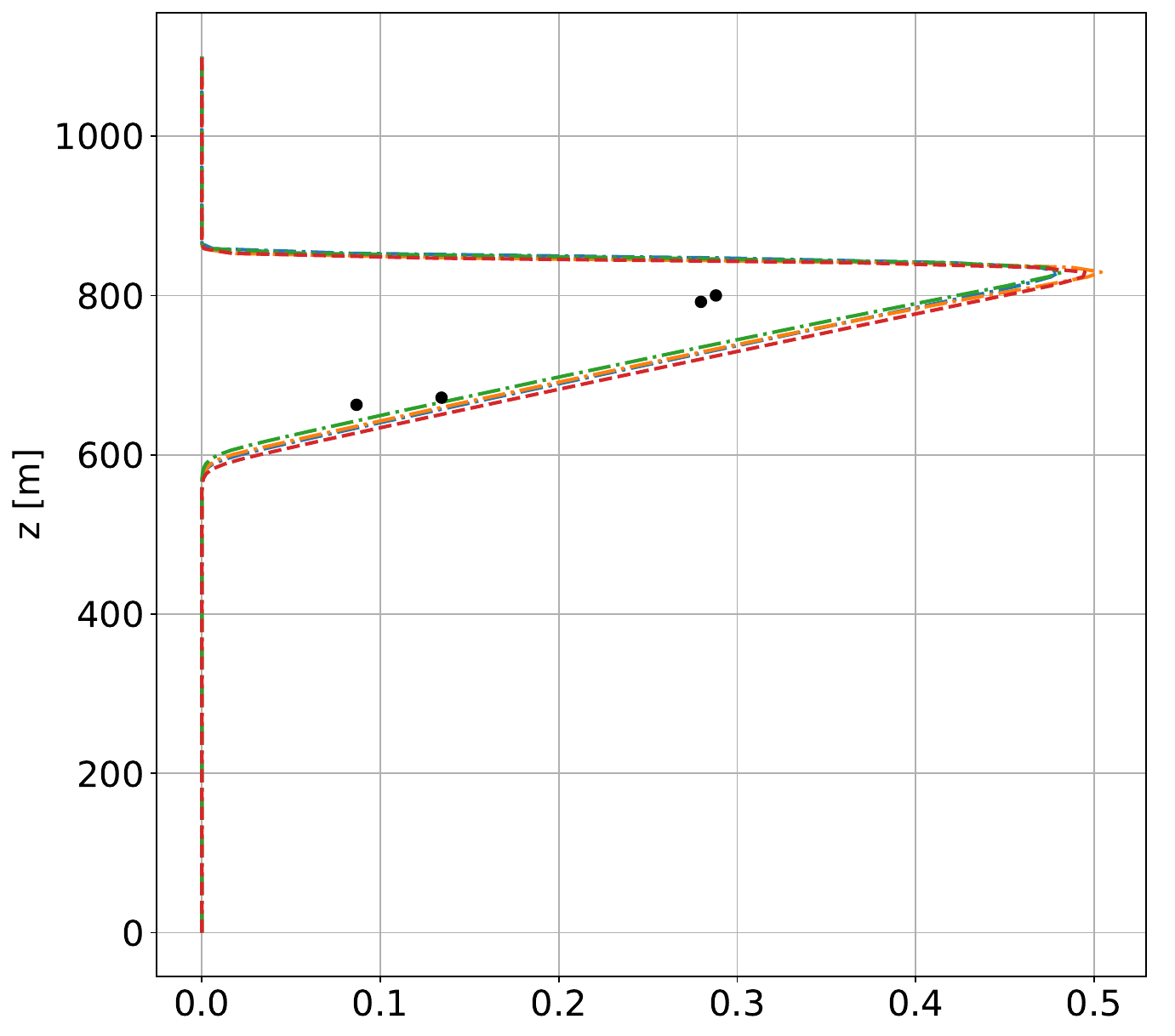}
        \caption{\label{fig:dycoms_q_l}$q_l\text{ [g/kg]}$}
    \end{subfigure}
    \begin{subfigure}{0.8\textwidth}
        \includegraphics[width=\textwidth]{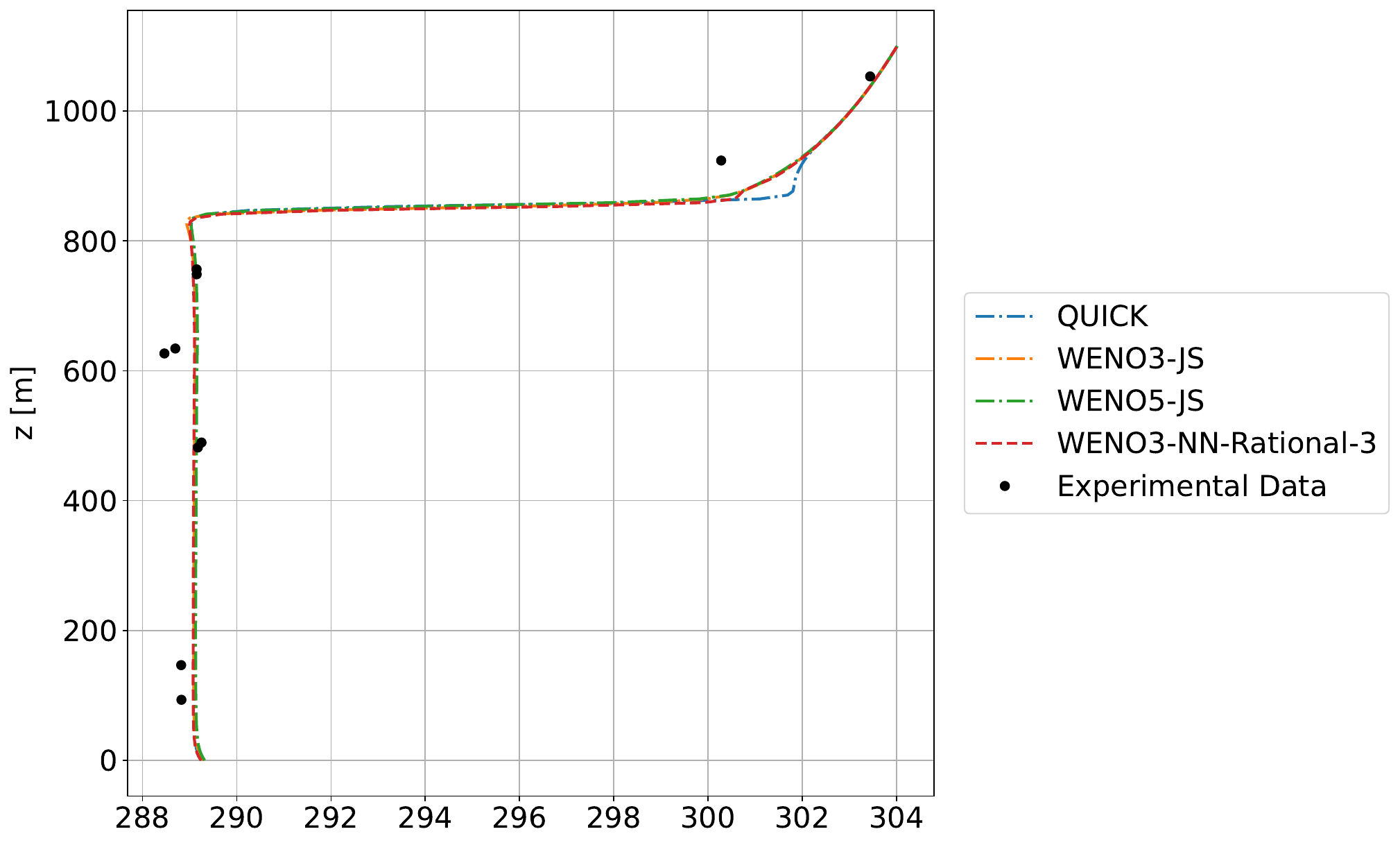}
        \caption{\label{fig:dycoms_thteta_li}$\theta_\text{li}\text{ [K]}$}
    \end{subfigure}
    \caption{\label{fig:dycoms_mean_profiles}Profiles of the mean humidity states and liquid-ice potential temperature. Solid dots indicate experimental observations.}
\end{figure}

\Cref{fig:dycoms_mean_profiles} shows the mean profiles of the humidity states and the liquid-ice potential temperature of simulation using four different interpolation schemes for the convective fluxes of both the momentum and scalar transport equations together with experimental observations. \Cref{fig:dycoms_mean_profiles} shows good agreement with the experimental results for all schemes. The two-layer structure for $\theta_\text{li}$ and $q_t$ are preserved with a well-mixed boundary layer below the cloud top. A stable cloud structure is represented in all cases, as shown in \cref{fig:dycoms_q_l}.
There is a spurious mixing layer above the inversion in \cref{fig:dycoms_thteta_li} with the QUICK scheme, which is attributed to numerical oscillations at the sharp interface as a result of the non-monotonicity of the scheme. This spurious structure is absent in simulations with WENO-based schemes, which indicates that the sub-grid scale mixing is well represented with these schemes.

\begin{figure}[H]
    \centering
    \begin{subfigure}{0.59\textwidth}
        \includegraphics[width=\textwidth]{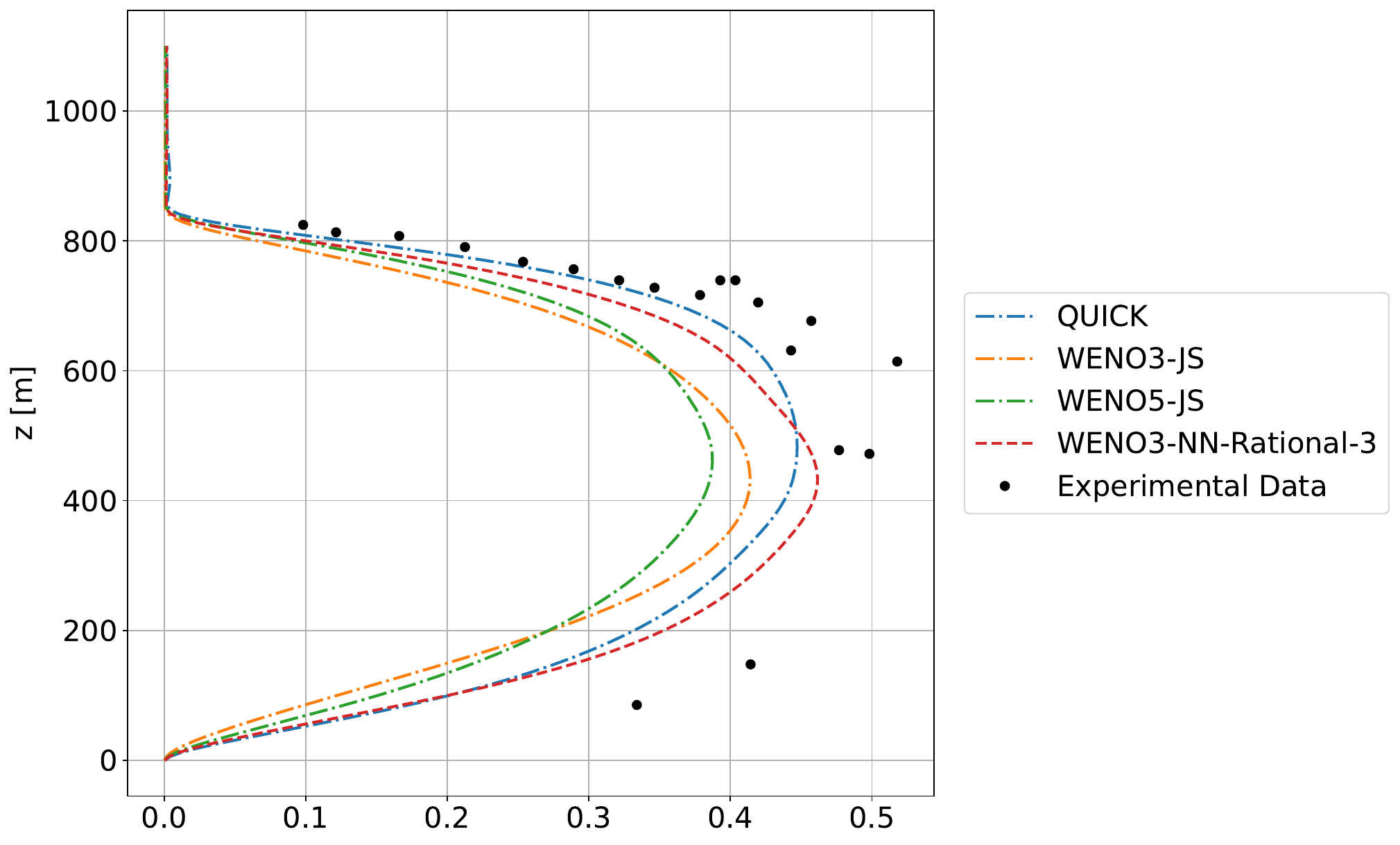}
        \caption{\label{fig:dycoms_ww}$\left<w^\prime\right>\ [\text{m}^2/\text{s}^2]$}
    \end{subfigure}
    \begin{subfigure}{0.4\textwidth}
        \includegraphics[width=\textwidth, clip, trim={0mm 0mm 110mm 0mm}]{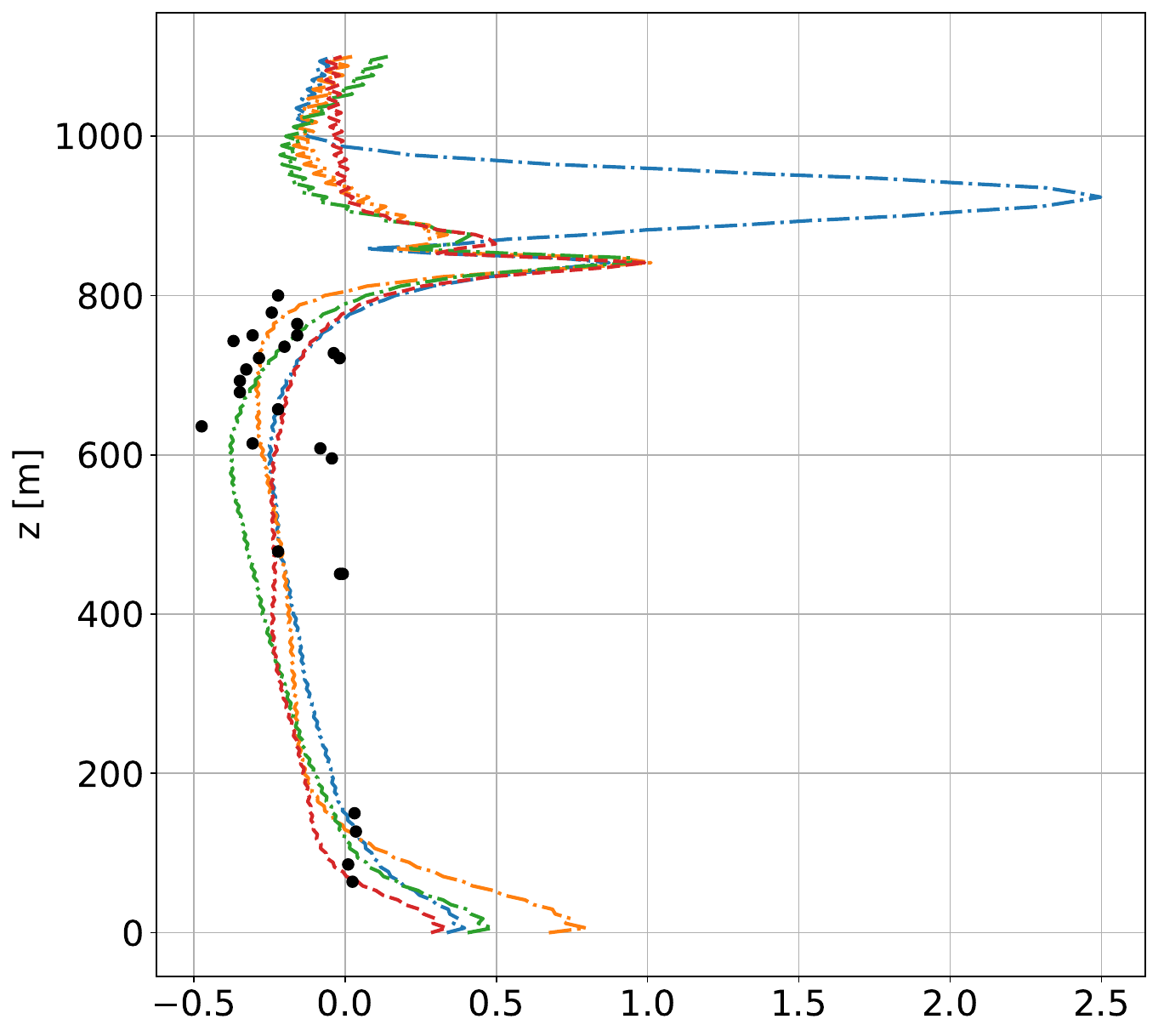}
        \caption{\label{fig:dycoms_s_w}Skewness}
    \end{subfigure}
    \caption{\label{fig:dycoms_turbulence}Profiles of the variance of the skewness of the vertical velocity. Solid dots indicate experimental observations.}
\end{figure}

To further assess the turbulence statistics of the simulations, we computed the variance and skewness of the vertical velocity, which are shown in \cref{fig:dycoms_s_w}. The peak of the variance around the cloud base is an indication of turbulence generation due to latent heat release. Excessive dissipation will lead to suppressed turbulent fluctuations that is shown as a diminished peak of the velocity variance. Positive vertical velocity skewness near the bottom surface is a sign of convection from the surface driven by the heat flux. The negative vertical velocity skewness near the cloud top suggests the presence of a downdraft as a result of radiative cooling by the cloud. Spurious oscillations will lead to an overprediction of turbulent mixing, therefore overheating the cloud region.
\Cref{fig:dycoms_s_w} shows that the convection behavior is well represented in all simulations. Spurious mixing observed in the simulation with the QUICK scheme in \cref{fig:dycoms_thteta_li} is insignificant so that the desired convection behavior is preserved. While all the WENO-based schemes make comparable predictions of the vertical velocity skewness, both WENO3-JS and WENO5-JS show a lower peak for the velocity variance. Results with the WENO-NN scheme show a significant improvement in the agreement with experimental observations.

Based on these results, we observe that the WENO-NN scheme preserves the monotonicity of the numerical flux without introducing excessive numerical dissipation. This balanced dissipation behavior makes it useful for turbulent flow simulations.

\section{Conclusions}\label{sec:conclusion}

This work introduced a novel approach for dynamically adapting WENO weights using rational neural networks. We showcased its effectiveness and versatility across a diverse range of test cases, from simple advection to complex stratocumulus simulations exhibiting three-dimensional turbulence. We demonstrated that data-driven approaches, using networks with better representation properties together with robust model selection criteria, can result in numerical methods with small errors in the low-resolution regime, while regaining their asymptotic optimality as the resolution increases. Their improved accuracy in the low-resolution regime is crucial to obtain phenomenologically correct simulations of large scale system with complex physics, which results in simulations better aligned with experimental observations, highlighting its ability to maintain the right level of numerical smoothing and prevent excessive mixing.

Our approach combines a multilayer perceptron with trainable rational functions as activations to estimate WENO weights from cell-averaged field values. Using the approximation properties of rational functions, we adapt the stencil weights to the local smoothness of the solution. This strategy mitigates the high dissipation and pre-asymptotic errors typically seen in low-resolution applications with conventional WENO3 schemes. As the network is trained offline on a class of analytical functions, we bypass the need for a differentiable solver while enabling efficient and precise validation of the network outputs against analytical benchmarks. The network can be applied across different grid resolutions without any retraining.

One possible direction is to couple this low-resolution accuracy enhancing techniques with closure models, enabling  more accurate simulation of multi-physics problems. Another direction is to further reduce the computational cost of the application of the network by replacing some of the layers by analytical expressions stemming from analytical regression.

\section*{Acknowledgements}
We thank Tapio Schneider, Matthias Ihme and Jeff Parker for valuable discussions and helpful comments. Authors at Technical University of Munich and Georgia Institute of Technology gratefully acknowledge Google Cloud for providing compute credit.

\bibliographystyle{bibsty}
\bibliography{references.bib}




\end{document}